\documentclass[preprint,authoryear]{elsarticle}

\let\elscorref\corref
\usepackage{lineno}
\modulolinenumbers[5]
\journal{Annual Reviews in Control}

\usepackage{ifthen}
\usepackage{pgffor}
\usepackage{etoolbox}
\usepackage{parselines}
\usepackage{environ}
\usepackage{keyval}
\usepackage{listofitems}
\usepackage{forloop}
\usepackage{calc}
\usepackage{xstring}
\usepackage{changepage}
\usepackage{scrextend}
\usepackage{expl3}

\makeatletter
\ifthenelse{\isundefined{\lefteqs}}{
\usepackage{amsmath}
\usepackage{mathtools}
}{
\usepackage[fleqn]{amsmath}
\usepackage[fleqn]{mathtools}
}
\usepackage{amssymb}
\@ifclassloaded{autart}{}{\usepackage{amsthm}}
\makeatother
\usepackage{thmtools}
\usepackage{ar}
\usepackage{empheq}
\usepackage{mathrsfs}
\usepackage{bm}
\usepackage{booktabs}

\makeatletter
\let\ltxxlabel\ltx@label
\makeatother

\usepackage{color}
\makeatletter
\@ifclassloaded{beamer}{}{
\usepackage[table]{xcolor}}
\makeatother

\usepackage{graphicx}
\makeatletter
\@ifclassloaded{UWThesis/uwthesis}{
\usepackage{subfig}
}{\usepackage{subcaption}}
\makeatother
\usepackage{caption}

\usepackage{tabularx}
\usepackage{adjustbox}
\usepackage{multicol}
\usepackage{diagbox} % \diagbox{rows}{columns} with slash in between
\usepackage{varwidth} % V{max_width} variable-width cell type
\usepackage{makecell} % \makecell[l|r|c]{line 1 \\ line 2}
\usepackage{arydshln}

\usepackage{tikz}
\usepackage{tikz-3dplot}
\usetikzlibrary{backgrounds}
\usetikzlibrary{arrows}
\usetikzlibrary{arrows.meta}
\usetikzlibrary{intersections}
\usetikzlibrary{math}
\usetikzlibrary{shapes}
\usetikzlibrary{positioning}
\usetikzlibrary{decorations.markings}
\usetikzlibrary{decorations.pathreplacing}
\usetikzlibrary{decorations.text}
\usetikzlibrary{tikzmark}
\usetikzlibrary{fadings}

\makeatletter
\@ifclassloaded{standalone}{}{\usepackage[tikz]{mdframed}}
\makeatother
\usepackage{tocloft}
\usepackage{enumitem}
\usepackage{listings}
\usepackage{algorithm, algorithmicx}
\usepackage{algpseudocode}%[noend]
\usepackage[many]{tcolorbox}
\usepackage[binary-units=true]{siunitx}
\usepackage{chemformula}
\usepackage[normalem]{ulem}
\usepackage{cancel}
\usepackage{accents}
\makeatletter
\@ifclassloaded{beamer}{\usepackage{multimedia}}{}
\makeatother
\usepackage{lipsum}
\usepackage{mfirstuc}
\usepackage{import}
\usepackage{nicematrix}
\usepackage{bibentry} % use \bibentry{...} for in-text citation
\makeatletter
\ifthenelse{\isundefined{\natbiboptions}}{
\def\natbiboptions{numbers,sort}
}{}
\ifthenelse{\isundefined{\nobibpkg}}{
\@ifclassloaded{ieeeconf}{}{
\@ifclassloaded{elsarticle}{}{\usepackage[\natbiboptions]{natbib}\nobibliography*}}
}{}
\makeatother

\usepackage{float}
\usepackage{afterpage}
\ifthenelse{\isundefined{\landscapemode}}{
\def\landscapemode{pdflscape}
}{}
\usepackage{\landscapemode} % writing (rotated PDF pages, don't need to twist head)
\usepackage{varwidth}
\usepackage{nomencl}
\usepackage{longtable}
\usepackage{marginnote}
\makeatletter
\@ifclassloaded{beamer}{}{
\@ifclassloaded{autart}{}{
\@ifclassloaded{standalone}{}{
\@ifclassloaded{ieeeconf}{}{
\usepackage[explicit]{titlesec}
}}}}
\makeatother
\usepackage{soul}
\usepackage{stackengine}
\usepackage{scalerel}
\usepackage{lastpage}
\usepackage{fancyhdr}
\usepackage{nameref}

\usepackage[breaklinks=true]{hyperref}
\usepackage{xspace}

\let\corref\relax
\newcommand{\union}{\cup}
\newcommand{\intersection}{\cap}

\newcommand{\set}[1]{\mathcal #1}
\newcommand{\setcomplement}[1]{#1^{\mathrm{c}}}
\newcommand{\eucledian}[1][]{\ifthenelse{\equal{#1}{}}{\mathbb E}{\mathbb #1}}
\DeclareMathOperator{\@convhull}{conv}
\newcommand{\reals}[1][]{%
  \ifthenelse{\equal{#1}{extended}}{\overline}{}{\mathbb R}%
}
\newcommand{\real}{\mathbb R}

\newcommand{\binary}{\mathbb I}

\newcommand{\dom}{\mathrm{dom}}

\newcommand{\algvar}[1]{\text{\IfSubStr{#1}{_}{%
    \StrSubstitute{#1}{_}{\textunderscore}}{#1}}}

\definecolor{listinggray}{gray}{0.9}
\definecolor{lbcolor}{rgb}{0.9,0.9,0.9}
\colorlet{Darkgreen}{green!60!black}
\newcommand{\grad}{\nabla}

\newcommand{\subdiff}[1][none]{%
  \ifthenelse{\equal{#1}{none}}{%
    \partial%
  }{%
    \partial_{#1}%
  }%
}

\makeatletter
\DeclareMathOperator{\@convhull}{conv}
\newcommand{\convhull}[1][]{%
  \ifthenelse{\equal{#1}{closed}}{\overline}{}\@convhull%
}
\DeclareMathOperator{\@closed@convex@envelope}{\overline{co}}
\newcommand{\cce}{\@closed@convex@envelope}
\makeatother

\newcommand{\support}[2][delta]{\ifthenelse{\equal{#1}{delta}}{%
    \delta^*_{#2}%
  }{\sigma_{#2}}}

\newcommand{\diag}[2][short]{\mathrm{diag}\ifthenelse{\equal{#1}{short}}{(#2)}{\left(#2\right)}}

\newcommand{\pinv}{^\dagger}

\newcommand{\Matrix}[2][]{%
  \ifthenelse{\equal{#1}{}}{}{\setlength\arraycolsep{#1}}%
  \begin{bmatrix}#2\end{bmatrix}}
\renewcommand{\Array}[2][]{%
  \ifthenelse{\equal{#1}{}}{}{\setlength\arraycolsep{#1}}%
  \begin{matrix}#2\end{matrix}}
\makeatletter
\newcommand{\toset}{%
  \def\arr@offset{0.15em}
  \def\arr@len{0.7em}
  \def\arr@height{0.3em}
  \tikz[minimum height=0ex,outer sep=0,inner sep=0]
  \path[-{Latex[length=0.8mm]}]
  node (a) at (0,0) {}
  node (b) at (0,\arr@height) {}
  (a) edge ++(\arr@len,0)
  (b) edge ++(\arr@len,0)
  (a) edge[draw=none] ++(0,-\arr@offset);%
}
\makeatother

\newcommand{\inv}[1][]{\ifthenelse{\equal{#1}{}}{^{-1}}{^{-#1}}}

\newcommand{\norm}[2][]{\left\|#2\right\|_{#1}}

\newcommand{\vertiii}[1]{{\left\vert\kern-0.25ex\left\vert\kern-0.25ex\left\vert #1
        \right\vert\kern-0.25ex\right\vert\kern-0.25ex\right\vert}}

\makeatletter
\define@key{projKeys}{set}{\def\set{#1}}
\define@key{proxKeys}{t}{\def\t{#1}}
\define@key{proxKeys}{f}{\def\f{#1}}
\makeatother

\DeclareMathOperator{\proximal}{prox}

\newcommand{\Prox}[2][]{%
  \setkeys{proxKeys}{t=t,#1}%
  \setkeys{proxKeys}{f=f,#1}%
  \mathchoice{\underset{\t\f}{\proximal}}%
  {\proximal_{\t\f}}{\proximal_{\t\f}}{\proximal_{\t\f}}%
  \ifthenelse{\equal{#2}{}}{}{\left(#2\right)}%
}
\DeclareMathOperator{\distance}{dist}
\newcommand{\dist}[2][]{%
  \distance_{#1}\ifthenelse{\equal{#2}{}}{}{(#2)}
}
\newcommand{\gauge}[2][]{%
  \gamma_{#1}\ifthenelse{\equal{#2}{}}{}{(#2)}
}

\newcommand{\transp}{{\scriptscriptstyle\mathsf{T}}}
\newcommand{\T}{^\transp}

\newcommand{\dd}{\mathrm{d}}
\newcommand{\dt}{\dd t}

\newcommand{\fun}[2][1]{%
  #2(%
  \foreach \index in {1, ..., #1} {%
    \ifthenelse{\equal{\index}{#1}}{%
      \cdot%
    }{%
      \cdot,%
    }%
  })}
\newcommand{\definedas}[1][tri]{%
  \ifthenelse{\equal{#1}{tri}}{\triangleq}{\coloneqq}}

\renewcommand{\implies}{\Rightarrow}
\renewcommand{\iff}{\Leftrightarrow}

\newcommand{\mgeq}{\succeq}

\newcommand{\mgr}{\succ}
\DeclareMathOperator*{\exptx}{exp}
\renewcommand{\exp}[2][exponent]{\ifthenelse{\equal{#1}{exponent}}{e^{#2}}{\exptx\left(#2\right)}}
\DeclareMathOperator*{\expsf}{\mathsf{exp}}
\newcommand{\expm}[2][exponent]{\ifthenelse{\equal{#1}{exponent}}{\mathsf{e}^{#2}}{\expsf\left(#2\right)}}

\makeatletter
\DeclareFontFamily{U}{tipa}{}
\DeclareFontShape{U}{tipa}{m}{n}{<->tipa10}{}
\newcommand{\arc@char}{{\usefont{U}{tipa}{m}{n}\symbol{62}}}%
\renewcommand{\arc}[1]{\mathpalette\arc@arc{#1}}
\newcommand{\arc@arc}[2]{%
  \sbox0{$\m@th#1#2$}%
  \vbox{
    \hbox{\resizebox{\wd0}{\height}{\arc@char}}
    \nointerlineskip
    \box0
  }%
}
\makeatother



\newcommand{\ev}[3][c]{%
  \ifthenelse{\equal{#1}{c}}{%
    \ifthenelse{\equal{#2}{}}{\forall[0,#3]}{\forall[#2,#3]}%
  }{%
    \ifthenelse{\equal{#1}{o}}{%
      \ifthenelse{\equal{#2}{}}{\forall(0,#3)}{\forall(#2,#3)}%
    }{%
      \ifthenelse{\equal{#1}{oc}}{%
        \ifthenelse{\equal{#2}{}}{\forall(0,#3]}{\forall(#2,#3]}%
      }{%
        \ifthenelse{\equal{#2}{}}{\forall[0,#3)}{\forall[#2,#3)}%
      }%
    }%
  }%
}

\DeclareMathOperator*{\argmin}{argmin}

\newcommand{\Behcet}{Beh\c{c}et}
\newcommand{\Acikmese}{A\c{c}{\i}kme\c{s}e}

\makeatletter
\@ifclassloaded{beamer}{}{}
\makeatother

\newcommand{\scvx}{SCvx\xspace}

\newcommand{\defintext}[1]{\textbf{#1}}

\newcommand{\GetLabel}[1]{\expandafter\csname #1Label\endcsname}
\makeatletter
\define@key{printRefKeys}{otherLabel}{\def\otherLabel{#1}}
\define@key{printRefKeys}{concatenate}{\def\concatenate{#1}}
\makeatother
\newcommand{\PrintRefs}[3][]{%
  \setkeys{printRefKeys}{otherLabel=,#1}%
  \setkeys{printRefKeys}{concatenate=false,#1}%
  \xdef\MyModifiedLabel{#2}%
  \ifthenelse{\equal{\otherLabel}{}}{}{\xdef\MyModifiedLabel{\otherLabel}}%
  \xdef\MyRefCount{0}%
  \foreach \i in {#3} {%
    \tikzmath{\MyRefCount=int(\MyRefCount+1);}%
    \xdef\MyRefCount{\MyRefCount}%
  }%
  \ifthenelse{\equal{\MyRefCount}{1}}{%
    #2~\ref{\GetLabel{\MyModifiedLabel}:#3}%
  }{%
    \xdef\MyCounter{0}%
    #2s~%
    \foreach \AlgRef in {#3} {%
      \tikzmath{\MyCounter=int(\MyCounter+1);}%
      \xdef\MyCounter{\MyCounter}%
      \ifthenelse{\equal{\concatenate}{true}}{%
        \pgfmathparse{\MyCounter==1 ? 1 : 0}%
        \ifthenelse{\pgfmathresult>0}{%
          \ref{\GetLabel{\MyModifiedLabel}:\AlgRef}-%
        }{%
          \pgfmathparse{\MyCounter==\MyRefCount ? 1 : 0}%
          \ifthenelse{\pgfmathresult>0}{%
            \ref{\GetLabel{\MyModifiedLabel}:\AlgRef}%
          }{}%
        }%
      }{%
        \pgfmathparse{\MyCounter<\MyRefCount-1 ? 1 : 0}%
        \ifthenelse{\pgfmathresult>0}{%
          \ref{\GetLabel{\MyModifiedLabel}:\AlgRef},~%
        }{%
          \pgfmathparse{\MyCounter<\MyRefCount ? 1 : 0}%
          \ifthenelse{\pgfmathresult>0}{%
            \ref{\GetLabel{\MyModifiedLabel}:\AlgRef}~and~%
          }{%
            \ref{\GetLabel{\MyModifiedLabel}:\AlgRef}%
          }%
        }%
      }%
    }%
  }%
}
\newcommand{\sref}[1]{\PrintRefs{Section}{#1}}
\newcommand{\ssref}[2][]{\PrintRefs[otherLabel=Subsection,#1]{Section}{#2}}
\newcommand{\sssref}[2][]{\PrintRefs[otherLabel=Subsubsection,#1]{Section}{#2}}
\newcommand{\cref}[1]{\PrintRefs{Chapter}{#1}}

\newcommand{\figref}[1]{\PrintRefs{Figure}{#1}}
\newcommand{\tabref}[2][]{\PrintRefs[#1]{Table}{#2}}

\makeatletter
\define@key{algKeys}{start}{\def\startline{#1}}
\define@key{algKeys}{end}{\def\endline{#1}}
\define@key{algKeys}{show}{\def\showalg{#1}}
\makeatother
\renewcommand{\algref}[2][]{%
  \setkeys{algKeys}{start=,#1}%
  \setkeys{algKeys}{end=,#1}%
  \setkeys{algKeys}{show=true,#1}%
  \ifthenelse{\equal{\showalg}{true}}{\PrintRefs{Algorithm}{#2}}{}%
  \ifthenelse{\equal{\startline}{}}{}{~L\ref{alg:#2:line:\startline}%
    \ifthenelse{\equal{\endline}{}}{}{-\ref{alg:#2:line:\endline}}%
  }%
}
\newcommand{\pref}[2][]{\PrintRefs[#1]{Problem}{#2}}

\newcommand{\corref}[2][]{\PrintRefs[#1]{Corollary}{#2}}

\input{optimization.tex}

\usepackage{doi}

\newcommand{\Nx}{{n_x}} % Number of states
\newcommand{\Nu}{{n_u}} % Number of inputs
\newcommand{\Nz}{{n_z}} % Number of binary variables
\newcommand{\Nb}{{n_b}} % Number of boundary conditions
\newcommand{\Nc}{{n_c}} % Number of triggered constraints
\newcommand{\NN}{N} % Number of temporal nodes

\newcommand{\RNx}{\reals^\Nx}
\newcommand{\RNu}{\reals^\Nu}
\newcommand{\RNz}{\reals^\Nz}

\newcommand{\RNg}{\reals^\Nc}
\newcommand{\RNb}{\reals^\Nb}
\newcommand{\RNc}{\reals^\Nc}

\renewcommand{\defintext}[1]{\textit{#1}}

\newcommand{\ones}{\bm{1}}

\graphicspath{{figures/}}

\hypersetup{
  colorlinks=true,
}

\allowdisplaybreaks

\def\uwaddress{William E. Boeing Department of Aeronautics and Astronautics, University of
  Washington, Seattle, WA 98195, USA}

\begin{document}
\hypersetup{
  linkcolor=red,
  urlcolor=blue,
  citecolor={blue!60}
}

\begin{frontmatter}

  \title{Advances in Trajectory Optimization for Space Vehicle Control}

  \author[uwaddress]{Danylo Malyuta\elscorref{correspondingauthor}}
  \cortext[correspondingauthor]{Corresponding author}
  \ead{danylo@uw.edu}
  \author[uwaddress]{Yue Yu}
  \ead{yueyu@uw.edu}
  \author[uwaddress]{Purnanand Elango}
  \ead{pelango@uw.edu}
  \author[uwaddress]{\Behcet{} \Acikmese}
  \ead{behcet@uw.edu}
  \address[uwaddress]{\uwaddress}

  \begin{abstract}
    Space mission design places a premium on cost and operational
    efficiency. The search for new science and life beyond Earth calls for
    spacecraft that can deliver scientific payloads to geologically rich yet
    hazardous landing sites. At the same time, the last four decades of
    optimization research have put a suite of powerful optimization tools at
    the fingertips of the controls engineer. As we enter the new decade,
    optimization theory, algorithms, and software tooling have reached a
    critical mass to start seeing serious application in space vehicle guidance
    and control systems. This survey paper provides a detailed overview of
    recent advances, successes, and promising directions for optimization-based
    space vehicle control. The considered applications include planetary
    landing, rendezvous and proximity operations, small body landing,
    constrained attitude reorientation, endo-atmospheric flight including
    ascent and reentry, and orbit transfer and injection. The primary focus is
    on the last ten years of progress, which have seen a veritable rise in the
    number of applications using three core technologies: lossless
    convexification, sequential convex programming, and model predictive
    control. The reader will come away with a well-rounded understanding of the
    state-of-the-art in each space vehicle control {application, and will be}
    well positioned to tackle important current open problems using convex
    optimization as a core technology.%
  \end{abstract}


\end{frontmatter}

\newpage
\section*{Contents}
\makeatletter
\renewcommand\tableofcontents{%
  \@starttoc{toc}%
}
\makeatother
\tableofcontents

%%%%%%%%%%%%%%%%%%%%%%%%%%%%%%%%%%%%%%%%%%%%%%%%%%%%
\section{Introduction}
\label{section:introduction}
%%%%%%%%%%%%%%%%%%%%%%%%%%%%%%%%%%%%%%%%%%%%%%%%%%%%

Improvements in computing hardware and maturing software libraries have made
optimization technology become practical for space vehicle control. The term
computational guidance and control (CGC) was recently coined to refer to
control techniques that are iterative in nature and that rely on the onboard
computation of control actions \citep{lu2017cgc,tsiotras2017toward}.

This paper surveys \textit{optimization-based} methods, which are a subset of
CGC for space vehicles. We consider applications for launchers, planetary
landers, satellites, and spacecraft. The common theme across all applications
is the use of an optimization problem to achieve a control objective. Generally
speaking, the goal is to solve:
\begin{optimus}[
  task=\min,
  variables={\bm{x}},
  objective={J(\bm{x})},
  plabel={intro_opt}]
  \bm{x}\in\set C,
\end{optimus}
where $J:\reals^n\to\reals$ is a \defintext{cost function},
$\set{C}\subseteq\reals^n$ is a \defintext{feasible set}, and $\bm{x}$ is an
$n$-dimensional vector of \defintext{decision variables}. Optimization is a
relevant area of study for modern space vehicle control for two reasons:
effectiveness of formulation, and the (emerging) existence of efficient
solution methods.
% {\color{red} Taylor: I would add that increasingly complex \& demanding
% missions (like TRN \& vision based nav) are starting to require satisfaction
% of challenging constraints. Optimization-based methods are a good way of
% finding feasible solutions that might elude classical methods.}

To answer why an optimization formulation is effective, consider the physical and operational constraints on the tasks that
recent and future space vehicles aim to perform. Future launchers and planetary
landers will require advanced entry, descent, and landing (EDL) algorithms to
drive down cost via reusability, or to access scientifically interesting sites
\citep{blackmore2016autonomous}. Instead of landing in open terrain, future
landers will navigate challenging environments such as volcanic vents and jagged
blades of ice
\citep{sanmartin2013development,robertson2017synopsis,europa2012study}. Meanwhile,
human exploration missions will likely be preceded by cargo delivery, requiring
landings to occur in close proximity
\citep{dwyercianciolo2019defining}. Motivated by the presence of water ice,
upcoming missions to the Moon will target its south pole
\citep{artemis2019moon}, where extreme light-dark lighting conditions call for
an automated sensor-based landing \citep{lroc2018lunar}. Even for robotic
missions, new onboard technology such as vision-based terrain relative
navigation requires the satisfaction of challenging constraints that couple
motion and sensing. Regardless of whether one achieves the lowest cost
\eqref{eq:intro_opt_a} or not, optimization is indeed one of the most compelling
frameworks for finding feasible solutions in the presence of challenging
constraints \citep{tsiotras2017toward}.

In orbit, foreseeable space missions will necessitate robotic docking for sample
return, debris capture, and human load alleviation
\citep{woffinden2007navigating}. Early forms of the capability have been shown
on the Japanese ETS-VII, European ATV, Russian Soyuz, US XSS-11, and US DART
demonstrators. Most recently, the human-rated SpaceX Crew Dragon performed
autonomous docking with the ISS, and the Orion Spacecraft is set to also feature
this ability \citep{stephens2013orion,dsouza2007orion}. Further development in
autonomous spacecraft rendezvous calls for smaller and cheaper sensors as well
as a reduction in the degree of cooperation by the target spacecraft. This will
require more flexibility in the chaser's autonomy, which is practically
achievable using onboard optimization.

% A planetary lander will have to satisfy a multitude of navigation sensor
% pointing constraints as it identifies hazardous geography while seeking to land
% on a dwindling fuel supply \citep{carson2019splice}.
The above mission objectives suggest that future space vehicle autonomy will
have to adeptly operate within a multitude of operational constraints. % To design
% a successful control strategy for such systems means to be \defintext{feasible}
% with respect to these constraints.
However, optimality usually stipulates operation near the boundary of the
set of feasible solutions. {In other words, the vehicle must activate
  its constraints (i.e., touch the constraint set boundary) at important or
  prolonged periods of its motion}. By virtue of the feasible set $\set C$ in
\eqref{eq:intro_opt_b}, optimization is one of the few suitable methods (and is
perhaps the most natural one) to directly impose system constraints
\citep{mayne2000constrained}.

\begin{figure}
  \centering
  \includegraphics[width=\textwidth]{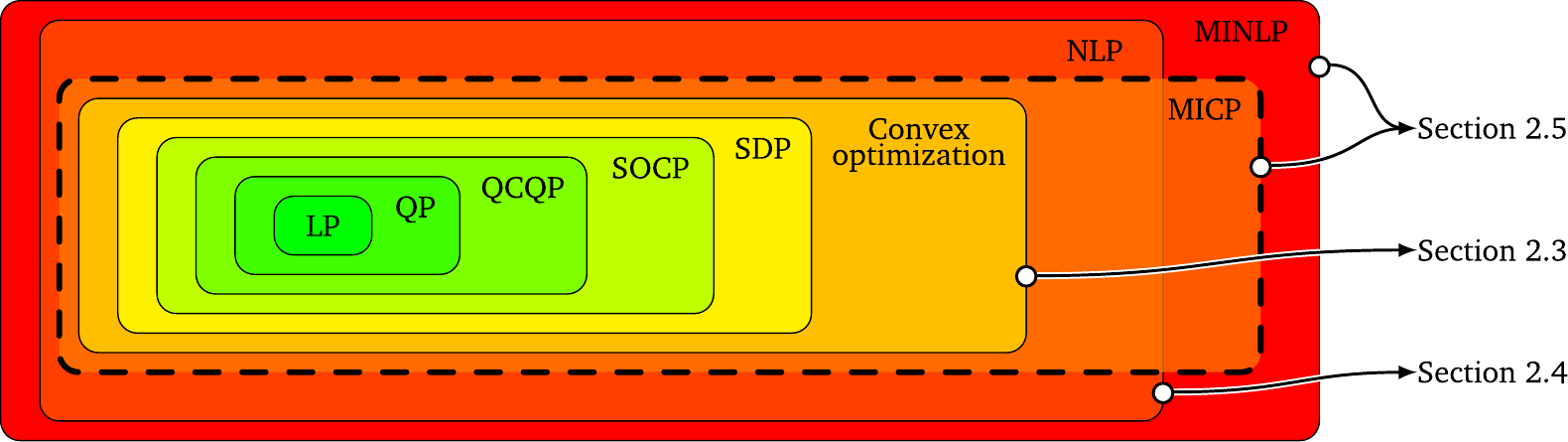}
  \caption{Taxonomy of optimization problems. Going from inside to outside, each
    class becomes more difficult to solve. Roughly speaking, SOCP is currently
    the most general class that can be solved reliably enough to be deployed on
    a space system in short order.}
  \label{fig:taxonomy}
\end{figure}

The benefit of an appropriate formulation, however, is limited if no algorithm
exists to solve \pref{intro_opt} \textit{efficiently}, which means quickly and
utilizing few computational resources. Convex optimization has been a popular
approach for formulating problems since it enables efficient solution
methods. \figref{taxonomy} illustrates a taxonomy of optimization problem
\textit{classes} or \textit{families}, of which convex optimization is a
part. The inner-most class in \figref{taxonomy} is the linear program (LP).  Next
comes the quadratic program (QP), followed by the second-order cone program
(SOCP). Troves of detail on each class may be found in many excellent
optimization textbooks
\citep{RockafellarConvexBook,NocedalBook,BoydConvexBook}. Roughly speaking, SOCP
is the most general class of problems that state-of-the-art algorithms can
solve with high reliability and rigorous performance guarantees
\citep{dueri2014automated,dueri2017customized,domahidi2013ecos}. Beyond SOCP,
the semidefinite program (SDP) class enables optimization over the space of
positive semidefinite matrices, which leads to many important robust control
design algorithms \citep{MIMOBook,LMIBook}. SDP is the most general class of
convex optimization for which off-the-shelf solvers are available, {and many
  advances have been made in recent years towards more scalable and robust SDP
  solvers \citep{majumdar2020recent}}.
% More general convex optimization problems than SDP exist, but they are
% lacking efficient numerical solver implementations.

Although convex optimization can solve a large number of practical engineering
problems, future space system requirements often surpass the flexibility of
``vanilla'' convex optimization. Solving nonconvex optimization problems will
be required for many foreseeable space vehicles \citep{carson2019splice}. Thus,
extending beyond SDP, we introduce three nonconvex problem classes.

First, one can abandon the convexity requirement, but retain function
continuity, leading to the nonlinear program (NLP) class. Here the objective and
constraint functions are continuous, albeit nonconvex. Alternatively, one could
retain convexity but abandon continuity. This leads to the mixed-integer convex
program (MICP) class, where binary variables are introduced to emulate discrete
\textit{switches}, such as those of valves, relays, or pulsing space thrusters
\citep{dueri2017customized,achterberg2007constrained,achterberg2013mixed}. Note
that the MICP and NLP classes overlap, since some constraints admit both forms
of expression. For example, the mixed-integer constraint:
\begin{equation}
  \label{eq:example_mip}
  x\in\{0\}\union\{x\in\reals:x\ge 1\},
\end{equation}
can be equivalently formulated as a nonlinear continuous constraint:
\begin{equation}
  \label{eq:example_nlp}
  x(x-1)\ge 0,~x\ge 0.
\end{equation}

In the most general case, nonlinearity and discontinuity are combined to form
the mixed-integer nonlinear program (MINLP) class. Since integer variables are
nowhere continuous and the corresponding solution methods are of a quite
different breed to continuous nonlinear programming, we reserve MINLP as the
largest and toughest problem class. Algorithms for NLP, MICP, and MINLP
typically suffer either from exponential complexity, a lack of convergence
guarantees, or both \citep{malyuta2019approximate}. Nevertheless, the
optimization community has had many successes in finding practical solution
methods even for these most challenging problems
\citep{achterberg2013mixed,szmuk2019successive}.
% Certainly more work is required, but each new conference edition and journal
% issue have brought the state-of-the-art closer to efficient methods for many
% of the constraints and objectives imposed on future space vehicles.

This paper stands in good company of numerous surveys on aerospace
optimization. \citep{betts1998survey} presents an eloquent, albeit somewhat
dated, treatise on trajectory optimization
methods. {\citep{trelat2012optimal} provides a comprehensive survey of
  modern optimal control theory and indirect methods for aerospace problems,
  covering geometric optimal control, homotopy methods, and favorable properies
  of orbital mechanics that can be leveraged for trajectory
  optimization. \citep{tsiotras2017toward} corroborate the importance of
  optimization in forthcoming space missions.}  \citep{liu2017survey} survey the
various appearances of lossless convexification and sequential convex
programming in aerospace guidance methods. \citep{eren2017model} cover
extensively the topic of model predictive control for aerospace applications,
where \pref{intro_opt} is solved recursively to compute control
actions. \citep{mao2018survey} survey three particular topics: lossless
convexification, sequential convex programming, and solver customization for
real-time computation. \citep{shirazi2018survey} provide a thorough discussion
on the general philosophy and specific methods and solutions for in-space
trajectory optimization.  Recently, \citep{song2020survey} surveyed optimization
methods in rocket powered descent guidance with a focus on feasibility, dynamic
accuracy, and real-time performance.

This paper contributes the most recent broad survey of convex
optimization-based space vehicle control methods. We consider rockets for
payload launch, rocket-powered planetary and small body landers, satellites,
{interplanetary spacecraft}, and atmospheric entry vehicles. {However, we do
  not cover some related topics like guidance of purely atmospheric vehicles
  (e.g., missiles and hypersonic aircraft), and control of satellite swarms, due
  to sufficiently unique distinctions. For a start in these areas, we refer the
  reader to
  \citep{palumbo2010modern,murillo2010fast,zarchan2019tactical,TewariBook} for
  hypersonic vehicle control, and
  \citep{rahmani2019swarm,morgan2012swarm,scharf2003survey,tillerson2002coordinate}
  for swarm control.}

From the algorithmic perspective, our focus is on convex optimization-based
methods for solving the full spectrum of optimization classes in
\figref{taxonomy}. The motivation for focusing on convex methods comes from the
great leaps in the reliability of convex solvers and the availability of flight
heritage, which gives convex optimization {a technology infusion
  advantage} for future onboard and ground-based algorithms
\citep{dueri2017customized,blackmore2016autonomous}. We nevertheless make side
references to other important, but not convex optimization-based, algorithms
throughout the text. {Lastly, this paper discusses algorithms at a high
  level, and chooses to cover a large number of applications and methods in
  favor of providing deep technical detail for each algorithm. The goal, in the end, is
  to expose the reader to dominant recent and future directions in convex
  optimization-based space vehicle control research.}

The paper is organized as follows. \sref{methods} covers general theory of
important optimization methods used throughout spaceflight
applications. \sref{applications} then surveys each space vehicle control
application individually. \ssref{pdg} surveys powered descent guidance for
planetary rocket landing. \ssref{rendezvous} discusses spacecraft rendezvous and
proximity operations, followed by a discussion in \ssref{smallbody} of its close cousin,
small body landing. Constrained attitude reorientation is covered in
\ssref{reorientation}. \ssref{endo} surveys endo-atmospheric flight, including
ascent and entry. Last but not least, orbit insertion and transfer are surveyed
in \ssref{orbit}. {We conclude the paper with a perspective on what lies
  ahead for computational guidance and control. As such, \sref{outlook}
  highlights some recent applications of machine learning to select
  problems. This final section also tabulates some of the optimization software
  tooling now available for getting started in optimization methods for
  spaceflight applications.}

\textit{Notation.} % Some esoteric notation definitions are in order.
Binary numbers belong to the set $\binary\definedas\{0,1\}$. Vectors are
written in bold, such as $\bm{x}\in\reals^n$ versus $y\in\reals$. The identity
matrix is generally written as $I$, and sometimes as $I_n\in\reals^{n\times n}$
in order to be explicit about size. The zero scalar, vector, or matrix is
always written as $0$, with its size derived from context. The vector of ones
is written as $\bm{1}$, with size again derived from context. Starred
quantities denote optimal values, for example $x^*$ is the optimal value of
$x$. We use $(\bm{a};\bm{b};\bm{c})$ to concatenate elements into a column
vector, like in MATLAB. The symbol $\otimes$ denotes the Kronecker matrix
product or quaternion multiplication, depending on context. The positive-part
function $[\bm{x}]^+\definedas\max\{0,\bm{x}\}$ saturates negative elements of
$\bm{x}$ to zero. Given a function $f(x(t),y(t),t)$, we simplify the argument
list via the shorthand $f[t]$. Throughout the paper, we interchangeably use the
terms ``optimization'' and ``programming'', courtesy of linear optimization
historically being used for planning military operations
\citep{WrightTalk}. When we talk about ``nonlinear programming'', we mean more
precisely ``nonconvex programming''. Convexity is now known to be the true
separator of efficient algorithms, however this discovery came after linear
programming already established itself as the dominant class that can be
efficiently solved via the Simplex method
\citep{rockafellar1993watershed}. Finally, ``guidance'' means ``trajectory
generation'', while ``navigation'' means ``state estimation''.

%%%%%%%%%%%%%%%%%%%%%%%%%%%%%%%%%%%%%%%%%%%%%%%%%%%%
\section{Background on Optimization Methods}
\label{section:methods}
%%%%%%%%%%%%%%%%%%%%%%%%%%%%%%%%%%%%%%%%%%%%%%%%%%%%

This section provides a broad overview of key algorithms for space vehicle
trajectory optimization. The main focus is on methods that exploit convexity,
since convex optimization is where state-of-the-art solvers provide the
strongest convergence guarantees at the smallest computational cost
\citep{NocedalBook,BoydConvexBook}.

Our algorithm overview proceeds as follows. First, \ssref{ocp} introduces the
general continuous-time optimal control problem. Then, \ssref{discretization}
describes how the problem is discretized to yield a finite-dimensional problem
that can be solved on a computer. Following this introduction,
\ssref[concatenate=true]{lcvx,mpc} overview important algorithms for space
vehicle trajectory optimization.

\subsection{Optimal Control Theory}
\label{subsection:ocp}

Optimal control theory is the bedrock of every trajectory optimization problem
\citep{PontryaginBook,berkovitz1974optimal}. The goal is to find an optimal
input trajectory for the following optimal control problem (OCP):
\begin{optimus}[
  task=\min,
  variables={t_f,\bm{u}},
  objective={L_f(\bm{x}(t_f),t_f)+\int_0^{t_f}L(\bm{x}(\tau),\bm{u}(\tau),\tau)\dd\tau},
  plabel={ocp}
  ]
  \dot{\bm{x}}(t) = \bm{f}(\bm{x}(t),\bm{u}(t),t),~\forall t\in [0,t_f], \#
  \bm{g}(\bm{x}(t),\bm{u}(t),t)\le 0,~\forall t\in [0,t_f], \#
  \bm{b}(\bm{x}(0),\bm{x}(t_f), t_f) = 0.
\end{optimus}

In \pref{ocp}, $\bm{x}:[0,t_f]\to\RNx$ is the state trajectory and
$\bm{u}:[0,t_f]\to\RNu$ is the input trajectory, while $t_f\in\reals$ is
the final time (i.e., the trajectory duration). The state evolves according to
the dynamics $\bm{f}:\RNx\times\RNu\times\reals\to\RNx$, and satisfies at all
times a set of constraints defined by
$\bm{g}:\RNx\times\RNu\times\reals\to\RNg$. At the temporal boundaries, the
state satisfies conditions provided by a boundary constraint
$\bm{b}:\RNx\times\RNx\times\reals\to\RNb$. The quality of an input trajectory
is measured by a cost function consisting of a running cost
$L:\RNx\times\RNu\times\reals\to\reals$ and a terminal cost
$L_f:\RNx\times\reals\to\reals$.

Two aspects differentiate \pref{ocp} from a typical parameter optimization
problem. First, the constraints encode a physical process governed by ordinary
differential equations (ODEs) \eqref{eq:ocp_b}.
%Second, time is a continuous variable.
Second, due to the continuity of time, the input trajectory has an infinite
{number of design parameters}. This makes \pref{ocp} a semi-infinite
optimization problem that cannot be directly implemented on a
computer. {In the following subsections, we provide a brief overview of
  two approaches for solving this problem, called the direct and indirect
  methods. Roughly speaking, direct methods discretize \pref{ocp} and solve it
  as a parameter optimization problem, while indirect methods attempt to
  satisfy the necessary conditions of optimality.}

\subsubsection{{Indirect Methods}}
\label{subsubsection:indirect}

The maximum principle, developed since the 1960s, extends the classical
calculus of variations and provides a set of necessary conditions of optimality
for \pref{ocp} \citep{PontryaginBook,hartl1995survey}. The maximum principle
has found numerous aerospace applications \citep{LonguskiBook}.

{The \defintext{indirect} family of optimization methods solves the
  necessary conditions of optimality, which involves a two-point boundary value
  problem (TPBVP) corresponding to the state and costate dynamics and their
  boundary conditions. Traditionally, this is solved by a single- or
  multiple-shooting method.}
%This requires one to pre-specify the time intervals over which the
%constraint \eqref{eq:ocp_c} is active \citep{Betts1998}.
One limitation of these methods is the requirement to specify in advance the time
intervals over which constraint \eqref{eq:ocp_c} is active
\citep{betts1998survey}. Other issues that hinder onboard implementation
include poor convergence stemming from a sensitivity to the initial guess, and
long computation time.

Despite these challenges, the indirect approach is often the only practical
solution method when aspects like numerical sensitivity and trajectory duration
rule out direct methods. Low-thrust trajectory optimization, discussed in
\ssref{orbit}, is a frequent candidate for the indirect approach since the low
thrust-to-weight ratios and long trajectory durations (from weeks to years)
create extreme numerical challenges when formulated as a parameter
optimization problem.

 Most indirect methods in aerospace literature solve only the
\textit{necessary} conditions of optimality for \pref{ocp}. However, nonlinear
optimization problems can have stationary points that are not local minima,
such as saddle points and local maxima. This has prompted interest in using
second-order conditions of optimality to ensure that the solution is indeed a
local minimum \citep{CesariBook}. At the turn of the century, researchers showed
how second-order information can be incorporated in orbit transfer applications
\citep{jo2000procedure}. In the last decade, further work used second-order
optimality conditions for orbit transfer and constrained attitude
reorientation problems \citep{caillau2012minimum,picot2012shooting}.

A promising modern indirect method family relies on homotopy in order to solve
the TPBVP
{\citep{pan2016double,Pan2019,Pan2020,Taheri2019,trelat2012optimal,rasotto2015multi,DiLizia2014high}}. {Homotopy
  aims to address the aforementioned challenges of slow convergence, initial
  guess quality, and active constraint specification.} The core idea is to
describe the problem as a family of problems parametrized by a \textit{homotopy
  parameter} $\kappa\in [0,1]$, such that the original problem is recovered for
$\kappa=1$, and the problem for $\kappa=0$ is trivially solved. For example,
consider solving a non-trivial root-finding problem:
\begin{equation}
  \label{eq:root_finding}
  \bm{F}(\bm{y})=0,
\end{equation}
where $\bm{y}\in\reals^n$ and $\bm{F}:\reals^n\to\reals^n$ is a smooth
mapping. A (linear) homotopy method will have us define the following homotopy
function:
\begin{equation}
  \label{eq:linear_homotopy}
  \bm{\Gamma}(\bm{y},\kappa)\definedas
  \kappa\bm{F}(\bm{y})+(1-\kappa)\bm{G}(\bm{y})=0,
\end{equation}
where $\bm{G}:\reals^n\to\reals^n$ is a smooth function that has a known or
easily computable root $\bm{y}_0\in\reals^n$. Popular choices are
$\bm{G}(\bm{y})=\bm{F}(\bm{y})-\bm{F}(\bm{y}_0)$, called \textit{Newton
  homotopy}, and $\bm{G}(\bm{y})=\bm{y}-\bm{y}_0$, called \textit{fixed-point
  homotopy}. In nonlinear homotopy, the function $\bm{\Gamma}(\bm{y},\kappa)$
is a nonlinear function of $\kappa$, but otherwise similar relationships
continue to hold.

% pan2018new              | probabillty 1 homotopy/continuous path existence
% watson2002probability   | probabillty 1 homotopy/continuous path existence

 The locus of points $(\bm{y},\kappa)$ where \eqref{eq:linear_homotopy}
holds is called a \textit{zero curve} of the root-finding problem. Success of
the homotopy approach relies on the zero curve being continuous in $\kappa$ on
the interval $\kappa\in [0,1]$, albeit possibly discontinuous in
$\bm{y}$. Unfortunately, the existence of such a curve is not guaranteed except
for a few restricted problems \citep{watson2002probability, pan2018new}. In
general, the loci of points satisfying \eqref{eq:linear_homotopy} may include
bifurcations, escapes to infinity, and limit points. Furthermore, the solution
at $\kappa=1$ may not be unique.

Nevertheless, homotopy methods have been developed to successfully traverse the
$\kappa\in [0,1]$ interval when a zero curve does exist. The essence of the
homotopy approach is to judiciously update an intermediate solution
$(\bm{y}_k,\kappa_k)$ so as to follow a $\kappa$-continuous zero curve from
$\bm{y}_0$ to $\bm{y}_K$, where $\bm{\Gamma}(\bm{y}_K,1)=\bm{F}(\bm{y}_K)=0$
and $K$ is the final iteration counter. At each iteration, some methods use a
Newton-based root finding approach \citep{pan2016double}, while others rely
solely on numerical integration \citep{caillau2012differential}. For further
details on the homotopy approach, we refer the reader to
\citep{pan2016double}.

% WRONG: Thanks to the continuity of $\bm{\Gamma}$ in both of its arguments,
% the loci of its roots trace continuous paths in the $(\bm{y},\kappa)$ space.

\subsubsection{{Direct Methods}}
\label{subsubsection:direct}

\defintext{Direct} methods offer a compelling alternative where one discretizes
\pref{ocp} and solves it {as a parameter optimization problem} via
numerical optimization. The resulting solution in the convex case is usually
very close to the optimal continuous-time one. {As discussed in the
  next section, the solution can satisfy \eqref{eq:ocp_b} exactly if an exact
  discretization method is used \citep{szmuk2018successive}}. The optimization
step is most often performed by a primal-dual interior point method (IPM), for
which a considerable software ecosystem now exists thanks to 40 years of active
development
\citep{NocedalBook,forsgren2002interior,wright2005interior}. {Some of
  this software is listed in \ssref{optimization_software}.}

{Thanks to this expansive software ecosystem, and the large research
  community actively working on numerical optimization algorithms, direct
  methods may be considered as the most popular approach today. Their ability
  to ``effortlessly'' handle constraints like \eqref{eq:ocp_c} makes them
  particularly attractive \citep{betts1998survey}. In the remainder of this
  paper, our main focus is on direct methods that use convex optimization.}

{Nevertheless, as mentioned in the previous section, indirect methods
  are still relevant for problems that exhibit peculiarities such as extreme
  numerical sensitivity. It must further be emphasized that some of the best
  modern algorithms have resulted from the \textit{combined} use of an indirect
  and a direct approach. Typically an indirect method, and in particular the
  necessary conditions of optimality, can be used to discover the solution
  structure, which informs more efficient customized algorithm design for a
  direct solution method. We will discuss how this ``fusion'' approach was
  taken for powered descent guidance and atmospheric entry applications in
  \ssref{lcvx} and \sssref{reentry} respectively. Last but not least, the maximum
  principle can sometimes be used to find the analytic globally optimal
  solution for problems where no direct method can do so (e.g., nonlinear
  problems). In this case, an indirect method can provide a reference solution
  against which one can benchmark a direct method's performance
  \citep{reynolds2020real,sundstrom2009generic}}. In summary, indirect and
direct methods play complementary roles: the former is a good ground-truth and
analysis tool, while the latter is preferred for real-time onboard
implementation.

\subsection{Discretization}
\label{subsection:discretization}

\def\subfigwidth{0.3\textwidth}
\def\sppattfigsep{4mm}
\begin{figure*}
  \centering
  \makebox[\textwidth]{\makebox[1.3\textwidth]{
      \begin{subfigure}[t]{\subfigwidth}
        \centering
        \includegraphics{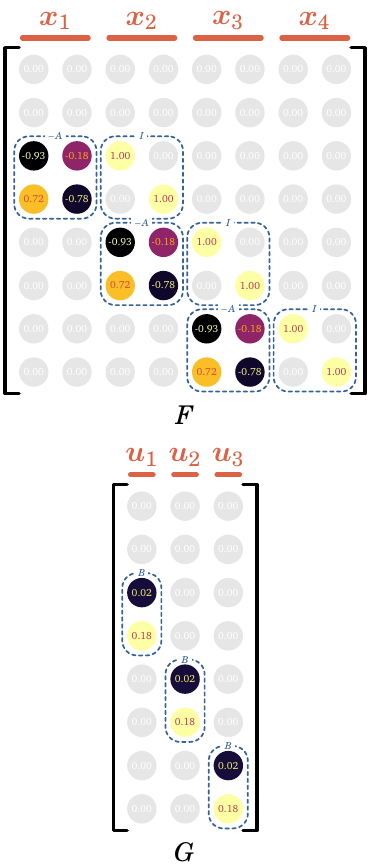}
        \caption{ZOH \eqref{eq:zoh_stacked_explicit}.}
        \label{fig:spy_zoh}
      \end{subfigure}%
      \hspace{\sppattfigsep}%
      \begin{subfigure}[t]{\subfigwidth}
        \centering
        \includegraphics{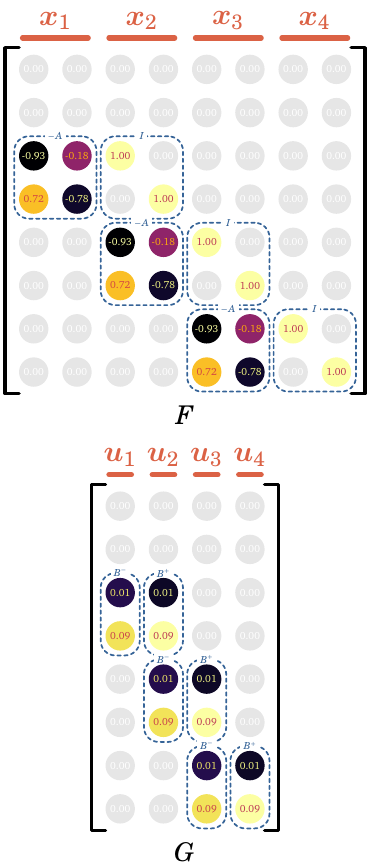}
        \caption{RK4 \eqref{eq:rk4_stacked_explicit}.}
        \label{fig:spy_rk4}
      \end{subfigure}%
      \hspace{\sppattfigsep}%
      \begin{subfigure}[t]{\subfigwidth}
        \centering
        \includegraphics{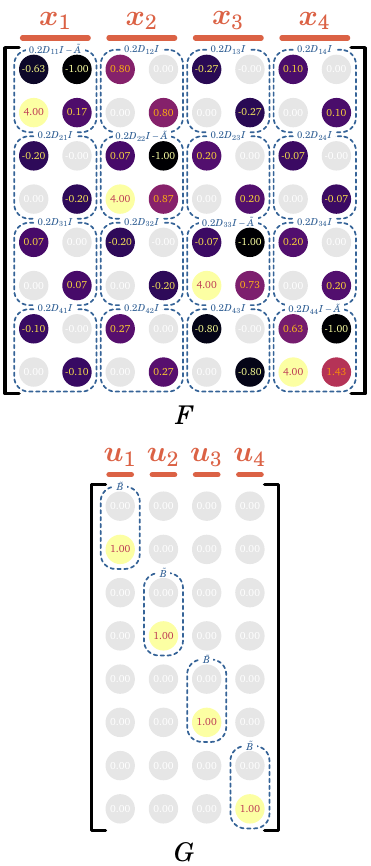}
        \caption{Global CGL \eqref{eq:global_pseudo_FG}.}
        \label{fig:spy_cgl}
      \end{subfigure}%
      \hspace{\sppattfigsep}%
      \begin{subfigure}[t]{0.24\textwidth}
        \centering
        \includegraphics{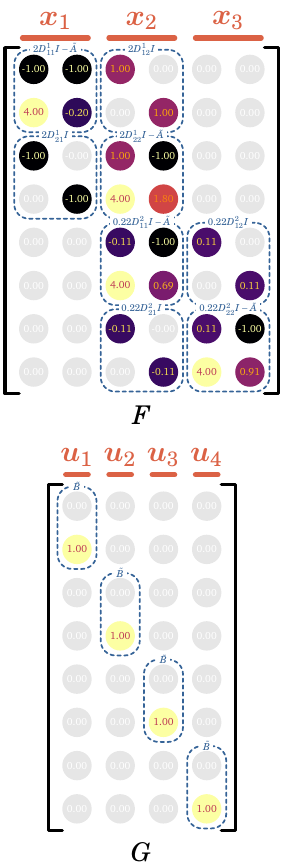}
        \caption{Adaptive CGL \eqref{eq:local_pseudo_FGl}.}
        \label{fig:spy_cgl_seg}
      \end{subfigure}%
    }}

  \caption{{Matrix sparsity patterns for the linear dynamics equation
      $F\bm{X}=G\bm{U}$ for the mass-spring-damper system in
      \sssref{toy_system_discretization}, using each of the discretization methods
      discussed in \ssref{discretization}. A salient feature of the ZOH and RK4
      methods is their relative sparsity compared to pseudospectral
      methods. The adaptive collocation in (\protect\subref{fig:spy_cgl_seg})
      increases sparsity by subdividing time into multiple intervals within
      which separate collocations are applied, and which are coupled only
      through continuity at their interface (in this figure, at
      $\bm{x}_2$). The non-zero element colors merely serve to visually
      separate the elements by their value (larger values correspond to warmer
      colors).}}
  \label{fig:spy}
\end{figure*}

 To be solvable on a computer, the semi-infinite \pref{ocp} must
generally be reduced to a finite-dimensional problem. This is done by the
process of \textit{discretization}, where the goal is to convert the
differential constraint \eqref{eq:ocp_b} into a finite-dimensional algebraic
constraint. This is especially important for the family of direct methods discussed
in \sssref{direct}, which rely on discretization to solve \pref{ocp} as a
parameter optimization problem.

Generally, discretization is achieved by partitioning time into a grid of $\NN$
nodes and fixing a basis for the state signal, the control signal, or both
\citep{malyuta2019discretization}. The following subsections discuss three
popular approaches: an exact discretization based on zeroth-order hold
(\sssref{zoh}), an approximate discretization based on the classic Runge-Kutta
method (\sssref{rk4}), and a pseudospectral discretization that is either global or
adaptive (\sssref{pseudospectral}). We frame the discussion in terms of three
salient features: 1) sparsity of the discrete representation of the dynamics
\eqref{eq:ocp_b}, 2) the mechanics of obtaining a continuous-time trajectory from
the discrete representation, and 3) the connection, if any, between the
discrete solution and the optimal costates of the original \pref{ocp} derived via
the maximum principle. Our goal is to give the reader enough insight into
discretization to appreciate the algorithmic choices for spaceflight
applications in \sref{applications}. For a more thorough discussion, we defer to
the specialized papers
\citep{betts1998survey,betts2010practical,kelly2017introduction,conway2011survey,
  ross2012review,rao2010survey,agamawi2020comparison,malyuta2019discretization,Phogat2018discreteAutomatica}.

\subsubsection{{Example Dynamical System}}
\label{subsubsection:toy_system_discretization}

 To ground our coverage of discretization in a concrete
application, let us restrict \eqref{eq:ocp_b} to a linear time-invariant (LTI)
system of the form:
\begin{equation}
  \label{eq:lti}
  \dot{\bm{x}}(t) = \tilde A\bm{x}(t)+\tilde B\bm{u}(t).
\end{equation}

Discretization for \eqref{eq:lti} is easily generalized to handle linearized
dynamics of a nonlinear system (like \eqref{eq:ocp_b}) about a reference state-input trajectory
$(\bar{\bm{x}},\bar{\bm{u}}):\reals\to\RNx\times\RNu$. To do so,
replace $\tilde A$ and $\tilde B$ with:
\begin{subequations}
  \label{eq:f_jacobians}
  \begin{align}
    \tilde A(t) &= \grad_{\bm{x}} \bar{\bm{f}}[t], \\
    \tilde B(t) &= \grad_{\bm{u}} \bar{\bm{f}}[t], \\
    \bar{\bm{f}}[t] &= \bm{f}(\bar{\bm{x}}(t),\bar{\bm{u}}(t),t),
  \end{align}
\end{subequations}
and add a residual term
$\bm{r}(t)=\bar{\bm{f}}[t]-\tilde A(t)\bar{\bm{x}}(t)-\tilde B(t)\bar{\bm{u}}(t)$
to the right-hand side of \eqref{eq:lti}. Note that in this case, \eqref{eq:lti}
generally becomes a linear time-varying (LTV) system. While the zeroth-order
hold method as presented in \sssref{zoh} requires linear dynamics, the
Runge-Kutta and pseudospectral methods in their full generality can in fact
handle the unmodified nonlinear dynamics \eqref{eq:ocp_b}.

As a particular example of \eqref{eq:lti}, we will consider a simple
mass-spring-damper system with $\Nx=2$ states and $\Nu=1$ input:
\begin{equation}
  \label{eq:msd_2ndOrder}
  \ddot r(t)+2\zeta\omega_n\dot r(t)+\omega_n^2 r(t)=m\inv f(t),
\end{equation}
where $r$ is the position, $m$ is the mass, $\zeta$ is the damping ratio,
$\omega_n$ is the natural frequency, and $f$ is the force (input). We set
$m=1~\si{\kilo\gram}$, $\zeta=0.2$, and
$\omega_n=2~\si{\radian\per\second}$. Furthermore, consider a staircase input
signal where $f(t)=1$ for $t\in [0,t_{step})$ and $f(t)=0$ for $t\ge t_{step}$.
We shall use $t_{step}=1~\si{\second}$. The initial condition is
$r(0)=\dot r(0)=0$. The simulation source code for this example is publically
available\footnote{Visit \url{https://github.com/dmalyuta/arc_2020_code}.}.

The dynamics \eqref{eq:msd_2ndOrder} can be written in the form \eqref{eq:lti} by
using the state $\bm{x}=(r;\dot r)\in\reals^2$, the input $u=f\in\reals$, and
the Jacobians:
\begin{equation}
  \label{eq:msd_lti}
  \tilde A = \Matrix{0 & 1 \\ -\omega_n^2 & -2\zeta\omega_n},~
  \tilde B = \Matrix{0 \\ m\inv}.
\end{equation}

\subsubsection{Zeroth-order Hold}
\label{subsubsection:zoh}

 Zeroth-order hold (ZOH) is a discretization method that assumes
the input to be a staircase signal on the temporal grid. ZOH is called an exact
discretization method because, if the input satisfies this staircase property,
then the discrete-time system state will exactly match the continuous-time system state
at the temporal grid nodes. % More on this after \eqref{eq:zoh_AB}.
In practice, ZOH is a highly relevant discretization type because off-the-shelf
actuators in most engineering domains, including spaceflight, output staircase
commands \citep{scharf2017implementation}.

Optimization routines that use ZOH typically consider a uniform temporal grid,
although the method generally allows for arbitrarily distributed grid nodes:
\begin{equation}
  \label{eq:uniform_temporal_grid}
  t_k=\frac{k-1}{N-1}t_f,~k=1,\dots,\NN.
\end{equation}

The input trajectory is then reduced to a finite number of inputs
$\bm{u}_k\in\RNu$, $k=1,\dots,\NN-1$, that define the aforementioned staircase
signal:
\begin{equation}
  \label{eq:input_staircase}
  \bm{u}(t) = \bm{u}_k,~\forall t\in [t_k,t_{k+1}),~k=1,\dots,\NN-1.
\end{equation}

It then becomes possible to find the explicit update equation for the state
across any $[t_k,t_{k+1}]$ time interval, using standard linear systems theory
\citep{AntsaklisBook}:
\begin{subequations}
  \label{eq:zoh_AB}
  \begin{align}
    \bm{x}_{k+1} &= A\bm{x}_k+B\bm{u}_k, \label{eq:zoh_update} \\
    A &= \Phi(\Delta t_{k},0),~%
        B = A\int_{0}^{\Delta t_k}\Phi(\tau,0)\inv\dd\tau \tilde B,
        \label{eq:AB_integrals}
  \end{align}
\end{subequations}
where $\Delta t_k=t_{k+1}-t_k$, $\bm{x}_k\equiv \bm{x}(t_k)$ and
$\Phi(\cdot,t_k):\reals\to\reals^{\Nx\times\Nx}$ is the \textit{state
  transition matrix}. Since we assumed the system to be LTI, we have
$\Phi(t,t_k)=\expm{\tilde A(t-t_k)}$ where $\expm{}$ is the matrix
exponential. If the system is LTV, the state transition matrix can be computed
by integrating the following ODE:
\begin{equation}
  \label{eq:stm_ode}
  \dot{\Phi}(t,t_k) = \tilde A(t)\Phi(t,t_k),~\Phi(t_k,t_k)=I.
\end{equation}

As it was said before, ZOH is an exact discretization method if the input
behaves according to \eqref{eq:input_staircase}. The reason behind this becomes
clear by inspecting \eqref{eq:zoh_AB}, which \textit{exactly} propagates the
state from one time step to the next. This is different from forward Euler
discretization, where the update is:
\begin{equation}
  \label{eq:forward_euler}
  \bm{x}_{k+1} = \bm{x}_{k}+\Delta t_k(
  \tilde A\bm{x}_k+\tilde B\bm{u}_k),
\end{equation}

For a general non-staircase input signal, however, there is a subtle connection
between ZOH and forward Euler discretization. The former does a zeroth-order
hold on the input signal, and integrates the state exactly, while the latter
does a zeroth-order hold on the output signal (i.e., the time derivative of the
system state). Thus,
unlike in forward Euler discretization, state propagation for ZOH
discretization cannot diverge for a stable system.

The incremental update equation \eqref{eq:zoh_update} can be written in
``stacked form'' to expose how the discrete dynamics are in fact a system of
linear equations. To begin, note that writing \eqref{eq:zoh_update} for
$k=1,\dots,\NN-1$ is mathematically equivalent to:%
\begingroup%
\begin{equation}
  \label{eq:zoh_stacked_explicit}
  \arraycolsep=1.5mm\def\arraystretch{1.2}
  \bm{X} = \left[
    \begin{array}{c|c}
      \begin{matrix}
        I_{\Nx} & 0
      \end{matrix}
      & 0 \\
      \hline
      I_{\NN-1}\otimes A & 0
    \end{array}
  \right]
  \bm{X}+
  \left[
    \begin{array}{c}
      0 \\
      \hline
      I_{\NN-1}\otimes B
    \end{array}
  \right]
  \bm{U}\definedas\mathsf{A}\bm{X}+\mathsf{B}\bm{U},%
\end{equation}%
\endgroup%
where $\bm{X}=(\bm{x}_1;\bm{x}_2;\dots;\bm{x}_{\NN})\in\reals^{\NN\Nx}$ is the
stacked state,
$\bm{U}=(\bm{u}_1;\bm{u}_2;\dots;\bm{u}_{\NN-1})\in\reals^{(\NN-1)\Nu}$ is the
stacked input, and $\otimes$ denotes the Kronecker product. Zeros in
\eqref{eq:zoh_stacked_explicit} denote blocks of commensurate dimensions. Clearly,
we can then write the discrete dynamics as:
\begin{equation}
  \label{eq:zoh_stacked}
  F\bm{X}=G\bm{U}~\text{where}~F\definedas I-\mathsf{A},~G\definedas\mathsf{B}.
\end{equation}

The sparsity pattern for \eqref{eq:zoh_stacked} using the mass-spring-damper
system \eqref{eq:msd_2ndOrder} with $\NN=4$ and $t_f=0.6~\si{\second}$ is shown
in \figref{spy_zoh}. Both $F$ and $G$ consist largely of zeros, which has
important consequences for customizing optimization routines that exploit this
sparsity to speed up computation
\citep{malyuta2019discretization,dueri2017customized}.

\begin{figure}
  \centering
  \includegraphics{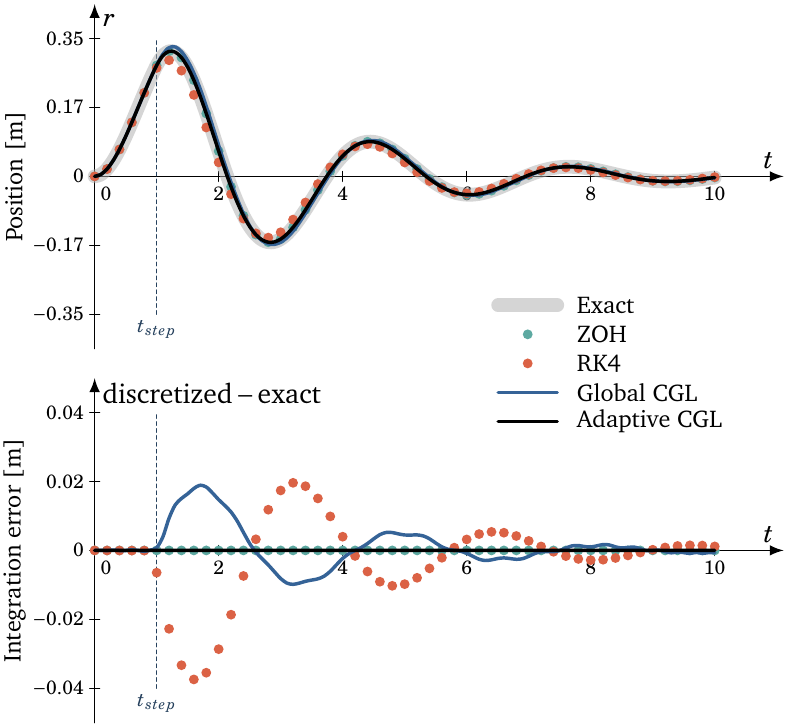}
  \caption{{Simulation of a mass-spring-damper system using four
      discretization types, whose matrix sparsity patterns are shown in
      \figref{spy}. In each case, the simulation final time is
      $t_f=10~\si{\second}$ and a temporal grid of $\NN=51$ is used. The
      pseudospectral simulations are drawn as continuous curves by using
      \eqref{eq:pseudo_state} to recover a continuous-time trajectory from a finite
      set of states at the temporal grid nodes.}}
  \label{fig:msd_simulation}
\end{figure}

An initial value problem (IVP) using a fixed $\bm{x}_1$ can be solved either by
recursively applying \eqref{eq:zoh_update}, or by solving \eqref{eq:zoh_stacked}:
\begin{equation}
  \label{eq:zoh_solution}
  \bm{X}_{2:} = F_{2:}\pinv\bigl(F_1\bm{x}_1+G\bm{U}\bigr),
\end{equation}
where $F_1$ represents the first $\Nx$ columns of $F$, $F_{2:}$ represents the
remaining columns, and $\bm{X}_{2:}=(\bm{x}_2;\dots;\bm{x}_{\NN})$. We use the
left pseudoinverse of $F_{2:}$, and note that the solution is unique since $F_{2:}$
has a trivial nullspace. \figref{msd_simulation} shows an example of applying
\eqref{eq:zoh_update} to the mass-spring-damper system \eqref{eq:msd_2ndOrder}
using $t_f=10~\si{\second}$ and $\NN=51$. Because the input step at
$t_{step}=1~\si{\second}$ falls exactly at a time node, the discretization is
exact.

To connect ZOH discretization back to \pref{ocp}, it is now possible to write
the problem as a finite-dimensional nonlinear parameter optimization:
\begin{optimus}[
  task=\min,
  variables={t_f,\bm{U}},
  objective={L_f(\bm{x}_{\NN},t_f)+\frac{t_f}{\NN-1}\sum_{k=1}^{\NN-1}L(\bm{x}_k,\bm{u}_k,t_k)},
  plabel={ocp_zoh}
  ]
  F\bm{X} = G\bm{U}, \#
  \bm{g}(\bm{x}_k,\bm{u}_k,t_k)\le 0,~\forall k=1,\dots,\NN-1, \#
  \bm{b}(\bm{x}_1,\bm{x}_{\NN},t_{\NN}) = 0.
\end{optimus}

A few remarks are in order about \pref{ocp_zoh}. First, the constraint
\eqref{eq:ocp_zoh_b} exactly satisfies the original dynamics \eqref{eq:ocp_b}
under the ZOH assumption. Second, the optimal solution is only approximately
optimal for \pref{ocp} due to an inexact running cost integration in
\eqref{eq:ocp_zoh_a} and a finite $(\NN-1)$-dimensional basis with which the ZOH
input is constructed. Third, the path constraints \eqref{eq:ocp_c} are satisfied
pointwise in time via \eqref{eq:ocp_zoh_c}, in other words there is a possibility
of inter-sample constraint violation, which can nevertheless sometimes be
avoided \citep{acikmese2008enhancements}. Finally, in many important special
cases \pref{ocp_zoh} is convex, and a multitude of algorithms exploit this fact,
as shall be seen throughout this paper.

\subsubsection{{Runge-Kutta Discretization}}
\label{subsubsection:rk4}

The classic Runge-Kutta (RK4) discretization method can be viewed as an
advanced version of the forward Euler method \eqref{eq:forward_euler}. Unlike ZOH,
which explicitly assumes a staircase input signal, RK4 is a general numerical
integration method that can be applied to any state and control signals. As
such, it is an inexact discretization method like forward Euler, albeit a much
more accurate one. In particular, if we use the uniform temporal grid
\eqref{eq:uniform_temporal_grid}, then the accumulated RK4 integration error
shrinks as $\mathsf{O}(N^{-4})$, whereas forward Euler integration error
shrinks at the much slower rate $\mathsf{O}(N^{-1})$
\citep{betts2010practical,ButcherNumericalBook}.

Directly from the definition of the RK4 method, we can write the following
state update equation:
\begin{subequations}
  \label{eq:rk4_raw}
  \begin{align}
    \bm{x}_{k+1} &= \bm{x}_k+\Delta t_k(
                   \bm{k}_1+2\bm{k}_2+2\bm{k}_3+\bm{k}_4)/6, \\
    \bm{k}_1 &= \tilde A\bm{x}_k+\tilde B\bm{u}_k,
               \label{eq:rk4_k1} \\
    \bm{k}_2 &= \tilde A(\bm{x}_k+0.5\Delta t_k\bm{k}_1)+\tilde B\bm{u}_{k+1/2},
               \label{eq:rk4_k2} \\
    \bm{k}_3 &= \tilde A(\bm{x}_k+0.5\Delta t_k\bm{k}_2)+\tilde B\bm{u}_{k+1/2},
               \label{eq:rk4_k3} \\
    \bm{k}_4 &= \tilde A(\bm{x}_k+\Delta t_k\bm{k}_3)+\tilde B\bm{u}_{k+1},
               \label{eq:rk4_k4}
  \end{align}
\end{subequations}
where $\bm{u}_{k+1/2}=0.5(\bm{u}_k+\bm{u}_{k+1})$. By reshuffling terms in
\eqref{eq:rk4_raw}, we can write the state update in a similar form to
\eqref{eq:zoh_update}:
\begin{equation}
  \label{eq:rk4_update}
  \bm{x}_{k+1} = A\bm{x}_k+B^-\bm{u}_k+B^+\bm{u}_{k+1},
\end{equation}
where $\{A,B^-,B^+\}$ are constructed from $\{I,\tilde A,\tilde B\}$ according
to \eqref{eq:rk4_raw}. Taking inspiration from \eqref{eq:zoh_stacked_explicit},
\eqref{eq:rk4_update} can be written in stacked form:
\begingroup%
\renewcommand*{\arraystretch}{1.2}
\begin{equation}
  \label{eq:rk4_stacked_explicit}
  \bm{X} =
  \left[
    \begin{array}{c|c}
      \begin{matrix}
        I_{\Nx} & 0
      \end{matrix}
      & 0 \\
      \hline
      I_{\NN-1}\otimes A & 0
    \end{array}
  \right]
  \bm{X}+
  \left[
    \begin{array}{c}
      0 \\
      \hline
      \bigl(I_{\NN-1}\otimes [B^-~B^+]\bigr)E
    \end{array}
  \right]
  \bm{U},%
\end{equation}%
\endgroup%
where the matrix $E$ combines columns in order to share the input values
$\bm{u}_k$ for $2\le k\le \NN-1$ (i.e., the internal time grid nodes):
\begin{equation}
  \label{eq:rk4_col_shift_matrix}
  E = \mathrm{blkdiag}\Bigl\{
  I_{\Nu},I_{\NN-2}\otimes\Matrix{I_{\Nu} \\ I_{\Nu}},I_{\Nu}
  \Bigr\}.
\end{equation}

Unlike ZOH, RK4 actually uses the input value at the $\NN$th time node, hence
there is one extra degree-of-freedom (DoF) leading to a slightly larger stacked
input, $\bm{U}=(\bm{u}_1;\bm{u}_2;\dots;\bm{u}_{\NN})\in\reals^{\NN\Nu}$. By
defining $\mathsf{A}$ and $\mathsf{B}$ according to
\eqref{eq:rk4_stacked_explicit}, the discrete dynamics take the same form as
\eqref{eq:zoh_stacked}. In this case, the sparsity pattern for the same
mass-spring-damper example, using $N=4$ and $t_f=0.6~\si{\second}$, is shown in
\figref{spy_rk4}. Like ZOH, RK4 yields a sparse representation of the dynamics.

Like ZOH, an IVP using the discretized dynamics can be solved either by
recursively applying \eqref{eq:rk4_update}, or via a pseudoinverse like
\eqref{eq:zoh_solution}. An example is shown in \figref{msd_simulation}. Clearly,
RK4 is an inexact discretization method. In this case, the interpolated input
$\bm{u}_{k+1/2}$ in \eqref{eq:rk4_k2}-\eqref{eq:rk4_k3} is erroneous just after
the input step at $t_{step}=1~\si{\second}$. Increasing the grid resolution
will quickly decrease the integration error, at the expense of a larger linear
system \eqref{eq:zoh_stacked}.

When discretized with RK4, \pref{ocp} looks much like \pref{ocp_zoh}, except
$\bm{u}_{\NN}$ is an extra decision variable and the user may also choose RK4
to integrate the running cost in \eqref{eq:ocp_a}. However, a subtle but
very important difference is that the solution to the discretized problem
generally no longer \textit{exactly} satisfies the original continuous-time
dynamics. Thus, although RK4 may be computationally slightly cheaper than ZOH
(especially for LTV dynamics, since it does not require computing integrals
like \eqref{eq:AB_integrals}), it is used less often than ZOH or
pseudospectral methods discussed next.

\subsubsection{Pseudospectral Discretization}
\label{subsubsection:pseudospectral}

A key property of ZOH discretization from \sssref{zoh} is that it does not
parametrize the state signal. As a result, numerical integration is required to
recover the continuous-time state trajectory from the solution of
\pref{ocp_zoh}. In trajectory optimization literature, this is known as explicit
simulation or \textit{time marching} \citep{rao2010survey}. An alternative to
this approach is to approximate the state trajectory upfront by a function that
is generated from a finite-dimensional basis of polynomials:
\begin{equation}
  \label{eq:pseudo_state}
  \bm{x}(t) = \sum_{i=1}^\NN \bm{x}_i\phi_i\bigl(\tau(t)\bigr),~t\in [0,t_f],
  \quad
  \phi_i(\tau) \definedas \prod_{\genfrac{}{}{0pt}{}{j=1}{j\neq i}}^{\NN}
  \frac{\tau-\tau_j}{\tau_i-\tau_j},
\end{equation}
where $\tau=2t_f\inv t-1$ and $\phi_i:[-1,1]\to\reals$ are known as Lagrange
interpolating polynomials of degree $\NN-1$. Note that the polynomials are
defined on a normalized time interval. Since the temporal mapping is bijective,
we can equivalently talk about either $t$ or $\tau$.

Given a temporal grid, the Lagrange interpolating polynomials satisfy an
\textit{isolation property}: $\phi_i(\tau_i)=1$ and $\phi_i(\tau_j)=0$ for all
$j\ne i$. Hence, the basis coefficients $\bm{x}_i$ correspond exactly to
trajectory values at the temporal grid nodes. Moreover, the trajectory at all
other times is known \textit{automatically} thanks to \eqref{eq:pseudo_state}. This
is known as implicit simulation or \textit{collocation}. In effect, solving for
the $\NN$ trajectory values at the temporal grid nodes is enough to recover the
complete (approximate) continuous-time trajectory. Because \eqref{eq:pseudo_state}
approximates the state trajectory over the entire $[0,t_f]$ interval, this
technique is known as a \textit{global} collocation.

Collocation conditions are used in order to make the polynomial obtained from
\eqref{eq:pseudo_state} behave according to the system dynamics \eqref{eq:lti}:
\begin{equation}
  \label{eq:collocation_condition}
  \dot{\bm{x}}(t_j)
  =2t_f\inv\sum_{i=1}^\NN \bm{x}_i\phi'_i\bigl(\tau(t_j)\bigr)
  =\tilde A\bm{x}(t_j)+\tilde B\bm{u}(t_j),~\forall j\in\set C,
\end{equation}
where the prime operator denotes differentiation with respect to $\tau$ (i.e.,
$\dd\phi_i/\dd\tau$) and $\set C\subseteq\{1,\dots,\NN\}$ is the set of
collocation points \citep{malyuta2019discretization}. Note that we have already
seen a disguised form of \eqref{eq:collocation_condition} earlier for the RK4
method. In particular, the well-known \eqref{eq:rk4_k1}-\eqref{eq:rk4_k4} are
essentially collocation conditions.

According to the Stone-Weierstrass theorem \citep{BoydChebyshevBook},
\eqref{eq:pseudo_state} approximates a smooth signal with arbitrary accuracy as
$\NN$ is increased. To avoid the so-called \textit{Runge's divergence
  phenomenon}, time is discretized according to one of several special
non-uniform distributions, known as \textit{orthogonal collocations}
\citep{deboor1973collocation}. In this scheme, the grid nodes $\tau_k$ are
chosen to be the roots of a polynomial that is a member of a family of
orthogonal polynomials. For example, Chebyshev-Gauss-Lobatto (CGL) orthogonal
collocation places the scaled temporal grid nodes at the roots of
$(1-\tau)^2c'_{\\N-1}(\tau)=0$, where $c_{\\N}(\tau)=\cos(\NN\arccos(\tau))$ is
a Chebyshev polynomial of degree $\NN$. This particular collocation admits an
explicit solution:
\begin{equation}
  \label{eq:cheb_temporal_grid}
  \tau_k = -\cos\left(\frac{k-1}{\\N-1}\pi \right),~k=1,\dots,\NN.
\end{equation}

A \textit{pseudospectral method} is a discretization scheme that approximates
the state trajectory using \eqref{eq:pseudo_state}, and selects a particular
orthogonal collocation for the collocation points
\citep{rao2010survey,ross2012review,kelly2017introduction}. In fact, the choice
of collocation points is so crucial that flavors of pseudospectral methods are
named after them (e.g., the CGL pseudospectral method). Given this choice,
if the dynamics and control are smooth, the approximation
\eqref{eq:pseudo_state} will converge spectrally (i.e., at an exponential rate in
$\NN$) to the exact state trajectory \citep{rao2010survey}.

Associated with any $\set C$ is a differentiation matrix
$D\in\reals^{|\set C|\times \NN}$ such that $D_{ji}=\phi'_i(\tau_j)$. Some
collocations (e.g., CGL) admit an explicit differentiation matrix, while for
others the matrix can be efficiently computed to within machine rounding error
via barycentric Lagrange interpolation \citep{berrut2004barycentric}. Having
$D$ available allows us to write the collocation conditions
\eqref{eq:collocation_condition} in stacked form:
\begin{equation}
  \label{eq:pseudo_stacked_explicit}
  (2t_f\inv D\otimes I_{\Nx})\bm{X}=(I_{|\set C|}\otimes\tilde A)\bm{X}+
  (I_{|\set C|}\otimes\tilde B)\bm{U},
\end{equation}
where the stacked state and input have the same dimensions as in RK4. We may
thus write the discrete dynamics in the form of \eqref{eq:ocp_zoh_b} by defining:
\begin{subequations}
  \label{eq:global_pseudo_FG}
  \begin{align}
    F &= 2t_f\inv D\otimes I_{\Nx}-I_{|\set C|}\otimes\tilde A, \\
    G &= I_{|\set C|}\otimes\tilde B.
  \end{align}
\end{subequations}

The sparsity pattern for the mass-spring-damper example, using $N=4$ and
$t_f=10~\si{\second}$, is shown in \figref{spy_cgl}. This time, due to $F$ the
dynamics constraint is not sparse. This has historically been a source of
computational difficulty and a performance bottleneck for pseudospectral
discretization-based optimal control
\citep{malyuta2019discretization,sagliano2019generalized}.

Unlike for ZOH and RK4, where an IVP can be solved by recursively applying an
update equation, pseudospectral methods require solving
\eqref{eq:pseudo_stacked_explicit} \textit{simultaneously}, which yields the state
values at the temporal grid nodes all at once. In general, the solution is once
again obtained via the pseudoinverse \eqref{eq:zoh_solution}. However, some
pseudospectral methods such as Legendre-Gauss (LG) and Legendre-Gauss-Radau
(LGR) produce a square and invertible $F_{2:}$ (furthermore,
$F_{2:}\inv F_1=\bm{1}\otimes I_{\Nx}$). This can be used to write
\eqref{eq:ocp_zoh_b} in an ``integral form'' that has certain computation
advantages \citep{francolin2014costate}:
\begin{equation}
  \label{eq:lgr_integral_form}
  \bm{X}_{2:} = -(\bm{1}\otimes I_{\Nx})\bm{x}_1+F_{2:}\inv G\bm{U}.
\end{equation}

Returning to our example of the mass-spring-damper system, \figref{msd_simulation}
shows a simulation using CGL collocation. Like RK4, pseudospectral
discretization is an inexact method, and only approaches exactness as $\NN$
grows large. In this case, the method struggles in particular due to the
discontinuous nature of the input signal, which steps from one to zero at
$t_{step}=1~\si{\second}$. The control trajectory is not smooth due to this
discontinuity, hence the aforementioned spectral convergence guarantee does not
apply. Indeed, it takes disproportionately more grid nodes to deal with this
discontinuity, than if we were to subdivide the simulation into two segments
$t\in [0,t_{step})$ and $t\in [t_{step},t_f]$, where the pre-discontinuity
input applies over the first interval and the post-discontinuity input applies
over the second interval \citep{darby2011hp}.

This idea is at the core of so-called \textit{adaptive}, or \textit{local},
collocation methods
\citep{darby2011hp,sagliano2019generalized,koeppen2019fast,zhao2018dynamics}. These
methods use schemes such as $hp$-adaptation ($h$ and $p$ stand for segment
width and polynomial degree, respectively) in order to search for points like
$t_{step}$, and to subdivide the $[0,t_f]$ interval into multiple segments
according to an error criterion. We defer to the above papers for the
description of these adaptation schemes. For our purposes, suppose that a
partition of $[0,t_f]$ into $S$ segments is available. The $\ell$th segment has
a basis of $N_\ell$ polynomials, a set of collocation points $\set C_\ell$, and
is of duration $t_{s,\ell}$ such that $\sum_{\ell=1}^S t_{s,\ell}=t_f$. Each
segment has its own version of \eqref{eq:global_pseudo_FG}:
\begin{subequations}
  \label{eq:local_pseudo_FGl}
  \begin{align}
    F_\ell &= 2t_{s,\ell}\inv D^\ell \otimes I_{\Nx}-I_{|\set C_\ell|}\otimes\tilde A, \\
    G_\ell &= I_{|\set C_\ell|}\otimes\tilde B,
  \end{align}
\end{subequations}
where the newly defined $F_\ell$ is not to be confused with the earlier $F_1$
and $F_{2:}$ matrices. The new matrices in \eqref{eq:local_pseudo_FGl} are now
combined to write a monolithic dynamics constraint \eqref{eq:ocp_zoh_b}. Doing so
is straightforward for the input, which can be discontinuous across segment
interfaces:
\begin{equation}
  \label{eq:local_pseudo_G}
  G = \mathrm{blkdiag}\{G_1,\dots,G_S\}.
\end{equation}

The state trajectory, however, must remain continuous across segment
interfaces. To this end, the same coefficient $\bm{x}_i$ is used in
\eqref{eq:pseudo_state} for both the final node of segment $\ell$, and the start
node of segment $\ell+1$. Understanding this, we can write:
\begin{equation}
  \label{eq:local_pseudo_F}
  F = \mathrm{blkdiag}\{F_1,\dots,F_S\}E,
\end{equation}
where $E$ combines the columns of $F$ in a similar way to
\eqref{eq:rk4_col_shift_matrix}. The net result is that the final $\Nx$ columns of
$F_{\ell}$ sit above the first $\Nx$ columns of $F_{\ell+1}$.

An example of the sparsity pattern for an adaptive collocation scheme is shown
in \figref{spy_cgl_seg} using $N_1=N_2=2$, $t_{s,1}=t_{step}=1~\si{\second}$, and
$t_f=10~\si{\second}$. One can observe that a second benefit of adaptive
collocation is that it results in a more sparse representation of the dynamics.
This helps to improve optimization algorithm performance
\citep{darby2011hp,sagliano2019generalized}.

Solving an IVP with adaptive collocation works in the same way as for global
collocation. An example is shown in \figref{msd_simulation}, where two segments are
used with the split occuring exactly at $t_{step}$. In this manner, two
instances of \eqref{eq:pseudo_state} are used to approximate the state trajectory,
which is smooth in the interior of both temporal segments. As such, the
approximation is extremely accurate and, for practical purposes, may be
considered exact in this case.

A final, and sometimes very important, aspect of pseudospectral discretization
is that certain collocation schemes yield direct correspondence to the maximum
principle costates of the original optimal control problem (\pref{ocp}). This is
known as the \textit{covector mapping theorem} \citep{gong2007connections}. One
example is the integral form \eqref{eq:lgr_integral_form} for LG and LGR
collocation \citep{francolin2014costate}. Roughly speaking, the Lagrange
multipliers of the corresponding parameter optimization \pref{ocp_zoh} can be
mapped to the costates of \pref{ocp}. Note that this requires approximating the
running cost integral in \eqref{eq:ocp_a} using quadrature weights
$\{w_k\}_{k=1}^{\NN}$ defined by the collocation scheme:
\begin{equation}
  \label{eq:quadrature_cost_integral}
  \int_0^{t_f}L(\bm{x}(\tau),\bm{u}(\tau),\tau)\dd\tau\approx
  \sum_{k=1}^{\NN} w_k L(\bm{x}_k,\bm{u}_k,t_k).
\end{equation}

This unique aspect of pseudospectral methods is why some of the optimal control
problem solvers in \tabref{software} at the end of this article, such as DIDO,
GPOPS-II, and SPARTAN, are listed as both direct and indirect solvers. In fact,
they all solve a discretized version of \pref{ocp}. Nevertheless, they are able to
recover the maximum principle costate trajectories from the optimal solution
\citep{ross2012review,patterson2014gpops,sagliano2019generalized}.

\subsection{Convex Optimization}
\label{subsection:lcvx}

We now come to the question of how to actually solve a finite-dimensional
optimization problem such as \pref{ocp_zoh}. As mentioned in the introduction, this
can be done relatively reliably using well established tools if the problem is
\textit{convex}. Convexity has pervaded optimization algorithm design due to
the following property. If a function is convex, global statements can be made
from local function evaluations. The ramifications of this property cannot be
understated, ranging from the guarantee of finding a global optimum
\citep{RockafellarConvexBook} to precise statements on the maximum iteration
count \citep{peng2002primal}. For the purposes of this review, it is sufficient
to keep in mind that a set $\set C\subseteq\reals^n$ is convex if it contains
the line segment connecting any two of its points:
\begin{equation}
  \label{eq:3}
  x,y\in\set C~\iff~[x,y]_\theta\in\set C~\forall\theta\in [0,1],
\end{equation}
where $[x,y]_\theta\definedas \theta x+(1-\theta)y$. Similarly, a function
$f:\reals^n\to\reals$ is convex if its domain is convex and it lies below the
line segment connecting any two of its points:
\begin{equation}
  \label{eq:3}
  x,y\in\dom(f)~\iff~f([x,y]_\theta)\le[f(x),f(y)]_\theta~\forall\theta\in [0,1].
\end{equation}

Countless resources cover the theory of convex optimization, among which are
the notable books by
\citep{BoydConvexBook,RockafellarConvexBook,NocedalBook}. After applying a
discretization technique akin to those in \ssref{discretization}, a trajectory
design convex optimization problem takes the following form (which is just
another way of writing \pref{ocp_zoh}):
\begin{optimus}[
  result={J^*(t_f)},
  task=\min,
  variables={\bm{U}},
  objective={J(\bm{X}, \bm{U}, t_f)},
  plabel={ocp_convex}
  ]
  \bm{x}_{k+1}=A_k\bm{x}_k+B_k\bm{u}_k+\bm{d}_k,~\forall k=1,\dots,\NN-1, \#
  \bm{g}(\bm{x}_k,\bm{u}_k,t_k)\le 0,~\forall k=1,\dots,\NN-1, \#
  \bm{b}(\bm{x}_1,\bm{x}_N)=0.
\end{optimus}

In \pref{ocp_convex}, $J:\reals^{N\Nx}\times\reals^{(N-1)\Nu}\times\reals\to\reals$ is a convex cost function,
$\bm{g}:\RNx\times\RNu\times\reals\to\RNg$ defines a convex set of constraints,
and $\bm{b}:\RNx\times\RNx\to\RNb$ is an affine function defining the trajectory
boundary conditions. If $t_f$ is a decision variable, a sequence of
\pref{ocp_convex} instances can be solved using a line search that computes
$\min_{t_f}J^*(t_f)$ \citep{blackmore2010minimum}. The sequence
$\{A_k,B_k,\bm{d}_k\}_{k=1}^{\NN-1}$ of matrices of commensurate dimensions
represents the linear time-varying discretized dynamics \eqref{eq:ocp_b}. In
numerous aerospace applications, including rocket landing and spacecraft rendezvous,
the dynamics are at least approximately of this form
\citep{acikmese2007convex,DeRuiterBook}.

The path constraints \eqref{eq:ocp_convex_c} are where nonconvexity appears
most often for a space vehicle trajectory optimization problems. Sometimes the
nonconvexity can be removed by a clever reformulation of the problem, a process
called \textit{convexification}. If the reformulation is exact, in other words
the solution set is neither reduced nor expanded, the convexification is
\textit{lossless}. One example of lossless
convexification that has pervaded rocket landing literature is a thrust
lower-bound constraint. Let $\bm{T}(t)\in\reals^3$ be a thrust
vector, then combustion stability and engine performance dictate the following
constraint:
\begin{equation}
  \label{eq:thrust_constraint}
  \rho_{\min}\le\norm[2]{\bm{T}(t)}\le\rho_{\max},~\forall t\in [0,t_f].
\end{equation}

The lower-bound is nonconvex, but it was shown that it admits the following
lossless convexification \citep{acikmese2007convex}:
\begin{equation}
  \label{eq:thrust_lcvx}
  \rho_{\min}\le\sigma(t)\le\rho_{\max},~\norm[2]{\bm{T}(t)}\le\sigma(t),~\forall t\in [0,t_f].
\end{equation}

\begin{figure}
  \centering
  \includegraphics{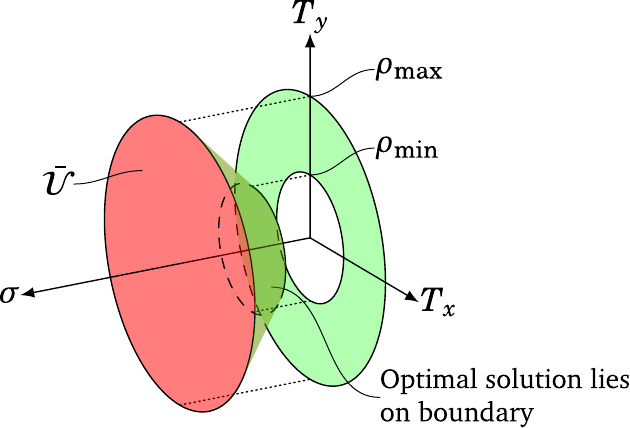}
  \caption{Illustration of the convex relaxation technique used throughout much
    of lossless convexification literature for powered descent guidance. Using
    the maximum principle, lossless convexification proves that the optimal
    solution $(\bm{T}^*(t); \sigma^*(t))$ lies on the green boundary of the set
    $\bar{\mathcal U}$.}
  \label{fig:lcvx}
\end{figure}

The convexification ``lifts'' the input space into an extra dimension, as
illustrated in \figref{lcvx}. Clearly the lifted feasible set $\bar{\set U}$ is
convex, and its projection onto the original coordinates contains the feasible
set defined by \eqref{eq:thrust_constraint}. The backbone of lossless
convexification is a proof via the maximum principle that the optimal solution
lies on the boundary of $\bar{\set U}$, as highlighted in \figref{lcvx}. Thus, it
can be shown that the solution using \eqref{eq:thrust_lcvx} is optimal for the
original problem which uses \eqref{eq:thrust_constraint}.

Another example of lossless convexification comes from the constrained
reorientation problem. Let \(\bm{q}(t)\in\reals^4\) with
\(\norm[2]{\bm{q}(t)}=1\) be a unit quaternion vector describing the attitude
of a spacecraft. During the reorientation maneuver, it is critical that
sensitive instruments on the spacecraft are not exposed to bright celestial
objects. This dictates the following path constraint:
\begin{equation}
  \label{eq:quaternion_constraint}
  \bm{q}(t)\T M\bm{q}(t)\leq 0,~\forall t\in [0,t_f],
\end{equation}
where $M\in\reals^{4\times 4}$ is a symmetric matrix that is not positive
semidefinite, making the constraint nonconvex. However, when considered
together with the implicit constraint \(\norm[2]{\bm{q}(t)}=1\),
\eqref{eq:quaternion_constraint} can be losslessly replaced with the following
convex constraint \citep{kim2004quadratically}:
\begin{equation}
  \label{eq:quaternion_constraint_lcvx}
  \bm{q}(t)\T (M+\mu I)\bm{q}(t)\leq \mu,~\forall t\in [0,t_f],
\end{equation}
where $\mu$ is chosen such that the matrix \(M+\mu I\) is positive
semidefinite. Instead of the maximum principle, the proof of this lossless
convexification hinges on the geometry of the unit quaternion.

\subsection{Sequential Convex Programming}
\label{subsection:scp}

Sequential convex programming (SCP) is an umbrella name for a family of
nonconvex local optimization methods. It is one of many available tools
alongside nonlinear programming, dynamic programming, and genetic algorithms
\citep{FloudasEncyclopedia}. If lossless convexification is a surgical knife to
remove acute nonconvexity, SCP is a catch-all sledgehammer for nonconvex
trajectory design. Clearly many aerospace problems are nonconvex, and SCP has
proven to be a competitive solution method for many of them
\citep{szmuk2018successive,liu2014solving,bonalli2017analytical}. This section
provides an intuition about how SCP algorithms work as well as their advantages
and limitations. The interested reader can find further information in
\citep{SCPTrajOptCSM2021} which provides a comprehensive tutorial on SCP.

\begin{figure}
  \centering
  \includegraphics{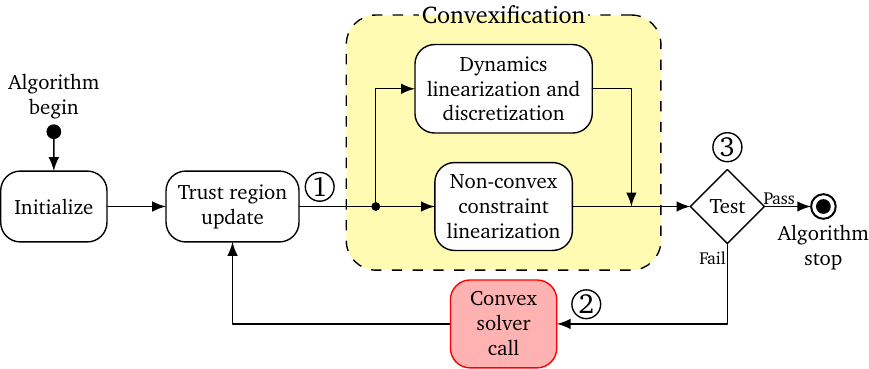}
  \caption{Block diagram illustration of a typical {SCP} algorithm. The forward
    path can be seen as a ``predictor'' step, while the reverse path that calls
    the convex optimization solver can be seen as a ``corrector''
    step. Although the test criterion can be guaranteed to trigger for certain
    SCP flavors, the solution may not be feasible for the original problem.}
  \label{fig:scvx_block_diagram}
\end{figure}

At the core, every SCP algorithm is based on the following idea: iteratively
solve a convex approximation of \pref{ocp}, and update the approximation as new
solutions are obtained. \figref{scvx_block_diagram} provides an illustration and
highlights how SCP can be thought of as a \textit{predictor-corrector}
algorithm. In the forward predictor path, the current solution is evaluated for
its quality. If the quality check fails, the reverse corrector path improves
the quality by solving a \textit{subproblem} that is a better convex
approximation. Examples of SCP algorithms include cases where the subproblem is
linear \citep{palaciosgomez1982nonlinear}, second-order conic
\citep{mao2018successive}, and semidefinite \citep{fares2002robust}.

\newcommand{\alglocation}[1]{\tikz[baseline=-0.8ex]{
    \node[
    circle,
    fill=white,
    draw=black,
    minimum size=0.6em,
    inner sep=0.1mm
    ] {#1};
  }}

Consider \pref{ocp_convex} for a simple concrete example of the SCP
philosophy. Assume that $\bm{g}$ is the only nonconvex element and that $t_f$
is fixed. At the location \alglocation{1} in \figref{scvx_block_diagram}, the SCP
method provides an iterate in the form of a current trajectory guess
$\{\bar{\bm{x}}_k,\bar{\bm{u}}_k\}_{k=1}^{\NN-1}$. In its most basic form, SCP
linearizes and relaxes the $\bm{g}$ function:
\begin{equation}
  \label{eq:10}
  \bar{\bm{g}}+
  \frac{\subdiff\bar{\bm{g}}}{\subdiff\bm{x}}\Delta\bm{x}_k+
  \frac{\subdiff\bar{\bm{g}}}{\subdiff\bm{u}}\Delta\bm{u}_k\le\bm{\alpha}_k,~
  \forall k=1,\dots,\NN-1,
\end{equation}
where $\bm{\alpha}_k\in\RNg$ is a \textit{virtual buffer zone} and we define
$\bar{\bm{g}}\definedas\bm{g}(\bar{\bm{x}}_k,\bar{\bm{u}}_k,t_k)$,
$\Delta\bm{x}_k\definedas\bm{x}_k-\bar{\bm{x}}_k$, and
$\Delta\bm{u}_k\definedas\bm{u}_k-\bar{\bm{u}}_k$. The subproblem solved at
location \alglocation{2} in \figref{scvx_block_diagram} is then given by:
\begin{optimus}[
  task=\min,
  variables={\substack{\bm{u}_1,\dots,\bm{u}_{\NN-1}\\\bm{\alpha}_1,\dots,\bm{\alpha}_{\NN-1}\\\eta_1,\dots,\eta_{\NN-1}}},
  objective={J(\bm{X}, \bm{U}, t_f)+w_{\text{vc}}\bm{1}\T\sum_{k=1}^{\NN-1}[\bm{\alpha}_k]^++
    w_{\text{tr}}\sum_{k=1}^{\NN-1}\eta_k},
  plabel={ocp_scvx}
  ]
  \bm{x}_{k+1}=A_k\bm{x}_k+B_k\bm{u}_k+\bm{d}_k,~\forall k=1,\dots,\NN-1, \#
  \bar{\bm{g}}+\frac{\subdiff\bar{\bm{g}}}{\subdiff\bm{x}}\Delta\bm{x}_k+
  \frac{\subdiff\bar{\bm{g}}}{\subdiff\bm{u}}\Delta\bm{u}_k\le\bm{\alpha}_k,~
  \forall k=1,\dots,\NN-1, \#
  \norm{\Delta\bm{u}_k}\le\eta_k,~\forall k=1,\dots,\NN-1, \#
  \bm{b}(\bm{x}_1,\bm{x}_N)=0.
\end{optimus}

\pref{ocp_scvx} introduces several new elements. The variables $\eta_k$ regulate
the size of \textit{trust regions} around the previous solution, and the
weights $w_{\text{tr}},~w_{\text{vc}}\in\reals$ are set to large positive
values that encourage convergence and zero constraint violation. The best
choice of $p$-norm $\norm{\cdot}$ in \eqref{eq:ocp_scvx_d} depends on the problem
structure. The stopping criterion used in \alglocation{3} of
\figref{scvx_block_diagram} may be, for example:
\begin{equation}
  \label{eq:11}
  \max_{k\in\{1,\dots,\NN-1\}}\eta_k\le\epsilon~\textnormal{and}~
  \max_{k\in\{1,\dots,\NN-1\}}\norm[\infty]{[\bm{\alpha}_k]^+}\le\epsilon,
\end{equation}
where $\epsilon$ is a user-chosen convergence tolerance constant that can be
interpreted as a small ``numerical error''.

SCP denotes a \textit{family} of solution methods and, as such, countless
variations of \pref{ocp_scvx} exist. Early versions of SCP for trajectory
generation focused on motion kinematics alone \citep{schulman2014motion} or
included dynamics but with few convergence guarantees
\citep{augugliaro2012generation}. Today, a family of methods is emerging with
stronger convergence guarantees, including \scvx \citep{mao2018successive},
GuSTO \citep{bonalli2019gusto,BonalliLewTAC2021}, and penalized trust region
(PTR) \citep{reynolds2020real}. \pref{ocp_scvx} exemplifies the PTR method, where
the trust region sizes $\eta_k$ are themselves optimization variables that are
kept small using a penalty in the cost \eqref{eq:ocp_scvx_a}. PTR is often the
fastest method, but its theoretical convergence properties are relatively
unexplored. In comparison, \scvx and GuSTO provide a guarantee that the
algorithm converges to a locally optimal solution, albeit with potentially
non-zero $\eta_k$ and $\bm{\alpha}_k$. When these variables are zero, however,
the solution is locally optimal for the original optimal control problem.

The main algorithmic differences across SCP methods lie in how the convex
approximations are formulated, what methods are used to update the intermediate
solutions and to measure progress towards optimality, and how all of this lends
itself to theoretical analysis. For example, \scvx uses a discrete-time
convergence proof while GuSTO uses the continuous-time maximum principle. The
main difference with the PTR method is that both \scvx and GuSTO update
$\eta_k$ outside of the optimization problem. Interestingly, the PTR method has
been observed to yield much faster convergence in practice, and a theoretical
explanation recently appeared \citep{reynolds2020crawling}.

\subsubsection{{Related Algorithms}}

In the general context of optimization, SCP belongs to the class of so-called
trust region methods \citep{NocedalBook,ConnTrustRegionBook}. However, SCP is
not to be confused with another popular trust region method, sequential
quadratic programming (SQP). First of all, SCP solves its subproblems to full
optimality. While this increases the number of iterations in the reverse path
of \figref{scvx_block_diagram}, it vastly reduces the number of forward
passes. Owing to the growing maturity of IPM solvers and the advent of solver
customization \citep{domahidi2013ecos,dueri2014automated}, iterations in the
reverse path are relatively ``cheap'', making the trade-off a good one. Second,
SCP requires only first-order problem information, since nonconvexities are
handled by a simple linearization such as in \eqref{eq:10}. On the other hand, SQP
is a second-order method that requires the factorization of a Hessian. This
raises concerns about matrix positive semidefiniteness and may require
computationally expensive techniques such as the BFGS update
\citep{gill2011sequential}.

Differential dynamic programming (DDP) is another family of algorithms that,
like SCP, is built around the idea of linearization
\citep{mayne1966ddp,jacobson1968new}. More precisely, DDP solves a discrete-time
optimal control problem with an additive cost function like the one in
\eqref{eq:ocp_zoh_a}. Although DDP falls outside the scope of this survey paper, we
will mention that it has major successful applications in space vehicle
trajectory optimization and provides an interesting variation of the
linearize-and-solve framework of \figref{scvx_block_diagram}. DDP is particularly
popular for low-thrust orbit transfer trajectory optimization. For example, the
NASA Mystic software used DDP for low-thrust trajectory optimization of the
Dawn Discovery mission to the Vesta and Ceres protoplanets of the asteroid belt
\citep{whiffen2001application,whiffen2006mystic}. Other appearances of DDP
include include multi-revolution and multi-target orbit transfer
\citep{lantoine2012hybrid1,lantoine2012hybrid2}, Earth to Moon transfer with an
exclusion zone \citep{pellegrini2020multiple1,pellegrini2020multiple2}, and
low-thrust flyby trajectory planning to near-Earth objects
\citep{colombo2009optimal}.

A disadvantage of the original DDP algorithm is that it is an unconstrained
optimization method. This means that while SCP naturally handles state and
input constraints like \eqref{eq:ocp_convex_c}, implementing such constraints is
still an active research area for DDP. Most attempts to incorporate constraints
make use of penalty, barrier, augmented Lagrangian, and active set methods
\citep{tassa2014control, xie2017differential}. Most recently, extensions of DDP
were proposed to handle general nonconvex state and input constraints using a
primal-dual interior point method
\citep{pavlov2020interior,aoyama2020constrained}.

Another disadvantage of DDP is that it is a second-order method. Like SQP, this
makes DDP more computationally expensive than SCP which only requires
first-order information. Nevertheless, there is numerical evidence that DDP is
faster than SQP \citep{lantoine2012hybrid2}. Furthermore, related algorithms
have been developed that only use first-order information, such as the
iterative linear quadratic regulator (iLQR). The ALTRO software is a popular
modern trajectory optimization toolbox based on the iLQR and augmented
Lagrangian methods \citep{howell2019altro}.

\subsection{Mixed-integer Programming}
\label{subsection:mip}

Mixed-integer programming (MIP) solves problems where some decision variables
are binary. Consider a concrete example in the context of \pref{ocp}. Without loss
of generality, suppose that the control input is partitioned into continuous
and binary variables:
\begin{equation}
  \label{eq:12}
  \bm{u}(t) = \Matrix{\bm{v}(t)\\\bm{\zeta}(t)}\in\RNu,~\bm{v}(t)\in\reals^{\Nu-n_\zeta},~%
  \bm{\zeta}(t)\in\binary^{n_\zeta}.
\end{equation}

Binary variables naturally encode discrete events such as the opening of valves
and relays, the pulsing of space thrusters, and mission phase transitions
\citep{bemporad1999control,sun2019multi}. Furthermore, binary variables can help
to approximate nonlinear gravity, aerodynamic drag, and other salient features
of a space vehicle trajectory optimization problem
\citep{blackmore2012lossless,marcucci2019mixed}.

 To formally discuss how MIP might be relevant for a spacecraft
trajectory optimization problem like \pref{ocp}, consider the space vehicle to
be an \textit{autonomously switched hybrid system}
\citep{saranathan2018relaxed}. In particular, suppose that the vehicle dynamics
\eqref{eq:ocp_b} and its constraints \eqref{eq:ocp_c} are continuous except for
the following extra ``if-then'' condition:
\begin{equation}
  \label{eq:ifthen}
  \bm{q}(\bm{z})<0~\implies~\bm{c}(\bm{z})=0,
\end{equation}
where $\bm{z}\in\RNz$ is some mixture of inputs, states, and time. In this
formulation, the \textit{constraint function} $\bm{c}:\RNz\to\RNc$ is activated
if the \textit{trigger function} $\bm{q}:\RNz\to\reals^{n_\zeta}$ maps to the
negative orthant $\reals^{n_\zeta}_{<0}$. The conditional statement
\eqref{eq:ifthen} can be formulated as the following set of mixed-integer
constraints:
\begin{subequations}
  \label{eq:mip_ifthen}
  \begin{align}
    -\zeta_i M
    &\le q_i(\bm{z})\le (1-\zeta_i)M,\quad i=1,\dots,n_\zeta,
      \label{eq:mip_q} \\
    -\bigl(1-\sigma(\bm{\zeta})\bigr)M\ones
    &\le \bm{c}(\bm{z}) \le \bigl(1-\sigma(\bm{\zeta})\bigr)M\ones,
      \label{eq:mip_c} \\
    \sigma(\bm{\zeta}) &= \prod_{i=1}^{n_\zeta}\zeta_i,
                         \label{eq:mip_activation_function}
  \end{align}
\end{subequations}
where $M\in\reals$ is a sufficiently large positive number. The function
$\sigma:\binary^{n_\zeta}\to\binary$ is called the \textit{activation
  function}, and it imposes the if-then logic of \eqref{eq:ifthen} through
\eqref{eq:mip_q}-\eqref{eq:mip_c}. When the binary variable $\zeta_i$ equals one, the
$i$th trigger is activated. Thus, when $\sigma(\bm{\zeta})=1$, the left-hand
side of \eqref{eq:ifthen} holds. In fact, $\bm{q}(\bm{z})= 0$ is also possible, but
this case is irrelevant since an optimal solution will not activate the
constraint function unnecessarily.

{Mixed-integer programming can be used to solve \pref{ocp} in the presence of the
  constraints \eqref{eq:mip_ifthen}}. Traditional MIP solvers are based on the
branch-and-bound method \citep{NemhauserIntegerBook,CookIntegerBook}. At their
core is a divide-and-conquer logic that often, though not always, speeds up the
solution process by eliminating large numbers of possible $\bm{\zeta}$
combinations. Modern MIP solvers also improve runtime through methods like
pre-solving (which reduces $n_\zeta$), cutting planes (which introduce new
constraints to tighten the feasible space), heuristics, parallelism, branching
variable selection, symmetry detection, and so on
\citep{achterberg2007constrained,achterberg2013mixed}. In the worst case,
however, MIP runtime remains exponential in $n_\zeta$. This is a large
hindrance to onboard implementation, since space vehicle hardware is often not
able to support the large MIP computational demand
\citep{malyuta2020fast,malyuta2019lossless,malyuta2019approximate}.

As usual in optimization, one can trade the global optimality of MIP for
solution speed by accepting local optimality or by approximating the precise
statement \eqref{eq:ifthen} with a more efficient formulation. In the following
subsections, we will introduce two popular approaches that have recently
emerged in both direct and indirect solution methods for solving MIP problems
without introducing integer variables.

\subsubsection{{State-triggered Constraints}}
\label{subsubsection:stc}

\begin{figure}
  \centering
  \includegraphics{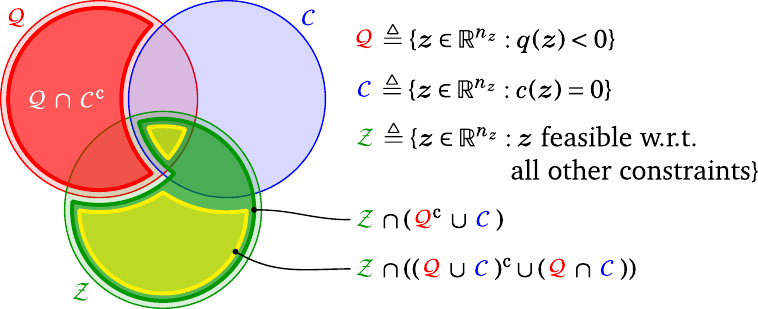}
  \caption{A Venn diagram of the sets of solution variables $\bm{z}\in\reals^{n_z}$
    that satisfy the trigger condition ${\color{red}\set Q}$, constraint
    condition ${\color{blue}\set C}$, and all other constraints
    ${\color{green!60!black}\set Z}$. STCs ensure that feasible solutions
    satisfy
    $\bm{z}\notin{\color{red}\set Q}\intersection\setcomplement{{\color{blue}\set C}}$
    (illustrated by the red set). The
    feasible set with an STC is shaded green, and the feasible set with the
    bidirectional constraint \eqref{eq:lcp_implication} is shaded in yellow. The
    sets with bold outlines are padded to help visual separation only.}
  \label{fig:stc}
\end{figure}

{\textit{State-triggered constraints} (STCs) take the direct solution
  approach, and are under active study using the SCP framework from \ssref{scp}
  \citep{szmuk2018successive,malyuta2020fast}}. Roughly speaking, STCs embed the
discrete if-then logic from \eqref{eq:ifthen} into a continuous direct
formulation with minimal penalty to the solution time
\citep{reynolds2019dual}. In its most basic form, an STC models \eqref{eq:ifthen}
for the scalar case $n_\zeta=1$ and $\Nc=1$. While there are several useful
theoretical connections between STCs and the linear complementarity problem
(LCP) \citep{LCPBook}, STCs encode a much larger feasible space than LCP
constraints \citep{szmuk2018successive}. Namely, STCs only encode forward
implications, and they are \textit{not} bidirectional statements of the
following form:
\begin{equation}
  \label{eq:lcp_implication}
  q(\bm{z})<0~\iff~c(\bm{z})=0.
\end{equation}

\figref{stc} illustrates the distinction between \eqref{eq:ifthen} and
\eqref{eq:lcp_implication}. The green set denotes the feasible set with the STC,
while the yellow set denotes the feasible set of the more restrictive
constraint \eqref{eq:lcp_implication}. {Clearly, the feasible space is
  larger when using the STC, and this can translate into a more optimal
  solution.}

Continuing our discussion for the scalar case, it can be shown that
\eqref{eq:ifthen} is equivalent to either one of the following two continuous
constraints:
\begin{subequations}
  \label{eq:stc_formulations}
  \begin{align}
    \label{eq:stc_nonunique_formulation}
    &q(\bm{z})+\sigma\ge 0,~\sigma\ge 0,~\sigma\cdot c(\bm{z})=0,~\text{or} \\
    &-\min(0,q(\bm{z}))\cdot c(\bm{z})=0,
  \end{align}
\end{subequations}
where $\sigma\in\reals$ is a slack variable {that plays the role of the
  activation function from \eqref{eq:mip_activation_function}}. Although both
constraints in \eqref{eq:stc_formulations} are nonconvex, they are readily ingested
by the SCP linearization process described in \ssref{scp}.

 A notable feature of STCs is that they readily extend to the
multivariable case of \eqref{eq:ifthen}, and have been shown to handle both AND
and OR combinations of triggers and constraints
\citep{szmuk2019real,szmuk2019successive}:
\begin{subequations}
  \label{eq:compound_stcs}
  \begin{align}
    \label{eq:compound_and_and}
    \bigwedge_{i=1}^{n_\zeta} q_i(\bm{z})<0~
    &\implies~ \bigwedge_{i=1}^{\Nc} c_i(\bm{z})=0, \\
    \label{eq:compound_or_and}
    \bigvee_{i=1}^{n_\zeta} q_i(\bm{z})<0~
    &\implies~ \bigwedge_{i=1}^{\Nc} c_i(\bm{z})=0, \\
    \label{eq:compound_and_or}
    \bigwedge_{i=1}^{n_\zeta} q_i(\bm{z})<0~
    &\implies~ \bigvee_{i=1}^{\Nc} c_i(\bm{z})=0, \\
    \label{eq:compound_or_or}
    \bigvee_{i=1}^{n_\zeta} q_i(\bm{z})<0~
    &\implies~ \bigvee_{i=1}^{\Nc} c_i(\bm{z})=0.
  \end{align}
\end{subequations}

In the general context of optimization, STCs do for trajectory optimization
what the $\set S$-procedure from linear matrix inequalities (LMIs) does for
stability analysis and controller synthesis \citep{LMIBook}, and what
sum-of-squares (SOS) programming does to impose polynomial non-negativity over
basic semialgebraic sets \citep{majumdar2017funnel}. {That is, STCs
  embed an otherwise difficult logic constraint into a tractable continuous
  formulation}. In particular, note that the scalar version of \eqref{eq:ifthen}
can be written as:
\begin{equation}
  \label{eq:stc_forall}
  c(\bm{z})=0\quad\forall \bm{z}~\text{s.t.}~q(\bm{z})<0.
\end{equation}

On the other hand, the $\set S$-procedure and SOS programming consider the
following constraints, respectively:
\begin{subequations}
  \label{eq:stc_alternatives}
  \begin{align}
    f_0(\bm{z})
    &\ge 0\quad\forall \bm{z}~\text{s.t.}~f_i(\bm{z})\ge 0,~i=1,\dots,p,
      \label{eq:s_procedure} \\
    p(\bm{z})
    &\ge 0\quad\forall \bm{z}~\text{s.t.}~p_{\mathrm{eq}}(\bm{z})=0,~
      p_{\mathrm{ineq}}(\bm{z})\ge 0,
      \label{eq:sos}
  \end{align}
\end{subequations}
where $f_i$, $i=0,\dots,p$, are quadratic functions, while $p$,
$p_{\mathrm{eq}}$, and $p_{\mathrm{ineq}}$ are polynomials. Comparing
\eqref{eq:stc_forall} with \eqref{eq:stc_alternatives} makes the connection to
STCs clear.

\subsubsection{{Homotopy Methods}}
\label{subsubsection:indirect_mip}

Homotopy methods, also known as \textit{numerical continuation schemes}, were
previously introduced in \sssref{indirect} in the context of solving standard
optimal control problems. It turns out that homotopy can also be used to encode
\eqref{eq:ifthen} in a continuous framework, and has been successfully embedded
into recent indirect trajectory optimization algorithms. In this section, we
briefly introduce the relaxed autonomously switched hybrid system (RASHS) and
composite smooth control (CSC) methods
{\citep{saranathan2018relaxed,taheri2020novel,arya2021composite}}.

To begin, let $\sigma$ denote the activation function from
\eqref{eq:mip_activation_function}. Using the third equation from
\eqref{eq:stc_nonunique_formulation}, we note that \eqref{eq:ifthen} is exactly
equivalent to the following constraint:
\begin{equation}
  \label{eq:ifthen_equivalent}
  \sigma(\bm{\zeta})\bm{c}(\bm{z})=0.
\end{equation}

\begin{figure}
  \centering
  \begin{subfigure}[t]{0.5\columnwidth}
    \centering
    \includegraphics{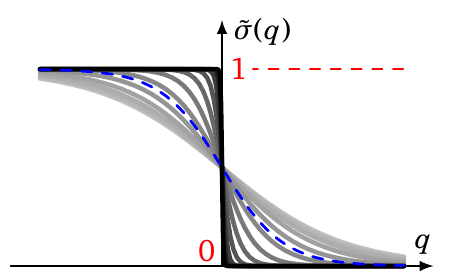}
    \caption{The RASHS sigmoid \eqref{eq:rashs_sigmoid}.}
    \label{fig:rashs_sigmoid}
  \end{subfigure}%
  \hfill%
  \begin{subfigure}[t]{0.5\columnwidth}
    \centering
    \includegraphics{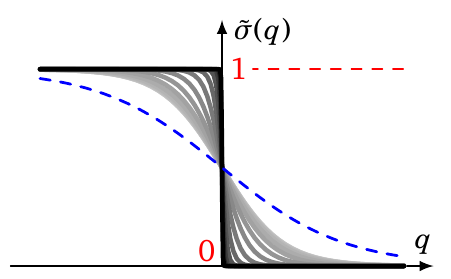}
    \caption{The CSC sigmoid \eqref{eq:csc_sigmoid}.}
    \label{fig:csc_sigmoid}
  \end{subfigure}
  \caption{{Comparison of the RASHS and CSC sigmoids, as defined in
      \eqref{eq:sigmoids}. A sweep is shown from homotopy parameter $\kappa=1$
      (lighter colors) to $\kappa=100$ (darker colors). The blue dashed curve
      shows the alternative sigmoid for $\kappa=1$ (i.e., CSC for
      (\protect\subref{fig:rashs_sigmoid}) and RASHS for
      (\protect\subref{fig:csc_sigmoid})). As $\kappa$ increases, the sigmoid
      quickly converges to an accurate approximation of the binary activation
      function in \eqref{eq:mip_activation_function}. While the nature of both
      sigmoids is similar, for a given $\kappa$ the CSC sigmoid is more
      localized around the $y$-axis, and hence is a closer approximation of a
      step signal.}}
  \label{fig:sigmoids}
\end{figure}

At the core of the RASHS and CSC methods is an approximation of the binary
function $\sigma$ by a continuous sigmoid function
$\tilde\sigma:\reals^{n_\zeta}\to [0,1]$:
\begin{subequations}
  \label{eq:sigmoids}
  \begin{align}
    \label{eq:rashs_sigmoid}
    \text{RASHS:}~
    &\tilde\sigma(\bm{q})=\prod_{i=1}^{n_{\zeta}}\bigl(1+
      \exp{\kappa q_i}\bigr)\inv, \\
    \label{eq:csc_sigmoid}
    \text{CSC:}~
    &\tilde\sigma(\bm{q})=\prod_{i=1}^{n_{\zeta}}\frac 12
      \bigl(1-\tanh(\kappa q_i)\bigr),
  \end{align}
\end{subequations}
where the latter equation uses the theory of hyperbolic tangent smoothing
\citep{taheri2018smoothing}. The homotopy parameter $\kappa\in [0,\infty)$
regulates the accuracy of the approximation, with increasing accuracy as
$\kappa$ grows, such that
$\lim_{\kappa\to\infty}\tilde\sigma(\bm{q})=\sigma(\bm{\zeta})$. \figref{sigmoids}
illustrates how both sigmoid functions evolve as $\kappa$ increases. The
core idea of RASHS and CSC is to begin with a small $\kappa$ where the optimal
control problem is continuous and ``easy'' to solve, and to judiciously
increase $\kappa$ to such a large value that the solution becomes
indistinguishable from its MIP counterpart.

It is worth noting the specific instances of \eqref{eq:ifthen} considered by
RASHS and CSC. The former method was developed to compute time- or
state-triggered multiphase trajectories, where vehicle dynamics change across
phases (e.g., stage separation during rocket ascent)
\citep{saranathan2018relaxed}. Such a system is also known as a
\textit{differential automaton} \citep{tavernini1987differential}. In this case,
we can have $m$ constraints of the form \eqref{eq:ifthen}, where the $k$th
constraint is:
\begin{equation}
  \label{eq:rashs_ifthen_time}
  \bm{q}^k(t)<0~\implies~\dot{\bm{x}}(t)-\bm{f}^k(\bm{x}(t),\bm{u}(t),t)=0,
\end{equation}
and $\bm{q}^k:\reals\to\reals^{n_{\zeta}^k}$ indicates the time interval where
the $k$th dynamics apply. Assuming that the time intervals do not overlap, we
can sum the smoothed versions of \eqref{eq:rashs_ifthen_time} to obtain a single
continuous system dynamics constraint:
\begin{equation}
  \label{eq:rashs_smooth_dynamics}
  \dot{\bm{x}}(t)=
  \biggl[
  \sum_{k=1}^m\tilde\sigma\bigl(\bm{q}^k(t)\bigr)
  \biggr]\inv
  \sum_{k=1}^m\tilde\sigma\bigl(\bm{q}^k(t)\bigr)\bm{f}^k[t].
\end{equation}

Note that the new dynamics \eqref{eq:rashs_smooth_dynamics} are a convex
combination of the individual dynamics over the $m$ time intervals. As $\kappa$
is increased, the approximation becomes more accurate, and the correct
$\bm{f}^k$ functions begin to dominate their respective intervals. Moreover,
using \eqref{eq:rashs_smooth_dynamics} instead of \eqref{eq:rashs_ifthen_time} has the
algorithmic advantage of replacing a multi-point BVP with a TPBVP.

The CSC method, on the other hand, considers systems with fixed dynamics but
multiple control inputs or constraints \citep{taheri2020novel}. In both cases,
the overall control input can be expressed as a function of $m$ ``building
block'' inputs $\bm{u}^k$, such that:
\begin{equation}
  \label{eq:csc_ifthen_input}
  \bm{q}^k(\bm{x}(t),\bm{u}(t),t)<0~\implies~\bm{u}(t)=\bm{u}^k(t),
\end{equation}
where $\bm{q}^k:\RNx\times\RNu\times\reals\to\reals^{n_\zeta^k}$ are mutually
exclusive indicators of when the $k$th building block input applies. Like for
RASHS, the following equation provides a smooth approximation of the control,
from which CSC derives its name:
\begin{equation}
  \label{eq:csc_smooth_control}
  \bm{u}(t)=
  \biggl[
  \sum_{k=1}^m\tilde\sigma\bigl(\bm{q}^k[t]\bigr)
  \biggr]\inv
  \sum_{k=1}^m\tilde\sigma\bigl(\bm{q}^k[t]\bigr)\bm{u}^k(t).
\end{equation}

Clearly, RASHS, CSC, and STCs are all approaching the same problem of
efficiently handling \eqref{eq:ifthen} from subtly different angles. It is worth
noting that for the moment, both RASHS and CSC can only handle the AND
combination \eqref{eq:compound_and_and} of trigger and constraint functions. Most
recently, \citep{MalyutaAcikmeseJGCD2021} showed that a similar homotopy
framework can handle OR combinations of trigger functions. This opens up an
interesting research avenue to develop a unifying homotopy method that handles
all the logical combinations in \eqref{eq:compound_stcs}.

\subsection{Model Predictive Control}
\label{subsection:mpc}

\begin{figure}
  \centering
  \includegraphics{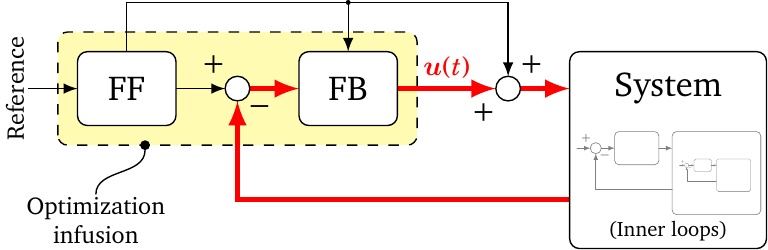}
  \caption{A typical control architecture consists of nested layers of
    feedforward (FF) and feedback (FB) elements. The execution frequency
    increases going from the outermost to the innermost layers. In particular,
    elements in the FB path (highlighted in red) have stricter execution time
    requirements than FF elements.}
  \label{fig:control_mental_model}
\end{figure}

The preceding sections focused on solving one instance of \pref{ocp}. We now place
ourselves in the context of a control system whose architecture is illustrated
in \figref{control_mental_model}. Two important algorithm categories that are part
of a control system are so-called \textit{feedforward} and \textit{feedback}
\citep{LurieBook}, and optimization-based methods can potentially be applied to
both. In the feedback path, the current state estimate of the system is used to
continually update the control signal, which means that \pref{ocp} must be
re-solved many times. This is the core idea of model predictive control (MPC).

In its most basic form, an MPC formulation of \pref{ocp_convex} can be
expressed as follows:
\begin{optimus}[
  result={\bm{u}^*_1},
  task=\argmin,
  variables={\bm{u}_1,\dots,\bm{u}_{\NN-1}},
  objective={\bm{x}_\NN\T Q_f\bm{x}_\NN+\sum_{k=1}^{\NN-1}
    \bm{x}_k\T Q\bm{x}_k+\bm{u}_k\T R\bm{u}_k},
  plabel={ocp_mpc}
  ]
  \bm{x}_{k+1}=A_k\bm{x}_k+B_k\bm{u}_k+\bm{d}_k,~\forall k=1,\dots,\NN-1, \#
  \bm{g}(\bm{x}_k,\bm{u}_k,t_k)\le 0,~\forall k=1,\dots,\NN-1, \#
  \bm{x}_1=\hat{\bm{x}},~\bm{b}(\bm{x}_N)=0.
\end{optimus}

\begin{figure}
  \centering
  \includegraphics{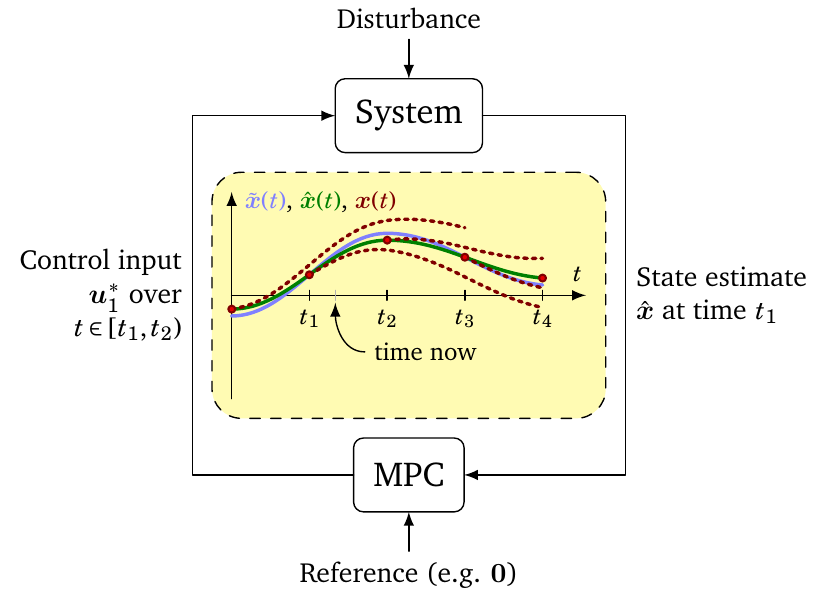}
  \caption{{Block diagram illustration of an MPC controller. At each
      time step $t_k$, MPC computes the optimal control input $\bm{u}_1^*$ by
      using a mathematical model of the system and solving \pref{ocp_mpc}, which
      is a receding horizon optimal control problem. Note the three states
      drawn in the diagram: the actual state $\tilde{\bm{x}}$, the estimated
      state $\hat{\bm{x}}$, and the internally propagated MPC state
      $\bm{x}$. Each state may be slightly different due to estimation error,
      model uncertainty, and disturbances.}}
  \label{fig:mpc_illustration}
\end{figure}

{\figref{mpc_illustration} illustrates how \pref{ocp_mpc} can be used to
  control a dynamical system}. Note that \pref{ocp_mpc} is a parametric
optimization problem because it depends on the current state estimate
$\hat{\bm{x}}\in\RNx$. The first optimal control input $\bm{u}_1^*$ for
\pref{ocp_mpc} becomes $\bm{u}(t)$ in \figref{control_mental_model}. The weight
matrices $Q\mgeq 0$ and $R\mgr 0$ in the running cost and the terminal weight matrix
$Q_f\mgeq 0$ are chosen to get a desired response. Together with the terminal
constraint \eqref{eq:ocp_mpc_d}, these choices must ensure stability and
recursive feasibility in closed-loop operation (i.e., the problem must be
feasible the next time that it is solved).

The main advantage of MPC is that it is arguably the most natural methodology
for handling system constraints in a feedback controller
\citep{mayne2000constrained}. However, because MPC operates in a feedback loop,
stability and performance are both critical and strongly dependent on
uncertainty robustness and execution frequency
\citep{LurieBook,MIMOBook}. Troves of information have been compiled on the
subject, which remains an active research area. Numerous surveys on MPC cover
general and future methods \citep{mayne2014model}, robustness
\citep{garcia1989model,bemporad2007robust,mayne2015robust}, computational
requirements \citep{alessio2009survey}, and industrial applications
\citep{eren2017model,mao2018survey,dicairano2018real,qin2003survey}. For space
vehicle applications in particular, where onboard computation is limited, we
single out so-called explicit MPC
\citep{bemporad2002explicit,RawlingsMPCBook,BorrelliMPCBook}. The concept is to
pre-compute a lookup table for the solution of \pref{ocp_mpc}. This turns out to be
possible to do exactly when the MPC problem is a QP, and approximately in more
general cases up to MICP \citep{malyuta2019approximate}. When onboard storage
and problem dimensionality permit, explicit MPC yields a much faster and
computationally cheaper algorithm in which onboard optimization is replaced by
a static lookup table.

%%%%%%%%%%%%%%%%%%%%%%%%%%%%%%%%%%%%%%%%%%%%%%%%%%%%
\section{Applications}
\label{section:applications}
%%%%%%%%%%%%%%%%%%%%%%%%%%%%%%%%%%%%%%%%%%%%%%%%%%%%

This section describes the application of optimization methods from the previous
section to state-of-the-art space vehicle control problems. The following
subsections cover the following key areas of spaceflight. \ssref{pdg} discusses
rocket powered descent for planetary landing. \ssref{rendezvous} covers spacecraft
rendezvous and \ssref{smallbody} covers a closely related problem of small body
landing. \ssref{reorientation} talks about attitude
reorientation. Endo-atmospheric ascent and entry are surveyed in
\ssref{endo}. Last but not least, orbit transfer is discussed in \ssref{orbit}.

\subsection{Powered Descent Guidance for Rocket Landing}
\label{subsection:pdg}

Powered descent guidance (PDG) is the terminal phase of EDL spanning the last
few kilometers of altitude. The goal is for a lander to achieve a soft and
precise touchdown on a planet's surface by using its rocket engine(s). PDG
technology is fundamental for reducing cost and enabling access to hazardous
yet scientifically rich sites
\citep{starek2015challenges,carson2019splice,steinfeldt2010guidance,braun2006mars,jones2018recent,robertson2017synopsis,europa2012study,artemis2019moon,lroc2018lunar}. The
modern consensus is that iteration-based algorithms within the CGC paradigm,
rather than closed-form solutions, are required for future landers
\citep{lu2017cgc,carson2019splice}. The survey of applications in this section
demonstrates that optimization offers a compelling iteration-based solution
method due to the availability of real-time algorithms that can enforce
relevant PDG constraints.

To place state-of-the-art PDG into context, let us briefly mention some key
heritage methods. Initial closed-form algorithms are known as
\defintext{explicit guidance}, which is characterized by directly considering
the targeting condition each time the guidance command is generated
\citep{lu2020theory}. Early algorithms solved a version of the following OCP:
\begin{optimus}[
  task={\min},
  variables={t_f,\bm{a}},
  objective={\int_0^{t_f}\bm{a}(t)\T\bm{a}(t)\dt},
  label={dsouza}
  ]
  \ddot{\bm{r}}(t) = \bm{g}+\bm{a}(t), \#
  \bm{r}(0)=\bm{r}_0,~\dot{\bm{r}}(0)=\dot{\bm{r}}_0,~
  \bm{r}(t_f)=\bm{r}_f,~\dot{\bm{r}}(t_f)=\dot{\bm{r}}_f.
\end{optimus}

Here, $\bm{r}(t)\in\reals^3$ denotes position, $\bm{a}(t)\in\reals^3$ is the
acceleration control input, $\bm{g}\in\reals^3$ is the gravitational
acceleration vector and $t_f$ is the flight duration. Position and velocity
boundary values are fixed. The optimal solution to \pref{dsouza} is known as the
E-Guidance (EG) law \citep{cherry1964general,dsouza1997optimal}:
\begin{equation}
  \label{eq:e_guidance}
  \bm{a}(t) = 6 t_{\text{go}}^{-2}\bm{ZEM}(t)-2 t_{\text{go}}\inv\bm{ZEV}(t),
\end{equation}
where $t_{\text{go}}\definedas t_f-t$ is the time-to-go and:
\begin{subequations}
  \begin{align}
    \bm{ZEM}(t) &\definedas \bm{r}_f-\left(
                  \bm{r}(t)+t_{\text{go}}\dot{\bm{r}}(t)+0.5 t_{\text{go}}^2\bm{g}
                  \right), \\
    \bm{ZEV}(t) &\definedas \dot{\bm{r}}_f-\left(
                  \dot{\bm{r}}(t)+t_{\text{go}}\bm{g}
                  \right),
  \end{align}
\end{subequations}
are the zero-effort-miss and zero-effort-velocity terms
\citep{furfaro2011nonlinear,song2020survey}. Nominally, \eqref{eq:e_guidance}
results in an affine acceleration profile. If instead one allows the
acceleration profile to be quadratic, an additional DoF appears, which is fixed
by setting the terminal acceleration $\bm{a}(t_f)=\bm{a}_f$. This results in
the Apollo guidance (APG) law, which flew on the historic Lunar missions
\citep{klumpp1974apollo}:
\begin{equation}
  \label{eq:4}
  \bm{a}(t) = 12 t_{\text{go}}^{-2}\bm{ZEM}(t)-6 t_{\text{go}}\inv\bm{ZEV}(t)+\bm{a}_f.
\end{equation}

The concept of an acceleration profile behind EG and APG has since been
extended and generalized, resulting in a \defintext{polynomial guidance} family
of algorithms. \citep{zhang2017collision} augment the cost \eqref{eq:dsouza_a} with
a surface collision-avoidance term. \citep{guo2013waypoint} formulate a QP to
solve for an intermediate waypoint that augments collision-avoidance
capabilites as well as enforces actuator saturation for thrust- and
power-limited engines. \citep{lu2019augmented} develops a general theory for
polynomial guidance laws that contains EG and APG as special cases. For one of
the best modern explanations of polynomial guidance methods, the reader should
consult \citep{lu2020theory}. Unfortunately, closed-form polynomial guidance is
unable to handle many operational constraints \citep{lu2018propellant} and is
not fuel optimal since the cost \eqref{eq:dsouza_a} rather penalized energy.

\begin{figure}
  \centering
  \includegraphics{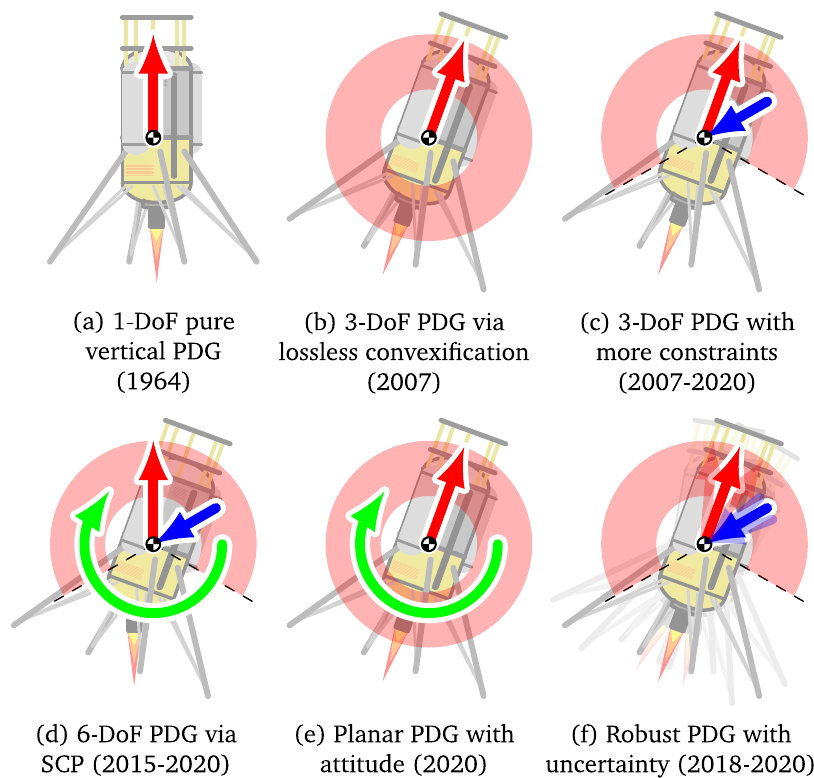}
  \caption{Progression of PDG problem complexity. The red, green, and blue arrows
    denote thrust, torque, and aerodynamic force respectively. The red region
    denotes the feasible thrust set.}
  \label{fig:pdg_cases}
\end{figure}

To overcome these limitations, research has long sought to characterize and
eventually solve the more general fuel-optimal PDG problem. The first milestone
towards fuel-optimal PDG was a closed-form single-DoF vertical descent solution
\citep{meditch1964problem}, illustrated in \figref{pdg_cases}a. Evidence
suggests that Apollo landings came close to this optimum
\citep{klumpp1974apollo,MindellBook}. The maximum principle
\citep{PontryaginBook} played a key role back then, and continues to do so in the
present day.

Seeking to generalize the single-DoF result, Lawden formulated the necessary
conditions of optimality for 3-DoF PDG \citep{LawdenBook,MarecBook}. However,
solving the necessary conditions requires shooting methods, which are typically
too computationally expensive and sensitive to the initial guess to allow
efficient onboard implementation \citep{betts1998survey}. More recently,
\citep{topcu2005fuel,topcu2007minimum} extended the results from
\citep{LawdenBook} to the case of angular velocity control, and compared the
solution quality of fuel-optimal 3-DoF PDG to the necessary conditions of
optimality. However, the aim of the work was not real-time onboard
implementation, so nonlinear programming (SQP) was used via the GESOP
solver.% which is computationally expensive and does not
% have convergence guarantees.

\begin{figure}
  \centering
  \includegraphics{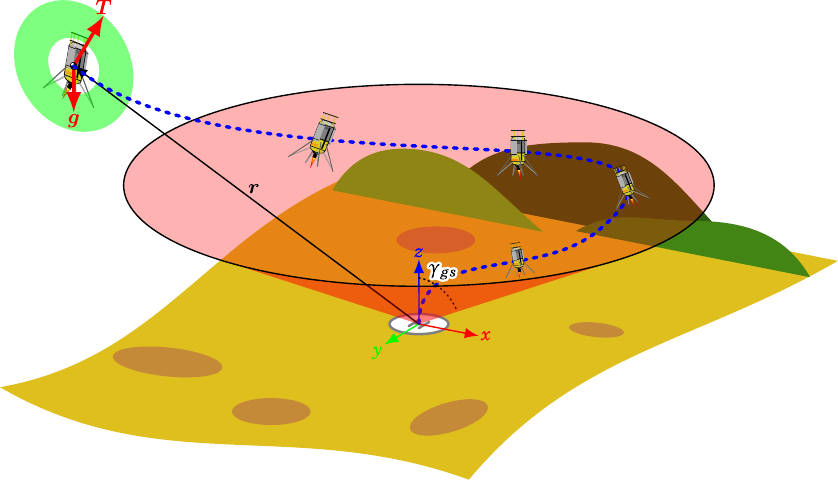}
  \caption{Illustration of the basic powered descent guidance solved by
    \pref{pdg2007} via lossless convexification. The goal is to safely bring the
    rocket lander to standstill on the landing pad while satisfying the thrust
    magnitude constraints and maintaining a minimum glideslope.}
  \label{fig:pdg_diagram}
\end{figure}

After decades of research into problem \textit{characterization}, a watershed
moment for problem \textit{solution} came in the mid 2000s with
the % in the landmark
papers \citep{acikmese2005powered,acikmese2007convex}. The authors solved the
following 3-DoF PDG problem, illustrated in \figref{pdg_diagram}, via the
process of lossless convexification described in \ssref{lcvx}:
\begin{optimus}[
  task=\min,
  variables={t_f,\bm{T}},
  objective={\int_0^{t_f}\norm[2]{\bm{T}(t)}\dt},
  label={pdg2007}
  ]
  \ddot{\bm{r}}(t) = \bm{g}+\bm{T}(t)m(t)\inv, \#
  \dot m(t) = -\alpha\norm[2]{\bm{T}(t)}, \#
  \rho_{\min}\le \norm[2]{\bm{T}(t)}\le\rho_{\max}, \#
  \bm{r}(t)\T \hat{\bm{e}}_z\ge\norm[2]{\bm{r}(t)}\cos(\gamma_{gs}), \#
  m(0)=m_0,~\bm{r}(0)=\bm{r}_0,~\dot{\bm{r}}(0)=\dot{\bm{r}}_0,~\bm{r}(t_f)=0,~\dot{\bm{r}}(t_f)=0.
\end{optimus}

\begin{figure}
  \centering
  \includegraphics{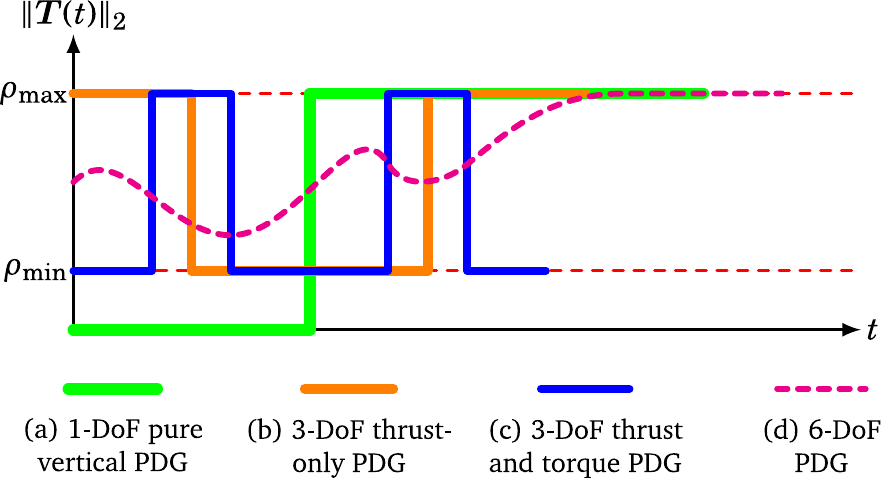}
  \caption{Optimal thrust profiles for several powered descent guidance
    formulations. (a) Corresponds to the classical single-DoF result by
    \citep{meditch1964problem}; (b) corresponds to {3-DoF translation-only}
    landing from \citep{LawdenBook,acikmese2007convex}; (c) corresponds to
    planar landing with rotation from \citep{reynolds2020optimal}. {The thrust
      profile for general 6-DoF PDG with translation and rotation is an open
      problem. In particular, there is no theory to guarantee that it should be
      bang--bang, thus (d) shows a profile with no clear structure.}}
  \label{fig:pdg_thrust}
\end{figure}

Unlike the classical \pref{dsouza}, \pref{pdg2007} readily handles several important
operational constraints, including thrust bounds \eqref{eq:pdg2007_d} and glide
slope \eqref{eq:pdg2007_e}. Through numerical simulations for a prototype Mars
lander, \citep{acikmese2007convex} confirmed that the optimal thrust has a
max-min-max profile as shown in \figref{pdg_thrust}. This profile was proven to be
optimal for 3-DoF PDG in \citep{LawdenBook,topcu2007minimum}.

Over the course of the next decade, the method was expanded to handle fairly
general nonconvex input sets \citep{acikmese2011lossless}, minimum-error
landing and thrust pointing constraints
\citep{blackmore2010minimum,carson2011lossless,acikmese2013lossless}, classes of
affine and quadratic state constraints
\citep{harris2013losslessacc,harris2013losslesscdc,harris2014lossless,harris2013maximum},
classes of nonlinear (mixed-integer) dynamics \citep{blackmore2012lossless},
certain binary input constraints \citep{malyuta2019lossless,HarrisTAC2021}, and conservative
conic obstacles \citep{bai2019optimal}.

The maturity of a method can be gauged by the availability of a precise
statement of its limits, similar to the role played by the Bode integral in
frequency-domain control \citep{MIMOBook,LurieBook}. Such a characterization
appeared for lossless convexification in the form of constrained reachable or
controllable sets \citep{eren2015constrained,dueri2016consistently,DueriThesis}
or ``access'' conditions \citep{song2020survey}. These sets, obtained
numerically and with arbitrarily high precision, define the boundary conditions for
which versions of \pref{pdg2007} are feasible.
% A simplified illustration of how the problem evolved is shown in
% Figures~\ref{fig:pdg_cases}b and \ref{fig:pdg_cases}c.

The practicality of lossless convexification-based PDG methods was demonstrated
through a series of flight tests conducted by the NASA Jet Propulsion
Laboratory, Masten Space Systems, and university partners. In a 3-year 7-flight
test campaign, the Masten Xombie sounding rocket demonstrated that robust
onboard real-time optimization is feasible on spaceflight processors
\citep{acikmese2013flight,scharf2014adapt,scharf2017implementation,flights2012a,flights2012b}. A
number of publications accompanied this flight test campaign, including a
comparison of lossless convexification to polynomial guidance
\citep{ploen2006comparison}, onboard computation time reduction via
time-of-flight interpolation \citep{scharf2015interpolation}, and complete
off-line lookup table generation \citep{acikmese2008enhancements}. The resulting
algorithm, G-FOLD \citep{GFOLDpatent,acikmese2012gfold}, solves a tailored
version of \pref{pdg2007} using a custom C-language SOCP solver called Bsocp
\citep{dueri2014automated}. G-FOLD is able to compute rocket landing
trajectories in 100~\si{\milli\second} on a 1.4~\si{\giga\hertz} Intel Pentium
M processor. Further evidence of real-time performance was presented by
\citep{dueri2017customized}, where Bsocp ran on a radiation-hardened BAE RAD750
PowerPC.

Despite the significant flight envelope expansion afforded by lossless
convexification
\citep{ploen2006comparison,wolf2012improving,carson2011capabilities}, an
inherent limitation of 3-DoF PDG is that the computed trajectory cannot
incorporate attitude constraints other than those on the thrust vector, which serves
as an attitude proxy. An extensive simulation campaign is required to validate
the 3-DoF trajectory to be executable by a fundamentally 6-DoF lander system
\citep{carson2019splice,KamathAAS2020}. Thus, recent PDG research has sought 6-DoF formulations that are able
to incorporate attitude dynamics and constraints.

The SCP family of methods, discussed in \ssref{scp}, has emerged as an effective
approach to transition from a fully convex 3-DoF problem to a 6-DoF problem
with some nonconvexity. Some of the popular SCP algorithms include \scvx
\citep{mao2018successive}, penalized trust region \citep{reynolds2020real}, and
GuSTO \citep{bonalli2019gusto,BonalliLewTAC2021}. Some other algorithms based
around similar ideas have also emerged, such as ALTRO which is based on
iterative LQR \citep{howell2019altro}.

A vast number of flavours of SCP exist, however, since it is a nonlinear
optimization technique that works best when tailored to exploit problem
structure. In certain cases, lossless convexification is embedded to remove some
nonlinearity.
% \citep{liu2019fuel,simplicio2019guidance,li2020online,wang2019optimal,szmuk2016successive}.
\citep{liu2019fuel} convexifies an angle-of-attack (AoA) constraint relating to
an aerodynamic control capability, \citep{simplicio2019guidance} solve a
version of \pref{pdg2007} in a first step and passes the solution to a second
step involving SCP, while
\citep{li2020online,wang2019optimal,szmuk2016successive} use the classical
convexification result for the thrust magnitude constraint
\eqref{eq:thrust_lcvx}.

Since the mid 2010s, SCP technology enabled the expression of quadratic
aerodynamic drag and thrust slew-rate constraints \citep{szmuk2016successive},
attitude dynamics \citep{szmuk2017successive}, variable time-of-flight
\citep{szmuk2018successive-conf}, and an ellipsoidal drag model that allows
aerodynamic lift generation along with variable ignition time
\citep{szmuk2018successive}. % Because SCP is a nonlinear programming algorithm
% without strong convergence guarantees, several
Several papers on SCP ``best practices'' have also appeared, including thrust
input modeling \citep{szmuk2017successive}, the effect of discretization on
performance~\citep{malyuta2019discretization}, and using dual quaternions to
alleviate nonconvexity in the constraints by off-loading it into the
dynamics~\citep{reynolds2019state}. Practical details on real-time
implementation are also available \citep{reynolds2020real}, where the SCP
solution is compared to the globally optimal trajectory for a planar landing
problem \citep{reynolds2020optimal}. Most recently, a comprehensive tutorial
paper with open-source code was published, and describes the algorithmic and
practical aspects of SCP methods and of lossless convexification
\citep{SCPTrajOptCSM2021}.

\figref{pdg_cases} summarizes the dominant directions of PDG development since
2005. Starting from the classical vertical-descent result by
\citep{meditch1964problem}, \figref{pdg_cases}a, the early breakthrough for
practical onboard PDG solution was achieved in 2007 by
\citep{acikmese2007convex}, \figref{pdg_cases}b. Since then, 3-DoF PDG methods
have been extended and flight tested, \figref{pdg_cases}c. {In
  particular, more complicated effects such as aerodynamic drag force were
  added by these extensions, which are listed in the preceding
  paragraph. Perhaps the biggest modern shift in PDG technology development has been
  to consider attitude dynamics, which is motivated by the inability to impose
  non-trivial attitude constraints in a 3-DoF formulation
  \citep{carson2019splice}. This has led to a family of so-called 6-DoF PDG
  algorithms, \figref{pdg_cases}d, that often rely on SCP methods}. To compare
how close SCP comes to the global optimum, recent work found optimal solutions
for {``planar'' PDG \citep{reynolds2020optimal},
  \figref{pdg_cases}e. This work restricts the landing trajectory to a 2D plane,
  but does include attitude dynamics. Therefore it represents both a
  generalization of \figref{pdg_cases}b and a restriction of \figref{pdg_cases}d,
  and provides new insight into the 6-DoF PDG optimal solution
  structure}. Today, PDG research evolves along the following broad directions:
guaranteeing real-time performance, convergence, and solution quality, handling
binary constraints, and incorporating uncertainty as shown in
\figref{pdg_cases}f.

One exciting development for SCP in recent years has been the advent of
state-triggered constraints, introduced in \sssref{stc}. This allows real-time
capable embedding of if-then logic into the guidance problem. To demonstrate
the capability, \citep{szmuk2018successive} imposed a velocity-triggered AoA
constraint, \citep{reynolds2019state} imposed a distance-triggered line-of-sight
constraint, \citep{szmuk2019successive} imposed a collision-avoidance
constraint, and \citep{reynolds2019dual} imposed a slant-range-triggered
line-of-sight constraint. In particular, the latter two works develop a theory
of \textit{compound} STCs that apply Boolean logic to combine multiple trigger
and constraint functions, as shown in \eqref{eq:compound_stcs}. The impact of STCs
on the ability to compute solutions in real-time is discussed in
\citep{szmuk2018successive,reynolds2019dual}.

Simultaneously with the development of SCP for PDG, the pseudospectral
discretization community has produced a rich body of work investigating the
solution quality benefits of that method. Building on foundational early work
\citep{rao2010survey,fahroo2002direct,garg2010unified,kelly2017introduction}, it
was demonstrated for a variant of \pref{pdg2007} that pseudospectral methods yield
greater solution accuracy with fewer temporal nodes
\citep{sagliano2018pseudospectral}. However, as discussed in
\sssref{pseudospectral}, pseudospectral methods traditionally yield slower solution
times because they generate non-sparse matrices for the discretized equations
of motion \citep{malyuta2019discretization}. By using an $hp$-adaptive scheme
inspired by the finite element method \citep{darby2011hp}, it was shown that
this can be somewhat circumvented
\citep{sagliano2018generalized-conf,sagliano2019generalized}. Furthermore, it
was shown that pseudospectral discretization within an SCP framework yields
solutions up to 20 times faster than using sequential quadratic programming
\citep{wang2018pseudospectral}.

As deterministic PDG algorithms mature, research is becoming increasingly
interested in making the trajectory planning problem robust to various sources
of uncertainty. One approach is to design a feedback controller to correct for
deviations from the nominal trajectory, such that the overall control input is
given by:
\begin{equation}
  \label{eq:2}
  \bm{u}(t) = \bar{\bm{u}}(t)+K(t)(\bm{x}(t)-\bar{\bm{x}}(t)),
\end{equation}
where $\bar{\bm{x}}(t)$ and $\bar{\bm{u}}(t)$ are the nominal state and control
respectively, and $K(t)\in\reals^{\Nu\times\Nx}$ is a feedback gain matrix. In
\citep{ganet2019optimal,scharf2017implementation}, the feedback controller is
designed separately from the nominal trajectory. However, incorporating
feedback law synthesis into the nominal trajectory generation problem can
achieve more optimal solutions \citep{garciasanz2019control}. This
``simultaneous'' feedback-feedforward design was done via multi-disciplinary
optimization in \citep{jiang2018computational}, desensitized optimal control in
\citep{shen2010desensitizing,seywald2019desensitizedoverview},
chance-constrained optimization in \citep{ono2015chance}, and covariance
steering in \citep{ridderhof2018uncertainty,ridderhof2019minimum}. Other work in
this domain includes open-loop robust trajectory design via Chebyshev interval
inclusion \citep{cheng2019uncertain}, and \textit{a posteriori} statistical
analysis through linear covariance propagation \citep{woffinden2019linear} and
Monte Carlo simulation \citep{scharf2017implementation}.

PDG methods based on lossless convexification and SCP are in most cases
\textit{implicit guidance} methods. In this setup, the targeting condition
(e.g., soft touchdown on the landing pad) is met by tracking a reference
trajectory that yields the correct terminal state. Functionally, PDG methods
are most often situated in the FF block of \figref{control_mental_model}, and they
generate a complete trajectory upfront that is tracked by a feedback
controller. From a systems engineering perspective, this has a clear advantage
of allowing heritage control methods to perform the intricate and critical
control of the actual vehicle. However, it was mentioned at the start of this
section and in \ssref{mpc} that continually re-solving for the PDG trajectory can
offer additional robustness. In contrast to traditional polynomial guidance,
some modern approaches aim to leverage this robustness and also satisfy system
constraints via model predictive control.

\citep{cui2012receding} show how to leverage MPC for landing with an uncertain
state and variable gravitational field, while \citep{wang2019optimal} show how
to ensure recursive feasibility and a bounded guidance error by executing a
nominal and relaxed optimization problem in parallel. In both methods, the full
trajectory optimization problem is solved from the current state to the final
landing location, thus the MPC horizon ``shrinks'' throughout the PDG
maneuver. A more traditional approach is taken by \citep{lee2017constrained},
where the prediction horizon extends for a finite duration beyond the current
state. The authors also show that difficult constraints on sensor line-of-sight
and spacecraft attitude are convex using a dual quaternion
representation. Numerical performance of MPC for PDG on an embedded ARM platform
was documented in \citep{pascucci2015model}.

\subsection{Rendezvous and Proximity Operations}
\label{subsection:rendezvous}

\begin{figure}
  \centering
  \includegraphics{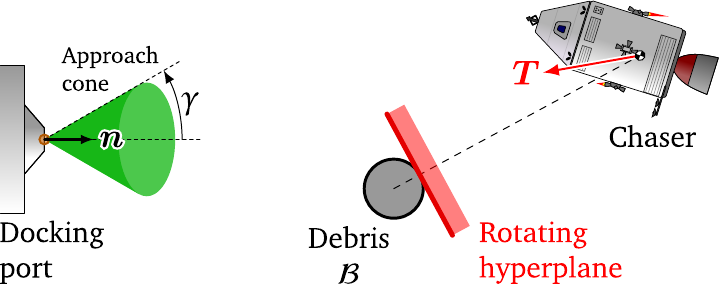}
  \caption{Illustration of a basic rendezvous scenario. Roughly speaking, the
    goal is for a chaser spacecraft to use the thrust $\bm{T}$ from its reaction
    control system to dock with a target while avoiding debris and respecting
    constraints such as staying within the approach cone.}
  \label{fig:rendezvous}
\end{figure}

Let us now switch contexts from the final stages of planetary landing to the
realm of orbital spaceflight. A key task for a spacecraft in orbit is to
perform rendezvous and proximity operations (RPO). The goal is to bring an
actively controlled chaser vehicle and a passively controlled target vehicle to
a prescribed relative configuration, in order to achieve mission objectives
such as inspection or docking. A detailed overview of RPO history and
technology development can be found in
\citep{FehseBook,goodman2006history,woffinden2007navigating,luo2014survey}. This
section focuses on the challenges and developments in RPO using convex
optimization-based solution methods.

Throughout this section we consider the following RPO trajectory optimization
problem, illustrated in \figref{rendezvous}:
\begin{optimus}[
  task={\min},
  variables={t_f,\bm{T}},
  objective={\int_0^{t_f}\norm[2]{\bm{T}(t)}\dt},
  label={rendezvous}
  ]
  \ddot{\bm{r}}(t)= -\mu\norm[2]{\bm{r}(t)}\inv[3]\bm{r}(t)+m(t)\inv\bm{T}(t), \#
  \dot{m}(t)=-\alpha\norm[2]{\bm{T}(t)}, \#
  \norm[2]{\bm{T}(t)}\leq \rho, \#
  \bm{r}(t)\notin \mathcal{B}(t),\#
  \norm[2]{\bm{r}(t)-\hat{\bm{r}}(t)}\cos\gamma\leq (\bm{r}(t)-\hat{\bm{r}}(t))\T \bm{n}(t), \#
  m(0)= m_0,~\bm{r}(0)=\bm{r}_0, ~ \dot{\bm{r}}(0)=\dot{\bm{r}}_0, \#
  \bm{r}({t_f})=\hat{\bm{r}}(t_f), ~ \dot{\bm{r}}(t_f)=\hat{\bm{r}}(t_f),
\end{optimus}
where ${\bm{r}}(t),~\hat{\bm{r}}(t)\in\reals^3$ denote
the positions of the chaser and target spacecraft in the inertial
frame. The basic objective in \eqref{eq:rendezvous_a} is to minimize fuel consumption
\citep{park2013model}. Other choices include sparsification of the control
sequence \citep{hartley2013terminal}, trading off flight duration with fuel
consumption \citep{hu2018trajectory}, encouraging smoothness of the control
sequence \citep{li2018model}, and reducing the sensitivity to sensing and
control uncertainties \citep{jin2020robust}. We note that \pref{rendezvous} only
characterizes the last phase of RPO. The reader is referred to
\citep{hartley2012model,sun2019multi} for examples of multi-phase RPO trajectory
optimization.

Since the quantity of interest in RPO is the relative motion between the chaser
and the target, it is commonplace to express the dynamics
\eqref{eq:rendezvous_b} in a different reference frame. Examples include the
local-vertical local-horizontal (LVLH) frame centered at the target, or a
line-of-sight polar reference frame \citep{li2017line}. Based on this choice,
different models of relative dynamics have been studied, and are surveyed in
\citep{sullivan2017comprehensive}. For near-circular orbits, linear
time-invariant Hill-Clohessy-Wiltshire (HCW) equations are the most popular
model \citep{clohessy1960terminal}. For elliptical orbits, the linear
time-varying Yamanaka-Ankerson (YA) state transition matrix is the usual choice
\citep{yamanaka2002new}. Perhaps a cleaner approach is to avoid relative
dynamics by working in the inertial frame, as done in
\eqref{eq:rendezvous_b}. \citep{lu2013autonomous,liu2014solving} showed that fast
and reliable trajectory optimization is still possible in this case, by
applying the same lossless convexification as in \pref{pdg2007} to the
constraints \eqref{eq:rendezvous_c} and \eqref{eq:rendezvous_d} {and
  successively linearizing the dynamics \eqref{eq:rendezvous_b}.
  \citep{benedikter2019convex} further proposed a filtering technique for
  updating the linearization reference point to improve the algorithm
  robustness}. The advantage of this approach is its compatibility with general
Keplerian orbits and perturbations like \(J_2\) harmonic and aerodynamic drag.

One key challenge in RPO is to avoid collision with external debris or part of
the target vehicle itself, which is described by constraint
\eqref{eq:rendezvous_e}. One approach to enforcing \eqref{eq:rendezvous_e} is to
pre-compute a so-called virtual net of trajectories that allows to avoid
obstacles in real-time via a simple graph search
\citep{frey2017constrained,weiss2015safe}. The pre-computation procedure,
however, may be prohibitively computationally demanding. In comparison, solving
\pref{rendezvous} directly can avoid virtual net construction altogether if an
efficient solution method is available. To this end, the keep-out zone
\(\mathcal{B}(t)\) is usually chosen to be a polytope, an ellipsoid, or the
union of a mix of both if multiple keep-out zones are considered
\citep{hu2018trajectory}. As shall be seen below, polytope approximation methods
yield better optimality, while ellipsoidal methods yield better computational
efficiency. The distinction goes back to \ssref{scp,mip}, because polytope methods
often rely on MIP programming while ellipsoidal methods tend to use SCP.

For the case where \(\mathcal{B}(t)\) is a polytope,
\citep{schouwenaars2001plume,richards2002spacecraft} first proposed to write
\eqref{eq:rendezvous_e} as a set of mixed-integer constraints defined by the
polytope facets. The resulting trajectory optimization can be solved using MIP
methods discussed in
\ssref{mip}. \citep{richards2003model,richards2003performance,richards2006robust}
apply this approach in the context of MPC with a variable horizon trajectory.

For the case where \(\mathcal{B}(t)\) is an ellipsoid, \eqref{eq:rendezvous_e}
is typically enforced by checking for collision using a conservative time-varying
halfspace inclusion constraint:
\begin{equation}
  \label{eq:inclusion}
  \bm{r}(t)\in\set H(t)~\implies~\bm{r}(t)\notin\set B(t),
\end{equation}
where $\set H(t)$ is a halfspace. Three methods belonging to this family have
been used. The first is a rotating hyperplane method, proposed by
\citep{park2011model,diCairano2012model}. Here, \eqref{eq:rendezvous_e} is replaced by a
pre-determined sequence of halfspaces that are tangent to the ellipsoid and
rotate around it at a fixed rate. This approach was first applied to a 2D
mission, and later extended to 3D \citep{weiss2012model,weiss2015model}. A
variation was introduced in \citep{park2016analysis} and further studied in
\citep{zagaris2018model}, where the rotating sequence is replaced by just two
halfspaces tangent to the obstacle and passing through the chaser and target
positions. This method requires to pre-specifying which of the two halfspaces the
chaser belongs to at each time instant.

Fixing the halfspace sequence enables the first two approaches to retain
convexity. However, a third and most natural approach is to impose
\eqref{eq:inclusion} directly by linearizing the ellipsoidal obstacle. This
approach is taken in \citep{liu2014solving}, and has also been applied to
multiple moving obstacles \citep{jewison2015model,wang2018model}. Because
convexity is not maintained, SCP solution methods are used as discussed in
\ssref{scp}.  \citep{zagaris2018model} provide a detailed comparison of the three
methods.

Another challenge in RPO is the thrust constraint \eqref{eq:rendezvous_d}. This
constraint allows the thrust magnitude to take any value in the continuous
interval \([0, \rho]\). In reality, however, control is often realized by a
reaction control system (RCS) that produces short-duration pulses of constant
thrust. Therefore, in many applications it makes more sense to consider an \textit{impulse}
constraint of the form:
\begin{equation}
  \label{eq:impulse}
  \Delta v(t) \in \{0\}\union [\Delta v_{\min},\Delta v_{\max}],
\end{equation}
where $\Delta v(t)\in\reals$ approximates the instantaneous change in the
chaser's velocity following a firing from the RCS jets. Realistic RCS thrusters
have a minimum impulse-bit (MIB) performance metric that governs the magnitude
of the smallest possible velocity change $\Delta v_{\min}>0$. Because
\eqref{eq:impulse} is a nonconvex disjoint constraint of the form
\eqref{eq:example_mip}, it has been historically challenging to handle. Indeed,
\citep{larsson2006fuel} suggest that MIP is necessary in general, but in certain
cases the LP relaxation $\Delta v(t)\in [0,\Delta v_{\max}]$ of \eqref{eq:impulse}
suffices. This happens, for example, when the velocity measurement noise
exceeds the MIB value.

More recently, it was shown that the impulsive rendezvous problem can be solved
via polynomial optimization \citep{arzelier2011using,arzelier2013new}. Using
results on non-negative polynomials, \citep{deaconu2015designing} showed that
impulsive rendezvous with linear path constraints can be solved as an SDP. This
formulation was further embedded in a glideslope guidance framework for RPO
\citep{ariba2018minimum} and in an MPC approach
\citep{arantes2019stable}. Distinct from polynomial optimization,
\citep{malyuta2019lossless} proved that in some special cases the constraint
\eqref{eq:impulse} can be losslessly convexified using techniques similar to those
in \ssref{lcvx}. For problems where lossless convexification is not possible,
\citep{MalyutaAcikmeseJGCD2021} showed that SCP with a numerical continuation
scheme is an effective solution method. Yet another approach was presented in
\citep{wan2019alternating}, where an alternating minimization algorithm was
proposed for the case $\Delta v_{\min}=\Delta v_{\max}$, in other words when
the control is bang--bang.

The impulsive rendezvous model \eqref{eq:impulse} considers an instantaneous firing
duration. The model's accuracy can be improved by explicitly considering the
finite firing duration, leading to a representation of the actual pulse-width
modulated (PWM) thrust signal. PWM rendezvous was first studied in
\citep{vazquez2011trajectory,vazquez2014trajectory}, where an optimization
similar to \pref{rendezvous} was first solved, then the optimal continuous-valued
thrust signal was discretized using a PWM filter and iteratively improved using
linearized dynamics. This approach was later embedded in MPC
\citep{vazquez2015model,vazquez2017pulse}. A subtly different approach is
presented in \citep{li2016state,li2018pulse}, called pulse-width pulse-frequency
modulation (PWPF). Instead of iteratively refining the thrust signal, PWPF
passes the continuous-valued thrust signal to a Schmitt trigger that converts
it into a PWM signal. It was shown that this can save fuel and that stability
is maintained. However, a potential implementation disadvantage is that the
duration of each period in the resulting PWM signal varies continuously, which
conflicts with typical hardware where this period is an integer multiple of a
fixed duration. An SCP approach was recently used to account for this via
state-triggered constraints from \sssref{stc} \citep{malyuta2020fast}.

Although RPO literature tends to focus on the relative chaser-target position
using a 3-DoF model, relative attitude control also plays an important role,
especially if the target is tumbling \citep{li2017model,dong2020tube}. Thanks to
advances in the speed and reliability of optimization solvers as mentioned in
\ssref{ocp}, there has been an increasing interest to optimize 6-DoF RPO
trajectories with explicit consideration of position-attitude coupling through
constraints such as plume impingement and sensor pointing
\citep{ventura2017fast,zhou2019receding}. The resulting 6-DoF RPO trajectory
optimization, however, is much more challenging to solve due to the presence of
attitude kinematics and dynamics. Nevertheless, a special case with field of
view and glideslope constraints was presented in \citep{lee2014dual}, where
6-DoF RPO was solved as a convex quadratically constrained QP by using a dual
quaternion representation of the dynamics, effectively establishing a
convexification.

For more general RPO problems, nonlinear programming software has been used
frequently. For example, \citep{ventura2017fast} used SNOPT \citep{gill2005snopt}
after parameterizing the desired trajectory using polynomials. A B-spline
parameterization was used in \citep{sanchez2020flatness}, and the resulting
nonlinear optimization was solved by the IPOPT software
\citep{wachter2005implementation}. MATLAB-based packages were also used in
\citep{malladi2019nonlinear,volpe2019optical}. Recently, SCP techniques
discussed in \ssref{scp} were applied to 6-DoF RPO trajectory
optimization. \citep{zhou2019receding} considered both collision avoidance and
sensor pointing constraints. \citep{malyuta2020fast} further considered integer
constraints on the PWM pulse width in order to respect the RCS MIB value, and
constraints on plume impingement, by using state-triggered constraints. The
algorithm was improved in \citep{MalyutaAcikmeseJGCD2021} by making the solution
method faster and more robust. The approach uses homotopy ideas from
\sssref{indirect_mip} to blend the PTR sequential convex programming method with
numerical continuation into a single iterative solution process.

The operation of two spacecraft in close proximity naturally makes RPO a
safety-critical phase of any mission. Thus, trajectory optimization that is
robust to modeling errors, disturbances, and measurement noise has been an
active research topic. MPC has been a popular approach in this context, as it
allows efficiently re-solving \pref{rendezvous} with online updated parameters
using hardware with limited resources
\citep{hartley2014field,goodyear2015hardware,park2013model}.
\citep{hartley2015tutorial} provides a tutorial and a detailed discussion.
Among the many different approaches that have been developed to explicitly
address robustness, we may count feedback corrections \citep{baldwin2013robust},
the extended command governor \citep{petersen2014model}, worst-case analysis
\citep{louembet2015robust,xu2019collision}, stochastic trajectory optimization
\citep{jewison2018probabilistic}, chance constrained MPC
\citep{gavilan2012chance,zhu2018robust}, sampling-based MPC
\citep{mammarella2018offline}, tube-based MPC
\citep{mammarella2018tube,dong2020tube}, and reactive collision avoidance
\citep{ScharfReactiveCollision2006}. In addition to various uncertainties,
anomalous system behavior such as guidance system shutdowns, thruster
failures, and loss of sensing, also poses unique challenges in RPO. In order to
ensure safety in the presence of these anomalies,
\citep{luo2007optimal_a,luo2007optimal_b,luo2008multi} used a safety performance
index to discourage collision with the target, and \citep{breger2008safe}
considered both passive and active collision avoidance constraints in online
trajectory optimization. \citep{zhang2015optimal} considered passive safety
constraints together with field of view and impulse constraints. Aside from
optimization-based methods, artificial potential functions
\citep{dong2017safety,li2018potential,liu2019artificial} and sampling-based
methods \citep{starek2017fast} have also been applied to achieve safety in RPO.

\subsection{Small Body Landing}
\label{subsection:smallbody}

\begin{figure}
  \centering
  \includegraphics{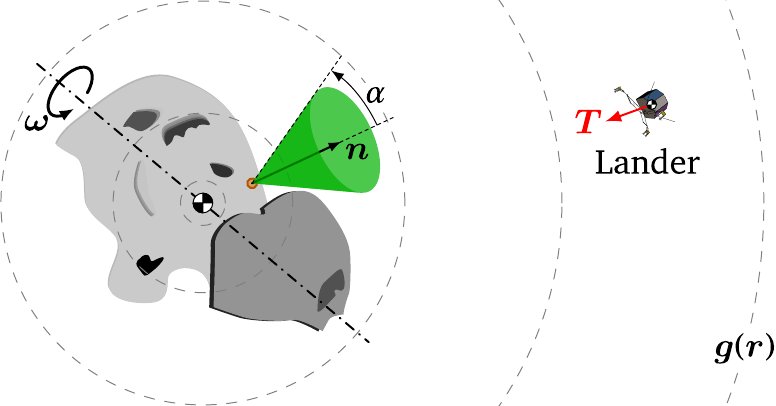}
  \caption{Illustration of a basic small body landing scenario. The basic
    concept is to use the thrust $\bm{T}$ to bring the lander spacecraft to a
    soft touchdown in the presence of rotational and gravitational
    nonlinearities, and operational constraints on glideslope, plume
    impingement, and collision avoidance.}
  \label{fig:smallbody}
\end{figure}

A maneuver similar to RPO is that of small body landing, where the target
spacecraft is replaced by a small celestial object such as an asteroid or a
comet. Trajectory optimization for small body landing has gathered increasing
levels of attention, spurred by recent high-profile asteroid exploration missions including Hayabusa
\citep{kawaguchi2008hayabusa}, Hayabusa2 \citep{crane2019final}, and OSIRIS-REx
\citep{berry2013osiris,lauretta2017osiris}. Unlike planetary rocket landing from
\ssref{pdg}, small body landing dynamics are highly nonlinear due to the irregular
shape, density, and rotation of the small body
\citep{werner1996exterior,scheeres1998dynamics}. Landing must furthermore ensure
a small touchdown velocity, possible plume impingement requirements, and
collision avoidance. These aspects pose unique challenges for trajectory
optimization. This section reviews recent developments in convex
optimization-based small body landing algorithms. Alternative trajectory
optimization methods have also been studied for this problem and which we do
not cover, such as indirect methods \citep{yang2015fuel,chen2019trajectory}.

The prototypical small body landing OCP is illustrated in \figref{smallbody} and
can be summarized as follows:
\begin{optimus}[
  task={\min},
  variables={t_f,\bm{T}},
  objective={\int_0^{t_f}\norm[2]{\bm{T}(t)}\dt},
  label={smallbody}
  ]
  \ddot{\bm{r}}(t)= -2\bm{\omega}\times
  \dot{\bm{r}}(t)-\bm{\omega}\times (\bm{\omega}\times
  \bm{r}(t))+m(t)\inv\bm{T}(t)+\bm{g}(\bm{r}(t)), \#
  \dot{m}(t)=-\alpha\norm[2]{\bm{T}(t)}, \#
  \rho_{\min}\leq \norm[2]{\bm{T}(t)}\leq \rho_{\max}, \#
    \norm[2]{\bm{r}(t)-\bm{r}_f}\cos\alpha\leq (\bm{r}(t)-\bm{r}_f)\T \bm{n},\#
  m(0)= m_0, ~ \bm{r}(0)=\bm{r}_0, ~ \dot{\bm{r}}(0)=\dot{\bm{r}}_0, ~
  \bm{r}({t_f})=\bm{r}_f, ~ \dot{\bm{r}}(t_f)=0.
\end{optimus}

Note the similarity between Problems~\ref{problem:pdg2007},
\ref{problem:rendezvous}, and \ref{problem:smallbody}. Compared to
\pref{pdg2007}, small body landing is expressed in the rotating frame of the
target. Thus, the main difference is in the dynamics \eqref{eq:smallbody_b}
that contain a general nonlinear gravity term \(\bm{g}(\bm{r}(t))\) and
inertial forces from the non-negligible angular velocity \(\bm{\omega}\) of the small body. The
glideslope constraint \eqref{eq:smallbody_e} is also shared with the approach
cone in RPO \eqref{eq:rendezvous_f}.

Early work by \citep{carson2006model,carson2008robust} ignored the mass dynamics
\eqref{eq:smallbody_c}, while \eqref{eq:smallbody_b} was linearized to solve for
acceleration rather than a thrust profile. The resulting tube MPC algorithm
includes a pre-determined feedback controller optimized using SDP and tracking
a feedforward trajectory from an SOCP in a robust and recursively feasible
manner. Some time later, \citep{pinson2015rapid} solved for a fixed-duration
trajectory by applying lossless convexification to \eqref{eq:smallbody_d} and
successive linearization to \eqref{eq:smallbody_b}, resulting in an SCP solution
method consisting of a sequence of SOCP
subproblems. \citep{pinson2018trajectory} further combine this solution
procedure with Brent's line search method to solve for the minimum-fuel flight
duration, which is similar to the use of golden-section search in the PDG context
\citep{blackmore2010minimum}. \citep{cui2017intelligent} combined
convexification with classic Runge-Kutta discretization to improve the
solution accuracy. \citep{yang2017rapid} showed how to solve the minimum-time
landing problem as a sequence of convex optimization problems. As a byproduct,
they showed that for time-optimal and short-duration minimum-landing-error
versions of \pref{smallbody}, the thrust stays at its maximum value, in which
case the lower bound in \eqref{eq:smallbody_d} can be removed and
\eqref{eq:smallbody_c} simplified.

Constraint \eqref{eq:smallbody_e} is the most basic type of collision avoidance
constraint. The heuristic reasoning behind \eqref{eq:smallbody_e} is that if the
lander stays approximately above a minimum glideslope, then it will avoid
nearby geologic hazards. An alternative two-phase trajectory optimization was
introduced in \citep{dunham2016constrained,liao2016model} by splitting the
landing maneuver into a circumnavigation and a landing phase. During
circumnavigation, the spacecraft is far away from the landing site and
\eqref{eq:smallbody_e} is replaced by collision avoidance constraint with the small
body. In the same manner as \ssref{rendezvous}, the small body is wrapped in an
ellipsoid and a rotating hyperplane constraint is used
\citep{dunham2016constrained,liao2016model,sanchez2018predictive}.
\citep{reynolds2017small} introduced an optimal separating hyperplane constraint
that also generates auxiliary setpoints for MPC tracking that converge to
the landing site. Once in close proximity to the landing site, the spacecraft
enters the landing phase where constraint \eqref{eq:smallbody_e} is enforced to
facilitate pinpoint landing.

Most small body landing work is 3-DoF in the sense that it considers point mass
translational dynamics. However, recently \citep{zhang2020trajectory} studied a
two-phase variable landing duration 6-DoF problem. The motivation was to impose a
field of view constraint for a landing camera. The resulting nonconvex
optimization trajectory problem was solved using SCP as covered in \ssref{scp}.

Parameters of the small body, such as $\bm{\omega}$ and $\bm{g}$, are often
subject to inevitable uncertainty, requiring judicious trajectory design. As a
result, many aforementioned works use MPC to cope with the uncertainty in small
body landing \citep{reynolds2017small,sanchez2018predictive}. Application
examples include tube MPC \citep{carson2006model,carson2008robust} and input
observers to compensate for gravity modeling errors
\citep{dunham2016constrained,liao2016model}. \citep{hu2016desensitized} also
proposed to jointly minimize fuel and trajectory dispersion described by
closed-loop linear covariance. For a detailed discussion on achieving
robustness in small body landing, we refer interested readers to the recent
survey \citep{simplicio2018review}.

\subsection{Constrained Reorientation}
\label{subsection:reorientation}

\begin{figure}
  \centering
  \includegraphics{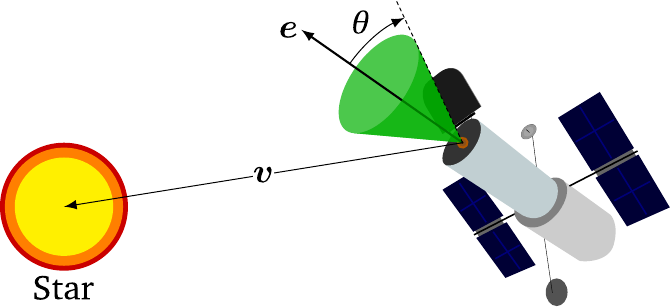}
  \caption{Illustration of a basic constrained reorientation scenario. The core
    challenge is to execute an attitude change maneuver while ensuring that the
    star vector $\bm{v}$ remains outside of the sensor keep-out cone.}
  \label{fig:reorientation}
\end{figure}

Scientific observation satellites commonly need to execute large angle
reorientation maneuvers while ensuring that their sensitive instruments, such
as cryogenically cooled infrared telescopes, are not exposed to direct sunlight
or heat. Famous examples include the Cassini spacecraft, the Hubble Space Telescope, and the
upcoming James Webb Space Telescope
\citep{singh1997constraint,JamesWebbConOps,HubbleConOps}. This section discusses
the challenges of constrained reorientation as a trajectory optimization
problem, and focuses on how convex optimization methods have been leveraged to
address these challenges.

A basic constrained reorientation OCP is illustrated in \figref{reorientation}
and can be formulated as follows:
\begin{optimus}[
  task={\min},
  variables={t_f, \bm{u}},
  objective={\int_0^{t_f}\norm[2]{\bm{u}(t)}\dt},
  label={reorientation}
  ]
  \dot{\bm{q}}(t)=\frac{1}{2} \bm{q}(t)\otimes\bm{\omega}(t), \#%~ \norm[2]{\bm{q}(t)}=1\#
  J\dot{\bm{\omega}}(t)=\bm{u}(t)- \bm{\omega}(t)\times(J\bm{\omega}(t)), \#
  \bm{q}(t)\T M_i\bm{q}(t)\leq 0,~ i=1, \ldots, n,\#
  \norm[\infty]{\bm{\omega}(t)}\leq \omega_{\max}, ~ \norm[\infty]{\bm{u}(t)}\leq u_{\max},\#
  \bm{q}(0)=\bm{q}_0, ~ \bm{\omega}(0)=\bm{\omega}_0,~ \bm{q}(t_f)= \bm{q}_f, ~ \bm{\omega}(t_f)= \bm{\omega}_f.
\end{optimus}

The set of constraints \eqref{eq:reorientation_d} encodes conical keep-out zones
for $n$ stars, similarly to the illustration in \figref{reorientation} for one
star. The parameters $M_i\in\reals^{4\times 4}$ are symmetric matrices that are
not positive semidefinite, as introduced in \ssref{lcvx}. The main challenge of
solving \pref{reorientation} stems from the fact that \eqref{eq:reorientation_d} and
the attitude dynamics \eqref{eq:reorientation_b}-\eqref{eq:reorientation_c} are
nonconvex.  \citep{kim2004quadratically} were the first to prove that
\eqref{eq:reorientation_d} can be losslessly replaced by convex quadratic
constraints, provided $\norm[2]{\bm{q}(t)}=1$. Based on this observation,
\citep{kim2004quadratically} proposed to greedily optimize one discretization
point at a time instead of the entire trajectory jointly.  The method was
further extended to the case of integral and dynamic pointing constraints in
\citep{kim2010convex}.

Although the method of \citep{kim2004quadratically} is computationally
efficient, it is inherently conservative and may fail to find a feasible
solution to \pref{reorientation} by greedily optimizing one discretization point at
a time.  As a result, several attempts have been made to improve its
performance. For example, \citep{tam2016constrained} propose to replace
constraint \eqref{eq:reorientation_d} with penalty terms in the objective function
in order to ensure that a feasible trajectory can be found. Binary logical variables were
also introduced in \eqref{eq:reorientation_d} to account for redundant
sensors. \citep{hutao2011rhc} showed how the convexification of constraints
\eqref{eq:reorientation_d} should be adjusted when optimizing an entire trajectory,
rather than a single time step as originally done in
\citep{kim2004quadratically}. Put into an MPC framework, the resulting
trajectory optimization yields less conservative performance. Alternatively,
\citep{eren2015mixed} proposed to first optimize a quaternion sequence without
kinematic and dynamic constraints, and then to compute the corresponding
angular velocity and torque using the resulting quaternions. A hyperplane
approximation of the unit sphere is used during quaternion optimization to
ensure dynamic feasibility, and is imposed via MIP. {Recently,
  \citep{McDonald2020} proposed an SCP method with a line search step that helps
  convergence, which provides a potential real-time solution to
  \pref{reorientation}}.

Aside from the quaternion representation in \pref{reorientation}, which is the
most popular choice, a
direction cosine matrix representation of attitude was also used by
\citep{walsh2018constrained} to solve an equivalent problem. The resulting
trajectory optimization can be approximated as an SDP using successive
linearization and relaxing \eqref{eq:reorientation_d}.

Due to its challenging nature, \pref{reorientation} has inspired many optimization
solutions other than those based on convex optimization. Pseudospectral methods
and NLP optimization software have all been used to solve \pref{reorientation}
directly \citep{xiaojun2010large,lee2013quaternion}. An indirect shooting method
was used in \citep{lee2017geometric,Phogat2018discreteJGCD}, and a differential
evolution method was used in \citep{wu2017time,wu2019energy}. Compared with
convex optimization based methods, these methods typically require more
computational resources to achieve real-time implementation.

\subsection{Endo-atmospheric Flight}
\label{subsection:endo}

\newcommand{\mtx}{\textcolor{magenta}}

Launching from or returning to a planet with an atmosphere are integral parts
of many space missions. These problems concern launch vehicles, missiles, and
entry vehicles such as capsules, reusable launchers, and hypersonic
gliders. Significant portions of launch and entry occur at high velocities and
in the presence of an atmosphere, making aerodynamics play a large
role. Aerodynamics and thermal heating are indeed the core differentiating
factors between endo-atmospheric flight and PDG from \ssref{pdg}. For the latter problem, small
velocities and thinness of the atmosphere make aerodynamic effects negligible
in many cases \citep{eren2015constrained}. This section summarizes recent
contributions to endo-atmospheric trajectory planning using convex
optimization-based methods. In particular, \sssref{ascent} discusses ascent and
\sssref{reentry} discusses entry.

\subsubsection{Ascent Flight}
\label{subsubsection:ascent}

\begin{figure}
  \centering
  \includegraphics{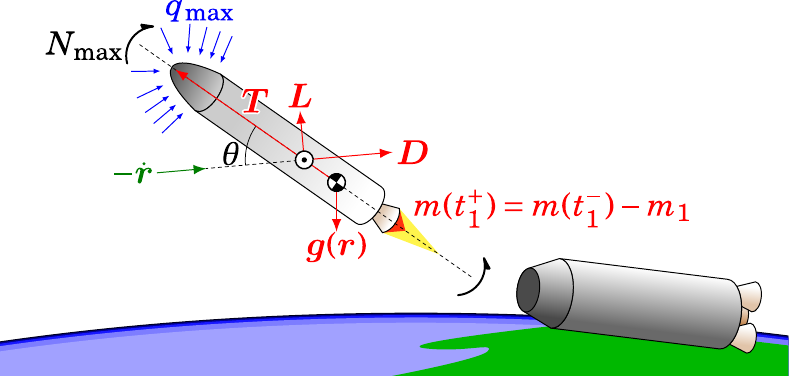}
  \caption{Illustration of a basic ascent scenario. The goal is to find an
    optimal angle-of-attack $\theta$ trajectory to transfer the launch vehicle's payload
    from the planet's surface to orbit, while minimizing fuel and satisfying
    structural integrity constraints.}
  \label{fig:ascent}
\end{figure}

The optimal ascent problem seeks to transfer a launch vehicle's payload from a
planet's surface to orbit while minimizing a quantity such as fuel. Naturally,
optimal control theory from \ssref{ocp} has found frequent applications in ascent
guidance, and we refer the reader to \citep{hanson1994ascent} for a
survey. Heritage algorithms date back to the iterative guidance mode (IGM) of
Saturn rockets
\citep{chandler1967development,horn1969iterative,adkins1970optimization,haeussermann1971saturn}
and the powered explicit guidance (PEG) of the Space Shuttle
\citep{mchenry1979space}.  A simple yet relevant optimal control problem
describing an orbital launch scenario is known as the \textit{Goddard rocket
  problem} \citep{betts2010practical, BrysonBook}. A version with variable
gravity and no atmospheric drag is stated as follows:
\begin{optimus}[
  task=\min,
  variables={t_f,\bm{T}},
  objective={-m(t_f)},
  label={goddard}
  ]
  \ddot{\bm{r}}(t) = -\mu\norm[2]{\bm{r}(t)}^{-3}\bm{r}(t)+
  m(t)\inv\bm{T}(t), \#
  \dot m(t) = -\alpha\norm[2]{\bm{T}(t)}, \#
  m(0)=m_0,~\bm{r}(0)=\bm{r}_0,~\dot{\bm{r}}(0)=\dot{\bm{r}}_0,~%
  \bm{\psi}(\bm{r}(t_f),\dot{\bm{r}}(t_f))=0.
\end{optimus}

\pref{goddard} models a three-dimensional point mass moving in a gravity field
under the influence of thrust. As such, it also applies to orbit transfer
problems which we discuss later in \ssref{orbit}. The vector
$\bm{r}(t)\in\real^3$ is the position vector, $\bm{T}(t)\in\reals^3$ is the
thrust vector, and $m(t)\in\reals$ is the vehicle mass. The vector function
$\bm{\psi}:\reals^6\to\reals^k$ imposes $k\le 6$ terminal conditions. In ascent
and orbit transfer applications, $\bm{\psi}$ usually acts to constrain $k$
orbital elements while leaving the other $6-k$ orbital elements free.

A key issue when solving \pref{goddard} using an indirect method is to resolve the
transversality conditions of the resulting TPBVP \citep{PontryaginBook,
  berkovitz1974optimal}:
\begin{equation}
  \label{eq:transversality_tpbvp}
  \bm{p}(t_f) = \big[\grad_{\bm{x}}\bm{\psi}\bigl(\bm{x}(t_f)\bigr)\big]\T\bm{\nu}_p,
\end{equation}
where $\bm{x}(t)\definedas(\bm{r}(t);\dot{\bm{r}}(t))$,
$\bm{p}(\cdot)\in\reals^6$ are the costates relating to the position and
velocity, and $\bm{\nu}_p\in\reals^k$ is a Lagrange multiplier
vector. Unfortunately, $\bm{\nu}_p$ has no physical or exploitable numerical
interpretation, and the magnitudes of its elements can vary wildly
\citep{pan2013reduced}. This causes a lot of difficulty for the solution process
in terms of numerics, robustness, and initial guess selection. Traditionally,
the problem is overcome by converting \eqref{eq:transversality_tpbvp} into a set of
$6-k$ so-called reduced transversality conditions, which are equivalent
\citep{lu2003closed}:
\begin{subequations}
  \begin{alignat}{3}
    \big[\grad_{\bm{x}}\bm{\psi}\bigl(\bm{x}(t_f)\bigr)\big]\bm{y}_i
    &= 0,\quad &&i=1,\dots,6-k, \\
    \label{eq:reduced_transversality}
    \bm{y}_i\T\bm{p}(t_f) &= 0,\quad &&i=1,\dots,6-k.
  \end{alignat}
\end{subequations}

The linearly independent vectors $\bm{y}_i\in\reals^6$ are known as the reduced
transversality vectors, and are a function of $\bm{x}(t_f)$. If they are known
analytically, then \eqref{eq:reduced_transversality} can replace
\eqref{eq:transversality_tpbvp}, which eliminates $\bm{\nu}_p$ from the problem and
simplifies the solution process considerably. However, solving for $\bm{y}_i$
symbolically is a difficult task, and the resulting expressions can be
complicated \citep{brown1967real}. An alternative approach was introduced in
\citep{pan2013reduced} where the authors provide an easy to use ``menu'' of the $6-k$
constraints in \eqref{eq:reduced_transversality} that are needed. This is achieved
by considering \pref{goddard} specifically and exploiting the structure offered by
the classical orbital elements. It is only assumed that the terminal constraint
function $\bm{\psi}$ fixes exactly $k$ of the 6 orbital elements, and leaves
the other orbital elements free.

The Goddard rocket problem in \pref{goddard} assumes no atmosphere.  When
there is an atmosphere, a popular classical method is the \textit{gravity turn}
maneuver, which maintains a low angle-of-attack so as to minimize lateral
aerodynamic loads. However, the general ascent problem with an atmosphere is
complicated due to strong coupling of aerodynamic and thrust forces
\citep{pan2010improvements}. Thus, ascent is typically performed via open-loop
implicit guidance, in the sense that feedback control is used to track a
pre-computed ascent trajectory stored onboard as a lookup table. However, this
approach cannot robustly handle off-nominal conditions, aborts, and
contingencies, which motivates research into closed-loop ascent techniques
\citep{brown1967real,lu2017cgc}.

A notable strategy in this context is to include aerodynamics in an onboard
ascent solution via a homotopy method starting from an optimal vacuum ascent
initial guess
\citep{calise1998design,gath2001optimization,calise2004generation}. Another
approach was developed in
\citep{lu2003closed,lu2005ascent,lu2010highly,pan2010improvements} using
indirect trajectory optimization. Here, a finite-difference scheme is proposed
to solve for the necessary conditions of optimality for ascent with an
atmosphere. In particular, fixed-point formulations were considered
\citep{lu2005ascent}, primer vector theory was invoked to determine trajectory
optimality \citep{lu2010highly}, and a generalization to arbitrary numbers of
burn and coast arcs was developed \citep{pan2010improvements}. {Finally,
  indirect methods relying on control smoothing via trigonometrization have
  been developed to address problems with bang--bang input and singular arcs
  \citep{mall2020uniform}. The Epsilon-Trig method \citep{mall2017epsilon}, which
  is an example of such an approach, was applied to the Goddard maximum
  altitude ascent problem to obtain its bang--singular--bang optimal
  solution. See \sssref{indirect} for a brief description of these approaches.}

Modern improvements in convex optimization have made direct optimization methods
attractive for ascent guidance. To this end, consider the following illustrative
ascent problem for a two-stage launch vehicle, as shown in \figref{ascent}:
\begin{optimus}[
  task=\min,
  variables={t_f,\theta},
  objective={-m(t_f)},
  label={ascent}
  ]
  \ddot{\bm{r}}(t) = -\mu\norm[2]{\bm{r}(t)}\inv[3]\bm{r}(t)+
  m(t)\inv(\bm{T}[t]+\bm{L}[t]+\bm{D}[t]), \#
  \dot m(t) = -\alpha\norm[2]{\bm{T}[t]}, \#
  \theta_{1}\le\theta(t)\le\theta_{2}, \#
  \rho[t]\norm[2]{\dot{\bm{r}}(t)}^2\leq q_{\max}, \#
  \rho[t]\norm[2]{\dot{\bm{r}}(t)}^2|\theta(t)| \leq N_{\max}, \#
  m(0)=m_0,~\bm{r}(0)=\bm{r}_0,~\dot{\bm{r}}(0)=\dot{\bm{r}}_0,~%
  \bm{\psi}(\bm{r}(t_f),\dot{\bm{r}}(t_f))=0, \#
  m(t_1^+)=m(t_1^-)-m_1.
\end{optimus}

\pref{ascent} is planar and formulated in an Earth-centered inertial (ECI)
frame. Control is performed using the angle-of-attack $\theta$, which
determines the direction of an otherwise pre-determined thrust profile
\citep{zhang2019rapid,li2020online,liu2014solving}. The major aerodynamic forces
are those of lift $\bm{L}$ and drag $\bm{D}$, each of which may be complex
expressions of state and control. Note that in \eqref{eq:ascent_b} we used the
shorthand $\bm{T}[t]$, $\bm{L}[t]$, and $\bm{D}[t]$ from \sref{introduction} to
abstract away the possible state and control arguments. The atmospheric density
is denoted by $\rho$, which varies during ascent as a nonlinear function of the
position $\bm{r}$. An example is given later in \eqref{eq:density}.  Important
constraints on the dynamic pressure \eqref{eq:ascent_e} and bending moment
\eqref{eq:ascent_f} are used to ensure the vehicle's structural integrity
\citep{lu2010highly}. The target orbit is prescribed by the vector function
$\bm{\psi}$ in \eqref{eq:ascent_g}, which is the same as in \eqref{eq:goddard_d} and
specifies some or all of the target orbital
elements. \citep{benedikter2020convex} chose boundary conditions based on the
radius and inclination of a circular target orbit. A final nuance is that, if
the rocket is assumed to be a two-stage vehicle, a stage separation event must
be scheduled at a pre-determined time $t_1$ via \eqref{eq:ascent_h}. At the
separation instant, the mass variable experiences a discontinuous decrease that
amounts to the {dry weight of the first stage
  \citep{benedikter2019convex_b,benedikter2020convex}}. A related constraint for
stage separation requires $\theta(t_1)=0$ in order to reduce lateral load
\citep{zhengxiang2018convex}. Furthermore, the splashdown location of a
burnt-out separated stage can also be constrained \citep{benedikter2020convex}.

Due to the presence of strong nonlinearities, convex optimization-based
solution algorithms for \pref{ascent} typically use SCP from \ssref{scp}. However, several
manipulations have been helpful to make the problem less nonlinear. Conversion
of the system dynamics \eqref{eq:ascent_b} to control-affine form, at times by
choosing an independent variable other than time, followed by the use of
lossless convexification within an SCP framework has been a common
approach. \citep{zhang2019rapid} obtained a control-affine form by assuming the
AoA to be small and defining $u_1=\theta$, $u_2=\theta^2$ as the new control
variables. This choice makes drag a linear function of the control, while the
constraint $u_1^2=u_2$ is relaxed to $u_1^2\le u_2$ via lossless
convexification. {Similarly, \citep{benedikter2019convex_b,benedikter2020convex}
  chose} thrust direction as input and losslessly convexified the unit norm
constraint on the thrust direction to a convex
inequality. \citep{cheng2017efficient} considered a 3D problem with AoA and bank
angle as control inputs, and applied lossless convexification to a constraint
of the form $u_1^2 + u_2^2 + u_3^2=1$. Furthermore, their choice of altitude as
the independent variable simplified the convexification of constraints
involving density, since it renders the density a state-independent
quantity. In particular, during collocation over a known grid within an
altitude interval, the density value is known at each node. This fails to be
the case when collocation is performed over time. The choice of altitude as
independent variable was also explored in \citep{liu2016exact}.

 The two-agent launch problem is an interesting and relevant
modern-day extension of \pref{ascent} \citep{ross2004rapid}. In this case, the
launch vehicle first stage is not just an idle dropped mass, but is a
controlled vehicle that must be brought back to Earth. This is the case for the
SpaceX Falcon 9 rocket, whose first stage is recovered by propulsive landing
after a series of post-separation maneuvers \citep{blackmore2016autonomous}. It
was shown in \citep{ross2005hybrid} how hybrid optimal control can be used to
solve the problem via mixed-integer programming. More generally, hybrid
optimal control has also found applications in low-thrust orbit transfer using
solar sails \citep{stevens2005preliminary,stevens2004earth}.

\subsubsection{Atmospheric Entry}
\label{subsubsection:reentry}

%\arccomment{@Purnanand}{Please find a good place to include citations\citep{lu2010rapid,johnson2020pterodactyl}}

\begin{figure}
  \centering
  \includegraphics{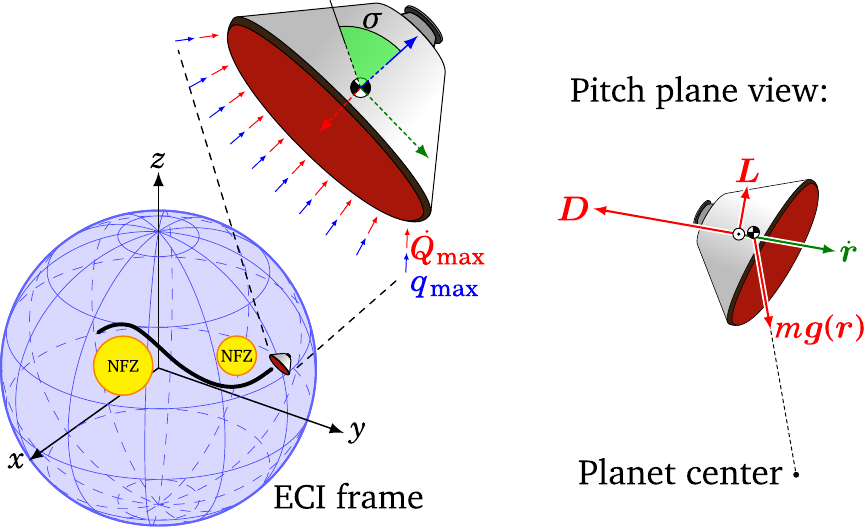}
  \caption{Illustration of a basic atmospheric entry scenario. The goal is to
    find an optimal bank angle $\sigma$ trajectory to dissipate energy while
    meeting structural integrity and targeting requirements.}
  \label{fig:reentry}
\end{figure}

Atmospheric entry, also known as reentry, is fundamentally a process of
controlled energy dissipation while meeting targeting and structural integrity
constraints
\citep{lu2014entry}. % Direct, skip, and lofted are prime examples of entry
% trajectories.
% Aerocapture is
% another related maneuver, where a pass through the atmosphere is used to insert
% a spacecraft into a parking orbit around a
% planet.
Computer-controlled entry guidance dates back to the Gemini and Apollo projects,
and \citep{sarigul2014survey} provide a comprehensive survey of existing
methods. Good documentation is available for Mars Science Laboratory's entry
guidance, which is based on Apollo heritage
\citep{way2007mars,mendeck2011entry,steltzner2014mars}.

A large body of work is available on predictor-corrector methods for entry
guidance
\citep{xue2010constrained,johnson2020pterodactyl,johnson2018entry,johnson2017automated,lu2014entry}
and for aerocapture \citep{lu2015optimal}. These methods are based on
root-finding algorithms, or variations thereof, and some versions are grounded
in solving the necessary conditions of optimality \citep{lu2018propellant}. We
refer the reader to \citep{lu2008predictor} for further details. In addition to
reentry trajectory generation, mission analysis tools for generating landing
footprints have also been developed \citep{lu2010rapid,eren2015constrained}.

Guidance schemes based on univariate root-finding, which are near-optimal for
reentry \citep{lu2014entry} and optimal for aerocapture \citep{lu2015optimal},
have also been developed. Reentry applications exploit the quasi-equilibrium
glide condition (QEGC), while aerocapture leverages the bang--bang nature of the
control solution obtained via the maximum principle. By using the knowledge of
where the input switches, univariate root-finding can approximate the optimal
solution in each phase to high accuracy. Such an approach, though based on an
indirect method, avoids directly solving the TPBVP. Recall that lossless
convexification, discussed in \ssref{lcvx}, is another approach where clever
reformulation of the optimal control problem and application of the maximum
principle yields an efficient solution strategy. This ties back to the last
paragraph of \sssref{direct}, which states that the fusion of indirect and direct
solution methods often yields more efficient solution algorithms than using any
one method in isolation. Because root-finding algorithms do not involve an
explicit call to an optimizer, we do not survey them here. Instead, this
section focuses on contributions by convex optimization-based methods to the
problem of entry trajectory computation.

Another methodology that simplifies the typical strategy in indirect methods
is the RASHS approach \citep{saranathan2018relaxed}. As discussed in
\sssref{indirect_mip}, RASHS converts a multi-phase optimal control problem into a
single-phase problem by using sigmoid functions of state-based conditions to
instigate smooth transitions between phases. As a consequence, the multi-point
BVP from Pontryagin's maximum principle is reduced to a TPBVP. The complete
entry, descent, and landing (EDL) problem is one example that can be solved
effectively via this technique.

Consider a basic entry guidance problem illustrated in \figref{reentry}, which
is formulated as follows:
\begin{optimus}[
  task=\min,
  variables={u},
  objective={~\max_{t\in [0,t_f]}\dot Q[t]},
  label={entry}
  ]
  \ddot{\bm{r}}(t) = -\mu\norm[2]{\bm{r}(t)}\inv[3]\bm{r}(t)+
  m\inv\bigl(\bm{L}[t]+\bm{D}[t]\bigr),
  \#
  |u(t)|\le 1, \#
  \rho[t]\norm[2]{\dot{\bm{r}}(t)}^2\le q_{\max}, \#
  \norm[2]{\bm{L}[t]+\bm{D}[t]}\le n_{\max}, \#
  \bm{r}(0) = \bm{r}_0,~\dot{\bm{r}}(0)=\dot{\bm{r}}_0.
\end{optimus}

\pref{entry} is planar and formulated in an ECI frame like \pref{ascent}. Aerodynamic
forces are governed by the lift, drag, and atmospheric density, which are
expressed as follows:
\begin{subequations}
  \begin{align}
    \bm{L}[t] &= R_{\pi/2}c_{L}\rho[t]\norm[2]{\dot{\bm{r}}(t)}\dot{\bm{r}}(t)u(t), \\
    \bm{D}[t] &= R_{\pi}c_{D}\rho[t]\norm[2]{\dot{\bm{r}}(t)}\dot{\bm{r}}(t), \\
    \rho[t] &= \rho_0\exp[inline]{-(\norm[2]{\bm{r}(t)}-r_0)/h_0}. \label{eq:density}
  \end{align}
\end{subequations}

The lift and drag coefficients are given by $c_L$ and $c_D$ while $\rho_0$,
$h_0$, and $r_0$ denote the reference density, reference altitude, and planet
radius. $R_\theta$ corresponds to a counter-clockwise rotation by $\theta$
radians. The control input $u(t)=\cos(\sigma(t))$ is the cosine of the bank
angle, and serves to modulate the projection of the lift vector onto the plane
of descent, known also as the \textit{pitch plane}. Entry is an extremely
stressful maneuver for much of the spacecraft's hardware, therefore structural
integrity constraints are placed on dynamic pressure \eqref{eq:entry_d} and
aerodynamic load \eqref{eq:entry_e}. The objective is to minimize the peak heating
rate, given by the Detra--Kemp--Riddell stagnation point heating correlation
\citep{detra1957heat,garrett1970heat}, which is appropriate for an insulative
reusable thermal protection system (TPS) such as on the Space Shuttle and
SpaceX Starship::% \citep{wang2018pseudospectral}:
\begin{equation}
  \label{eq:heating_rate}
  \dot Q[t]\propto \sqrt{\rho[t]}\norm[2]{\dot{\bm{r}}(t)}^{3.15}.
\end{equation}

\pref{entry} is a simple example that gives a taste for the reentry problem. We now
survey variants of this problem that have been explored in the
literature. First of all, many other objectives have been proposed in place of
\eqref{eq:entry_a}. These include a minimum heat load
\citep{wang2017constrained,han2019rapid}, minimum peak normal load
\citep{wang2019maximum,wang2019optimal-normal}, minimum time-of-flight
\citep{wang2019rapid,han2019rapid}, minimum terminal velocity
\citep{wang2017constrained}, maximum terminal velocity \citep{wang2019improved},
minimum phugoid oscillation \citep{liu2016rapid}, and minimum cross-range error
\citep{fahroo2003footprint,fahroo2003modeling}. In the problem of aerocapture,
where a spacecraft uses the planet's atmosphere for insertion into a parking
orbit, minimum velocity error \citep{zhang2015convex} and minimum impulse,
time-of-flight, or heat load \citep{han2019rapid} were studied. Minimizing the
total heat load, which is equivalent to the average heating rate, is
particularly relevant for ablative TPS that work by carrying heat away from the
surface through mass loss. This has been the method of choice for Apollo,
SpaceX Crew Dragon, and almost all interplanetary entry systems, because it can
sustain very high transient peak heating rates \citep{hicks2009entry}.

\pref{entry} is expressed in the pitch plane and without regard for planetary
rotation. To account for rotation and aspects like cross-range tracking, other
formulations have been explored. This includes a pitch plane formulation with
rotation \citep{chawla2010suboptimal}, a 3D formulation without rotation
\citep{zhao2017reentry}, and a 3D formulation with rotation
\citep{wang2017constrained,wang2019optimal,wang2019maximum,wang2018autonomous,
  wang2019rapid,liu2016rapid,liu2016entry,liu2015solving,han2019rapid}.

The two main path constraints present in \pref{entry} are on the dynamic pressure
\eqref{eq:entry_d} and aerodynamic load \eqref{eq:entry_e}. The heating rate is also
indirectly constrained since \eqref{eq:entry_a} must achieve a lower value than the
maximum heating rate $\dot Q_{\max}$, otherwise the computed trajectory melts
the spacecraft. Since these three constraints are critical for structural
integrity, they permeate much of reentry optimization literature
\citep{wang2017constrained,wang2019improved,wang2019optimal,wang2018autonomous,
  wang2019maximum,zhao2017reentry,wang2019rapid,liu2015solving,liu2016entry,
  han2019rapid,liu2016rapid,sagliano2018optimal}. Some researchers have also
included no-fly zone (NFZ) constraints, as illustrated in \figref{reentry}
\citep{zhao2017reentry,liu2016entry}. A bank angle reversal constraint has also
been considered in
\citep{zhao2017reentry,han2019rapid,liu2015solving,liu2016entry,liu2016rapid}. This
is a nonconvex constraint of the form:
\begin{equation}
  \label{eq:bank_angle_constraint}
  0<\sigma_{\min}\le |\sigma(t)| \le\sigma_{\max}.
\end{equation}

A common approach to handle \eqref{eq:bank_angle_constraint} is to define
$u_1(t)\definedas\cos(\sigma(t))$ and $u_2(t)\definedas\sin(\sigma(t))$, and to
impose:
\begin{equation}
  \label{eq:bank_angle_rewrite}
  \cos(\sigma_{\max})\le u_1(t)\le\cos(\sigma_{\min}),~
  u_1(t)^2+u_2(t)^2=1,
\end{equation}
where the nonconvex equality constraint is subsequently losslessly
convexified to $u_1(t)^2+u_2(t)^2\le 1$
\citep{liu2016rapid,liu2016entry,liu2015solving}.

The bank angle with a prescribed AoA profile is a popular control input choice
for reentry, dating back to Apollo \citep{rea2007comparison}. Some works have
considered bank angle rate as the input, which improves control smoothness
\citep{wang2017constrained,wang2019improved,wang2019maximum,
  wang2019optimal,wang2018autonomous,wang2019rapid}. However, banking is not the
only possible control mechanism for reentry, and several other choices have
been explored. \citep{chawla2010suboptimal} use the AoA as input and omit bank
and heading. \citep{fahroo2003footprint} use AoA, bank angle, and altitude,
assuming the aforementioned QEGC with a small flight-path angle between the
velocity vector and the local horizontal. \citep{zhao2017reentry} use bank angle
and a normalized lift coefficient as inputs.

High frequency oscillation in the control signal, known as jitter, is a common
issue in entry trajectory optimization. Several works explicitly
address this issue \citep{szmuk2017successive,wang2017constrained,
  liu2016entry,liu2015solving}. Jitter is believed to be caused by the nonlinear
coupling of state and control constraints \citep{wang2017constrained}, and it
appears to be reduced by a control-affine reformulation of the dynamics
\citep{liu2016entry}. Other strategies to remove jitter have been to apply the
reparametrization \eqref{eq:bank_angle_rewrite} or to filter the control
signal. The latter approach includes the aforementioned use of bank angle rate
as the control, or using a first-order low-pass filter \citep{liu2015solving}.

Aside from fixing jitter, efforts have been devoted to simplifying the
SCP-based solution methods, and to improving their convergence
properties. Reformulating the dynamics using energy as the independent
variable, in a similar way to how altitude was used for optimal ascent, is one
tactic that achieves the former
\citep{lu2014entry,liu2015solving,liu2016rapid,liu2016entry}.  Such a
parametrization eliminates the differential equation for airspeed, and instead
yields an algebraic approximation for airspeed in terms of
energy. \citep{fahroo2003footprint,fahroo2003modeling} applied a related
elimination process by considering energy as state variable. Apart from this,
it is worth noting the heuristics proposed for improving the convergence of the
SCP-based
approaches. \citep{liu2015solving,wang2019improved,wang2019optimal-normal} used
backtracking line search at each SCP iteration to reduce constraint
violation. It was found that with the line search, the number of iterations
required for convergence reduced by half \citep{liu2015solving}.
\citep{zhang2015convex} constrained the SCP iterates to form a Cauchy
sequence. \citep{wang2019rapid} proposed a dynamic trust region update scheme
that is tailored for hypersonic reentry. In particular, the trust region update
accounts for the linearization error due to each state instead of the typical
approach of considering the average linearization error.
% \arccomment{@Purnanand}{For the last three sentences, can you provide
% additionally what the authors have found -did convergence improve? By how
% much and with respect to what baseline?}

Aside from using SCP to optimize the entire entry trajectory, another popular
approach for entry guidance is via MPC from \ssref{mpc}. Some important recent
MPC-based developments are the dynamic control allocation scheme
\citep{luo2007model} and the application of model predictive static programming
(MPSP). The approach by \citep{luo2007model} centers around posing an SQP
problem as a linear complimentarity problem, while the principle behind MPSP is
to combine MPC and approximate dynamic programming though a parametric
optimization formulation
\citep{halbe2014robust,halbe2010energy,chawla2010suboptimal}.
\citep{van2006combined,recasens2005robust} corroborate the effectiveness of
MPC-based approaches by comparing the performance of constrained MPC with that
of PID control applied to feedback-linearized reentry flight.

Last but not least, we conclude by discussing pseudosectral discretization from
\sssref{pseudospectral} as a popular methodology in a variety of reentry problem
formulations. The method is appealing for its ability to yield accurate
solutions with a relatively sparse temporal collocation grid, and recent
results on the estimation of costates with spectral accuracy provide a strong
theoretical grounding \citep{francolin2014costate,gong2010costate}.
% \arccomment{@Purnanand}{Can you add a sentence on \textit{why}people looked
% at pseudospectral methods? Is it only for higher accuracy?}.
\citep{rea2003launch,fahroo2003footprint,fahroo2003modeling} applied direct
Legendre collocation to generate an entry vehicle footprint by solving a
nonconvex NLP. \citep{sagliano2018optimal} used Legendre-Lobatto collocation and
lossless convexification to generate an optimal profile for the heritage
drag-energy guidance scheme. In addition to these approaches, which rely solely
on the direct method, a combination of direct and indirect methods was
discussed in \citep{josselyn2002rapid} for verifying optimality of reentry
trajectories using the DIDO solver (see \tabref{software}). In
\citep{tian2011optimal}, a feedback guidance law through an indirect Legendre
pseudospectral method was developed to track a reference generated using a
direct pseudospectral method. Finally, akin to explicit MPC,
\citep{sagliano2016onboard} developed a pre-computed interpolation-based
multivariate pseudospectral technique that is coupled with a subspace selection
algorithm to generate nearly optimal trajectories in real-time for entry
scenarios in the presence of wide dispersions at the entry interface.

\subsection{Orbit Transfer and Injection}
\label{subsection:orbit}

\begin{figure}
  \centering
  \includegraphics{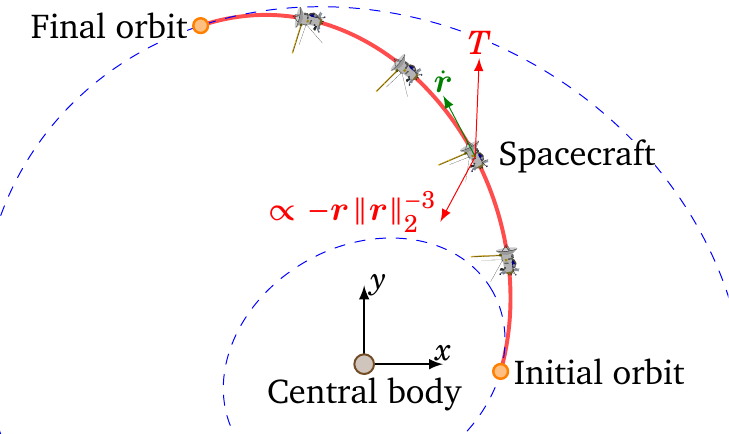}
  \caption{Illustration of a basic orbit transfer scenario. The goal is to use
    thrust $\bm{T}$ to transfer the spacecraft state from an initial orbit to a
    target orbit under the influence of gravity, while minimizing a quantity
    such as fuel or time.}
  \label{fig:orbit}
\end{figure}

A usual task in a space mission is to attain a certain orbit, or to change
orbits. The goal of the so-called orbit transfer and injection (OTI) problem is
to transport a low-thrust space vehicle from an initial to a target orbit while
minimizing a quantity such as time or fuel. Unsurprisingly, the problem is as
old as spaceflight itself, with the earliest bibliographic entry dating to the
late 1950s \citep{faulders1958low}. Traditionally, the problem has been solved
using optimal control theory from \ssref{ocp}, and for this we can cite the books
\citep{LonguskiBook,LawdenBook,BrysonBook,KirkBook,ConwayBook}. Numerous
solution methods have been studied, {including methods based on primer vector
  theory \citep{russell2007primer,petropoulos2008low,restrepo2017shadow}, direct
  methods based on solving an NLP
  \citep{betts2000very,arrieta2005real,ross2007low,starek2014nonlinear,
    graham2015minimum,graham2016minimum}, and indirect methods
  \citep{alfano1994circle,fernandes1995optimum,kechichian1995optimal,
    haberkorn2004low,gong2008spectral,gil2010practical,zimmer2010reducing,pan2012coast,
    pontani2014optimal,cerf2016low,taheri2016enhanced,taheri2017costate,rasotto2015multi,DiLizia2014high}. Some
  recent advances for indirect methods include homotopy methods {
    \citep{Pan2019,Pan2020,cerf2011continuation}}, optimal switching surfaces
  \citep{Taheri2019}, the RASHS and CSC approaches from \sssref{indirect_mip}
  \citep{saranathan2018relaxed,taheri2020novel,Taheri2020b}}, {and simultaneous
  optimization (also known as co-optimization) of the trajectory and the
  spacecraft design parameters \citep{arya2021composite}}.

In a similar way to the previous sections, improvements in convex optimization
technology have prompted an increased interest in applying the direct family of
methods to OTI. {For example, \citep{betts2003optimal} solved a minimum-fuel
  Earth to Moon transfer using a solar electric propulsion system, which is a
  complex problem with a transfer duration of over 200 days. The problem is
  highly nonconvex, and the optimization algorithm is based on SQP.} This
section discusses some of the recent developments for solving OTI using convex
optimization-based methods, and their extensions to optimal exo-atmospheric
launch vehicle ascent.

A basic optimal OTI problem is illustrated in \figref{orbit} and can be
formulated as follows:
\begin{optimus}[
  task={\min},
  variables={t_f,\bm{T}},
  objective={\int_0^{t_f} \norm[2]{\bm{T}(t)}\dt},
  label={orbit}
  ]
  \ddot{\bm{r}}(t)=-\mu\norm[2]{\bm{r}(t)}\inv[3]\bm{r}(t)+m(t)\inv\bm{T}(t), \#
  \dot{m}(t)=-\alpha\norm[2]{\bm{T}(t)}, \#
  \norm[2]{\bm{T}(t)}\leq \rho,\#
  m(0)= m_0,~\bm{r}(0)=\bm{r}_0, ~ \dot{\bm{r}}(0)=\dot{\bm{r}}_0, \#
  \bm{\psi}(\bm{r}(t_f), \dot{\bm{r}}(t_f))=0.
\end{optimus}

{Just like in \pref{goddard}, the vector function $\bm{\psi}$ in \eqref{eq:orbit_f}
  describes the final orbit insertion constraints}, usually in the form of
orbital elements. Note that \eqref{eq:orbit_b}-\eqref{eq:orbit_d} are identical to
\pref{rendezvous}. Naturally, we may hope that previously developed lossless
convexification and SCP techniques from \citep{lu2013autonomous} apply for
OTI. The main novelty is the nonlinear insertion constraint
\eqref{eq:orbit_f}. \citep{liu2014solving} showed that \eqref{eq:orbit_f} can be
linearized with a second-order correction term, and \pref{orbit} can be solved via
SCP as a sequence of SOCPs. The method is efficient and reliable, even for
extremely sensitive cases like McCue's orbit transfer problem
\citep{mccue1967quasilinearization}. Using similar convexification techniques, a
3D minimum-fuel OTI problem was considered in
\citep{wang2018minimum}. Similarly, a 2D minimum-time OTI problem was studied in
\citep{wang2018optimization}, where the dynamics were parametrized by transfer
angle (i.e., orbit true anomaly) instead of time as the independent
variable. Both works consider circular orbits, where \eqref{eq:orbit_f} can be
linearized using spherical or polar coordinates. \citep{tang2018fuel} solved a
minimum-fuel orbit transfer problem by combining SCP with lossless
convexification and pseudospectral discretization. \citep{song2019solar} studied
a minimum-time interplanetary solar sail mission, where the thrust is replaced
by solar radiation force, and optimized the trajectory via SCP as a sequence of
SOCP problems.

The above paragraph mentions works that deal mainly with orbit transfer. A
companion problem is that of orbit injection, where the vehicle is taken from a
non-orbiting state to a target orbit. This occurs, for instance, in the last
stage of rocket ascent. \citep{liu2014solving} showed that \pref{orbit} can also
model the optimal exo-atmospheric ascent flight of a medium-lift launch
vehicle. In this case, the initial condition \eqref{eq:orbit_e} typically denotes
burnout of the launch vehicle's previous stage. In \citep{liu2014solving},
constraint \eqref{eq:orbit_f} denotes the radius and velocity at the perigee of the
target circular orbit. \citep{li2019optimal} considered a similar optimal ascent
problem where the thrust magnitude is constant, and constraint \eqref{eq:orbit_f}
describes the orbital elements of a general elliptical orbit. Using
pseudospectral discretization and SCP, this optimal ascent problem is solved as
a sequence of SOCPs. \citep{li2020online} further considered optimal ascent
flight in the case of a power system fault. In this case, depending on the
severity of the fault, \eqref{eq:orbit_f} describes progressively relaxed insertion
constraints. constraints. Once the spacecraft is in orbit, an adjacent task is
to avoid debris crossing its path. \citep{Armellin2021collision} develops a
real-time collision avoidance algorithm based on lossless convexification and
SCP, and provides a detailed statistical analysis corroborating the method's
effectiveness. % In addition to SCP and
% pseudospectral discretization, this work used slack variables to relax
% constraints \eqref{eq:orbit_b} and \eqref{eq:orbit_d}.

Mission planning often sits one layer above the OTI problem. For example, a
mission plan may consist of a series of planetary flyby and gravity assist
maneuvers. A mission, then, can be viewed as a sequence of OTI solutions that
minimizes an overall objective such as fuel usage or travel time. A modern
approach to mission planning is through hybrid optimal control, and some
methods were already mentioned at the end of \sssref{ascent}
\citep{ross2005hybrid,stevens2005preliminary,stevens2004earth}. Evolutionary
optimization using genetic algorithms offers an alternative solution for
mission planning \citep{ConwayBook}. This approach was used to plan several
complex missions: a Galileo-type mission from Earth to Jupiter, a Cassini-type
mission from Earth to Saturn, and an OSIRIS-REx type mission from Earth to the
asteroid Bennu \citep{englander2012automated}. The Saturn mission is almost
identical to that used by the actual NASA/ESA Cassini mission, but is obtained
fully automatically at a fraction of the time and cost. The algorithm, known as
the evolutionary mission trajectory generator (EMTG), has been made available
by NASA Goddard as an open-source software package \citep{englanderEMTG}.

%%%%%%%%%%%%%%%%%%%%%%%%%%%%%%%%%%%%%%%%%%%%%%%%%%%%
\section{Outlook}
\label{section:outlook}
%%%%%%%%%%%%%%%%%%%%%%%%%%%%%%%%%%%%%%%%%%%%%%%%%%%%

This paper surveyed promising convex optimization-based techniques for next
generation space vehicle control systems. We touched on planetary rocket
landing, small body landing, spacecraft rendezvous, attitude reorientation,
orbit transfer, and endo-atmospheric flight including ascent and reentry. The
discussion topics were chosen with a particular sensitivity towards
computational efficiency and guaranteed functionality, which are questions of
utmost importance for spaceflight control. We conclude by listing in
\ssref{optimization_software} some of the most popular optimization software now
available to the controls engineer, and outlining in \ssref{outlook_future} future
research directions to which the reader may wish to contribute.

\subsection{Optimization Software}
\label{subsection:optimization_software}

\begin{table*}[!t]
  \centering

  \makebox[\textwidth]{\makebox[1.1\textwidth]{
      \includegraphics[width=1.1\textwidth]{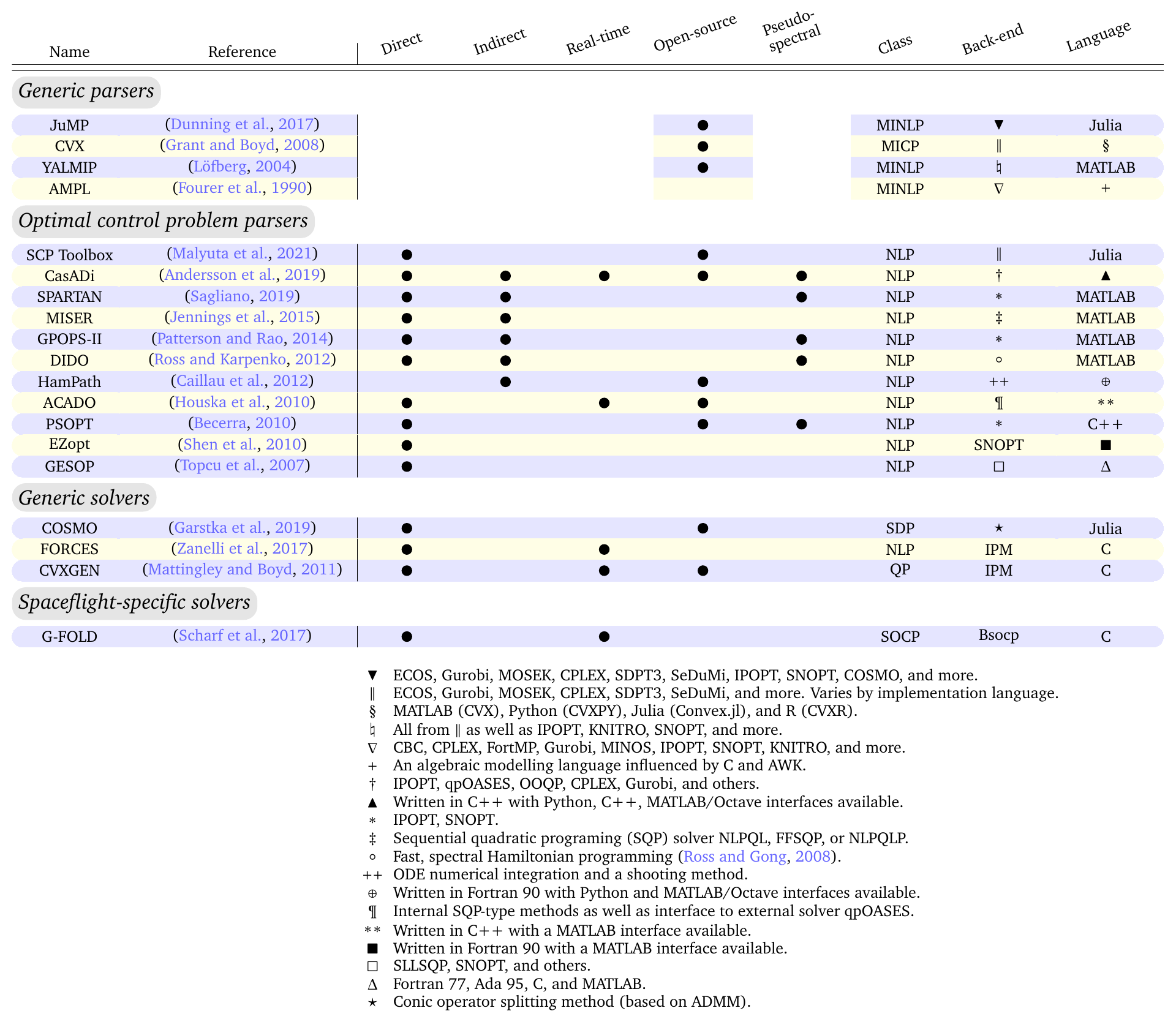}%
    }}

  \caption{{Summary of popular optimization software packages. The columns
      \textit{Direct} and \textit{Indirect} specify which solution method the
      software uses, as discussed in \ssref{ocp}. The column \textit{Real-time}
      denotes if the software is destined for real-time onboard
      use. \textit{Open-source} identifies free-to-download packages with
      viewable source code. \textit{Pseudospectral} identifies software that is
      compatible with pseudospectral discretization. \textit{Class} describes
      the most general class of problems that the software can
      solve. \textit{Back-end} lists which low-level optimizers are used, and
      \textit{Language} lists the implementation and front-end interface
      languages of the package. {Certain classifications that do not apply to
        the ``generic parsers'' software category are indicated by an empty
        cell background.}}}
  \label{table:software}
\end{table*}

\nocite{ross2008guess}
\nocite{scharf2017implementation}
\nocite{Mattingley2011}
\nocite{zanelli2017forces}
\nocite{garstka2019cosmo}
\nocite{topcu2007minimum}
\nocite{shen2010desensitizing}
\nocite{becerra2010solving}
\nocite{houska2010acado}
\nocite{caillau2012differential}
\nocite{ross2012review}
\nocite{patterson2014gpops}
\nocite{jennings2015miser}
\nocite{sagliano2019generalized}
\nocite{andersson2019casadi}
\nocite{SCPTrajOptCSM2021}
\nocite{fourer1990modeling}
\nocite{lofberg2004yalmip}
\nocite{grant2008graph}
\nocite{DunningHuchetteLubin2017}

Success in any computational engineering discipline owes in large part to the
availability of good software. \tabref{software} lists modern optimization
software packages that facilitate the implementation of algorithms discussed in
\sref{applications}. {This list is by no means complete, and should be
  understood to merely indicate some of the popular optimization software
  packages that are quite mature and already available today.}

\subsection{Future Directions}
\label{subsection:outlook_future}

We conclude this survey paper by listing some interesting and important future
directions for optimization-based space vehicle control.

\subsubsection{Guaranteed Performance}
\label{subsubsection:guarantees}

When proposing a new control algorithm for a real system, it is sobering to
remember that the vehicle's survival, along with that of its occupants,
literally hangs in the balance \citep{stein2003respect}. The modern controls
engineer has immense responsibility both to mission success and to upholding the
foundation of trust created by the high reliability of traditional control
methods. If we cannot guarantee an equal or greater level of reliability, then
new optimization-based control methods will quite certainly be relegated
to a ground support role \citep{ploen2006comparison}.

Consequently, a direction of great importance for optimization-based space
vehicle control is to rigorously certify that optimization-based algorithms
converge to solutions that yield safe and robust operation in the real
world. Active research is being done in the area, but general results are
limited and many promising optimization-based methods lack proper
guarantees. Today, researchers are looking at real-time performance
\citep{reynolds2020real,malyuta2020fast}, optimality \citep{reynolds2020optimal},
and convergence rates \citep{mao2018successive,bonalli2019gusto}. Perhaps the
most important yet difficult guarantee is that the algorithm terminates in
finite time, which is imperative for control. In the convex setting, algorithms
with guaranteed convergence are available and have been flight-tested
\citep{dueri2017customized,scharf2017implementation}, so one direction to
explore is how to convexify more general types of nonlinearity
\citep{malyuta2019lossless,liu2014solving,lee2014dual}. For more difficult
nonlinearities that are not convexifiable, an emerging subject of
\textit{funnel libraries} is being {investigated
  \citep{majumdar2017funnel,reynolds2020temporally,ReynoldsThesis,acikmese2008enhancements}}. The
idea, akin to explicit MPC, is to pre-compute a lookup table of trajectories
and invariant controllers in order to replace onboard optimization with a
search algorithm followed by, in some cases, numerical integration. This can
result in a substantially simpler onboard implementation at the expense of a
higher storage memory footprint.

\subsubsection{Machine Learning}
\label{subsubsection:learning}

Impressive advances in machine learning, and particularly in reinforcement
learning (RL), could not side-step space vehicle control without due
consideration \citep{tsiotras2017toward}. The main advantage of RL is that it is
able to optimize over a stochastic data stream rather than assuming a particular
description of a dynamic model
\citep{buoniu2018reinforcement,arulkumaran2017deep}. As an optimization tool for
nonlinear stochastic systems, it is not surprising that the RL method is
attractive for aerospace control.

Although RL for space vehicle control is less than a decade old, a certain
amount of literature is now available that addresses almost all of the
applications presented in \sref{applications}. The reader is referred to
\citep{izzo2019survey} for a dedicated survey. In powered descent guidance,
\citep{cheng2019real} use deep RL (DRL) for lunar landing,
\citep{gaudet2020adaptive} improve ZEM/ZEV guidance via DRL, and
\citep{gaudet2019adaptive} use recursive RL for Mars landing. For spacecraft
rendezvous, \citep{scorsoglio2019actor} use actor-critic RL (ACRL) in
near-rectilinear orbits, \citep{gaudet2018spacecraft} consider cluttered
environments, and \citep{doerr2020space,linares2018physically} use inverse RL to
learn the target's behavior. In reentry guidance, \citep{shi2020deep} aim for
real-time computation by training a deep neural network (DNN) to learn the
functional relationship between state-action pairs obtained from a
high-fidelity optimizer. {Alternatively, \citep{cheng2020multi} use a DNN to
  provide a numerical predictor-corrector guidance algorithm with a range
  prediction based on the current vehicle state. This method improves runtime
  performance by replacing traditional propagation-based trajectory prediction
  with a neural network. A different line of work is presented
  in\citep{jin2017neural}, where} the attitude of a reentry vehicle with model
uncertainty and external disturbances is controlled by a robust adaptive fuzzy
PID-type sliding mode controller designed using a radial basis function neural
network. For small body landing,
\citep{gaudet2020terminal,gaudet2020six,gaudet2019seeker} use RL meta-learning
for greater adaptability, and \citep{cheng2020real} train several DNNs to
approximate a nonlinear gravity field as well as the optimal solution obtained
using an indirect method. {Another interesting approach was proposed in
  \citep{cheng2020fast}, where DNNs are used to supply good costate initial
  guesses, while an accurate trajectory is obtained by a downstream shooting
  method and a homotopy process.} In orbit insertion and transfer applications,
\citep{cheng2019transfer} develop a multiscale DNN architecture to approximate
the optimal solution for a solar sail mission, \citep{holt2020low} use ACRL for
low-thrust trajectory optimization under changing dynamics,
\citep{lafarge2020guidance} use RL for libration point transfer in lunar
applications, and \citep{miller2019low,miller2019interplanetary} use proximal
policy optimization.

A promising modern direction for spacecraft trajectory RL is to learn a small
number of ``behind the scenes'' parameters (called \textit{solution
  hyperparameters}) that govern the optimal solution, instead of directly
learning the high-dimensional optimal state-input map. Most importantly, the
relationship between these parameters and the control policy is much more
predictable, and hence can be learned more easily and with less training
data. This survey paper makes it clear that most if not all spaceflight
trajectory generation problems can be formulated as a variant of the optimal
control \pref{ocp}. Hence, the solution hyperparameters are often the maximum
principle costates, or combinations thereof, that completely define the optimal
control policy. Among these, we find aforementioned concepts of a
\textit{primer vector} \citep{acikmese2007convex,lu2010highly,LawdenBook}, and
switching functions for bang--bang control \citep{taheri2020novel}. This RL
approach was shown to be effective for 3-DoF PDG in
\citep{sixiong2020learning,sixiong2020real}, where the authors learned 10
hyperparameters instead of the map from a 7D state to a 3D input. Most
importantly, only $\approx 10^3$ training trajectories were required. In
comparison, the state-input map learning approach of
\citep{sanchezsanchez2018dnn} also achieved good results, but required
$\approx 10^7$ training samples. A slightly different approach was taken for
3-DoF small body landing in \citep{cheng2020fast}, where homotopy and coordinate
transforms were used to learn a 5D costate vector instead of the map from a 7D
state to a 3D input. The DNN's output was then used to provide accurate initial
guesses and to improve the convergence of a downstream shooting method. To
summarize, the fact that learning hyperparameters works better than learning
the optimal state-input mapping is just an observation of the fact that
application domain knowledge can go a long way towards improving learning
performance \citep{tabuada2020data}. In the case of spacecraft trajectory
optimization, this knowledge often comes from applying Pontryagin's maximum
principle.

As discussed in \sssref{guarantees}, performance guarantees for an RL-based
controller will have to be provided before serious onboard consideration. This
may be harder to achieve for RL, since controllers are typically based on
neural networks whose out-of-sample performance is still very difficult to
characterize. Nevertheless, RL and other machine learning approaches are
appealing for adaptive control systems. Future research will likely see the
aerospace control community search for the right opportunities where RL can be
embedded to improve traditional control systems.

%%%%%%%%%%%%%%%%%%%%%%%%%%%%%%%%%%%%%%%%%%%%%%%%%%%%
% Acknowledgements
%%%%%%%%%%%%%%%%%%%%%%%%%%%%%%%%%%%%%%%%%%%%%%%%%%%%

\section*{Acknowledgments}

This work is funded in part by the National Science Foundation (Grants
CMMI-1613235 and ECCS-1931744) and the Office of Naval Research (Grant
N00014-20-1-2288). The authors would like to extend their gratitude to Yuanqi
Mao for his helpful inputs on sequential convex programming algorithms, and to
Taylor P. Reynolds for his invaluable review of an early draft and for his
knowledge of space vehicle control at large.

%%%%%%%%%%%%%%%%%%%%%%%%%%%%%%%%%%%%%%%%%%%%%%%%%%%%
% Bibliography
%%%%%%%%%%%%%%%%%%%%%%%%%%%%%%%%%%%%%%%%%%%%%%%%%%%%

\bibliographystyle{elsarticle-harv}
\bibliography{new_bib.bib}

\end{document}